\documentclass[sn-mathphys-num]{sn-jnl}% Math and Physical Sciences Numbered Reference Style 
\usepackage[resetlabels]{multibib}
\newcites{SM}{Online Supplement References}

\usepackage{silence}
\usepackage{graphicx}%
\usepackage{multirow}%
\usepackage{amsmath,amssymb,amsfonts}%
\usepackage{amsthm}%
\usepackage{mathrsfs}%
\usepackage[title]{appendix}%
\usepackage{textcomp}%
\usepackage{manyfoot}%
\usepackage{chngcntr}

\usepackage{booktabs}%
\usepackage{algorithm}%
\usepackage{algorithmicx}%
\usepackage{algpseudocode}%
\usepackage{listings}%
\usepackage{thmtools,thm-restate}
\usepackage{lscape} % Add this to your preamble

\usepackage{subcaption} % for subfigure environment

\usepackage{braket}
\usepackage{setspace}
\usepackage{blkarray}
\usepackage{changepage}
\usepackage[table,xcdraw]{xcolor}
\usepackage{mathtools}
\usepackage{tabularx}
\usepackage{hyperref}
\usepackage{ragged2e}
\usepackage{enumitem}

\DeclarePairedDelimiterX{\norm}[1]{\lVert}{\rVert}{#1}

\newcommand{\bs}[1]{\boldsymbol{#1}}
\newlength{\littlespace}
\makeatletter
\def\EGT{\footnotesize}% Define a font size for notes
\def\FigureNoteStyle{\EGT\rm} % Style for the figure notes
\def\FigureNoteName{Note.} % Text label for notes
\def\HD#1#2{\vrule height #1pt depth #2pt width 0pt\relax} % Define spacing for layout

\newcommand{\FIGURE}[3]{%
  \begingroup % Start a group to contain temporary definitions
  \long\def\fignote##1{\relax}% Default: do nothing
  \if\relax\detokenize{#3}\relax % Check if the third argument is empty
  \else
    \long\def\fignote##1{\par\noindent
      \FigureNoteStyle\noindent\HD{11}{0}{\it\FigureNoteName}\enskip \justifying ##1\par}%
  \fi
  \begin{minipage}{\textwidth} % Allow the figure to span the full text width
    \caption[#2]{#2}% Second argument: the caption
    \vskip4pt
    \centering % Center the figure
    #1 % First argument: the figure itself
    \vskip4pt
    \fignote{#3}% Third argument: the note (conditionally included)
  \end{minipage}
  \endgroup % End the temporary group
}
\makeatother

\makeatletter
\def\EGT{\footnotesize}% Define a font size for notes
\def\TableNoteStyle{\EGT\rm} % Style for the table notes
\def\TableNoteName{Note.} % Text label for notes
\def\HD#1#2{\vrule height #1pt depth #2pt width 0pt\relax} % Define spacing for layout
\newcommand{\TABLE}[3]{%
  \begingroup % Start a group to contain temporary definitions
  \long\def\tabnote##1{\relax}% Default: do nothing
  \if\relax\detokenize{#3}\relax % Check if the third argument is empty
  \else
    \long\def\tabnote##1{\par\noindent
      \TableNoteStyle\noindent\HD{11}{0}{\it\TableNoteName}\enskip \justifying ##1\par}%
  \fi
  \begin{minipage}{\textwidth} % Allow the table to span the full text width
    \caption[#1]{#1}% First argument: the caption
    \centering % Center the table
    #2 % Second argument: the table content
    \tabnote{#3}% Third argument: the note (conditionally included)
  \end{minipage}
  \endgroup % End the temporary group
}
\makeatother

\newlist{caseenum}{enumerate}{1}
\setlist[caseenum,1]{
    label=\textbf{Case \arabic*:},
    left=0pt,
    itemsep=1ex,
    topsep=1ex
}

\theoremstyle{thmstyleone}%
\newtheorem{theorem}{Theorem}%  meant for continuous numbers
\newtheorem{proposition}[theorem]{Proposition}% 
\newtheorem{lemma}{Lemma}
\newtheorem{corollary}{Corollary}

\theoremstyle{thmstyletwo}%

\theoremstyle{thmstylethree}%
\newtheorem{definition}{Definition}%

\newtheorem{assumption}{Assumption}

\begin{document}

\title[Contextual Scenario Generation for Two-Stage Stochastic Programming]{Contextual Scenario Generation for Two-Stage Stochastic Programming}


\author[1]{\fnm{David} \sur{Islip}}\email{ryan.islip@mail.utoronto.ca \orcid{https://orcid.org/0009-0001-6857-097X
}}\dgr{}

\author*[1]{\fnm{Roy H.} \sur{Kwon}}\email{rkwon@mie.utoronto.ca \orcid{https://orcid.org/0000-0002-0502-1607
}}

\author[2]{\fnm{Sanghyeon} \sur{Bae}}\email{azureharry@kaist.ac.kr \orcid{https://orcid.org/0000-0002-8862-7363
}}

\author[2]{\fnm{Woo Chang} \sur{Kim}}\email{wkim@kaist.ac.kr \orcid{https://orcid.org/0000-0001-8385-9598
}}

\affil[1]{\orgdiv{Department of Mechanical and Industrial Engineering}, \orgname{University of Toronto}, \orgaddress{\street{5 King's College Rd}, \city{Toronto}, \postcode{M5S 3G8}, \state{Ontario}, \country{Canada}}}

\affil[2]{\orgdiv{Department of Industrial and Systems Engineering}, \orgname{Korea Advanced Institute of Science and Technology (KAIST)}, \orgaddress{\street{291 Daehak-ro}, \city{Yuseong-gu, Daejeon}, \postcode{34141}, \country{Republic of Korea}}}


\abstract{Two-stage stochastic programs (2SPs) are widely used for decision-making under uncertainty, but their practical deployment is often limited by the large number of scenarios needed to approximate the conditional distribution of uncertain outcomes. We study contextual scenario generation: given contextual information, learn to produce a small, user-specified set of surrogate scenarios that, when used as input into the 2SP, lead to high-quality 2SP decisions. Existing scenario generation methods either ignore contextual information or are computationally burdensome in this setting. We propose contextual scenario generation (CSG), which learns a mapping from context to a set of surrogate scenarios. We develop two complementary methodologies: (i) a distributional approach that learns a mapping from context to scenarios by minimizing a kernel-based distance to the conditional distribution, and (ii) a task-based approach that selects the mapping to optimize decision quality via differentiating through a learned surrogate of the downstream 2SP objective. Both approaches are broadly applicable and require only repeated solution of the underlying subproblems and 2SPs defined on the generated scenarios. We provide finite-sample generalization guarantees and demonstrate strong empirical performance across multiple 2SP classes.}




\maketitle

\section{Introduction}\label{Introduction}

Two-stage stochastic programming is a widely adopted decision modelling approach that enables decision-makers to account for uncertainty that unfolds over two stages. The structure is as follows: decisions are made in the first stage before some uncertainty is realized, after which the decision-maker has recourse and responds to the realized uncertainty (second stage). There is a cost for the first-stage decision, and for each specific realization of uncertainty, there is a cost linked to the recourse actions. Consequently, a common objective is to minimize the combined cost of the first-stage decision and the expected cost of recourse actions. \citet{Ntaimo2024} provides a recent survey of various applications of two-stage stochastic programming in logistics, portfolio management, and manufacturing, among others.

In the two-stage setting, the decision-maker selects a first-stage decision $\boldsymbol{y}$ by solving a two-stage stochastic program (2SP):
\begin{equation}
    \min_{\boldsymbol{y}} h(\boldsymbol{y}) + \mathbb{Q}(\boldsymbol{y}) \quad \textit{s.t.} \quad \boldsymbol{y} \in \mathcal{Y}, \quad \boldsymbol{y} \in \mathbb{R}^{s_1},  \label{eq:Stage I} \tag{2SP}
\end{equation}
where $h: \text{dom}(h) (\subseteq \mathbb{R}^{s_1}) \rightarrow \text{codom}(h) (\subseteq \mathbb{R})$ models the cost of the first-stage decision, with $\text{dom}(h)$ and $\text{codom}(h)$ denoting the domain and codomain of $h$ respectively. Additionally, $\mathcal{Y}$ is the feasible set for first-stage decisions, $\bs{\omega} \in \Omega \subseteq \mathbb{R}^p$ represents the uncertainty distributed according to the probability measure $\mathbb{P}$, and $\mathbb{Q}(\boldsymbol{y}) = \mathbb{E}_{\mathbb{P}}[Q(\boldsymbol{y}, \bs{\omega})]$ is the expected recourse cost with $Q(\boldsymbol{y}, \bs{\omega})$ denoting the recourse cost of the first-stage decision $\boldsymbol{y}$ when uncertainty $\bs{\omega}$ is realized. After the first-stage decision is made and uncertainty is realized, the decision-maker has recourse $\boldsymbol{z}$ and determines the best recourse action by minimizing the cost. Determining the optimal recourse action corresponds to the following optimization:
\begin{equation}
Q(\boldsymbol{y}, \bs{\omega}) = \min_{\boldsymbol{z}} q(\boldsymbol{z}, \bs{\omega}) \quad  \textit{s.t.} \quad \boldsymbol{z} \in \mathcal{Z}(\boldsymbol{y}, \bs{\omega}), \quad \boldsymbol{z} \in \mathbb{R}^{s_2}, \label{eq:stage2} \tag{Stage II}
\end{equation} 
where $q: \text{dom}(q) (\subseteq \mathbb{R}^{s_2} \times \mathbb{R}^{p}) \rightarrow \text{codom}(q) (\subseteq \mathbb{R})$ models the recourse cost for uncertainty $\bs{\omega} \in \Omega \subseteq \mathbb{R}^p$, the feasible set of recourse actions given the first-stage decision $\boldsymbol{y}$ and uncertainty $\bs{\omega}$ is denoted by $\mathcal{Z}(\boldsymbol{y}, \bs{\omega})$. The following assumption is often employed for 2SPs, and we also adopt it.
\begin{assumption}\label{a1}
    It is assumed that \eqref{eq:stage2} can be solved efficiently.
\end{assumption}

Solving \eqref{eq:Stage I} involves evaluating a multidimensional integral that is usually not analytically tractable or amenable to numerical integration due to high dimensionality. Thus, the distribution is often represented by a finite set of scenarios $\widehat{\Omega} = \{\bs{\omega}^{(j)}\}_{j=1}^M \subseteq \Omega$. Sample average approximation (SAA) is a common approach for generating scenarios, which proceeds by sampling $\widehat{\Omega}$ from $\mathbb{P}$ \citep{kleywegt2002sample}. Using $\widehat{\Omega}$, the following problem is then solved in place of \eqref{eq:Stage I}:
\begin{equation}
\min_{\boldsymbol{y} \in \mathcal{Y} \cap \mathbb{R}^{s_1}} \quad   h(\boldsymbol{y}) + \frac{1}{|\widehat{\Omega}|}\sum_{j = 1}^{|\widehat{\Omega}|}  \min_{\boldsymbol{z} \in \mathcal{Z}(\boldsymbol{y}, \bs{\omega}^{(j)})} q(\boldsymbol{z}, \bs{\omega}^{(j)}) .
\label{eq:2SP-SAA} \tag{2SP-SAA}
\end{equation}
Although SAA is asymptotically optimal, solving \eqref{eq:2SP-SAA} on a small sample of scenarios can lead to suboptimal first-stage decisions. Consequently, large sample sets are required to obtain high-quality solutions. One possible approach to solve \eqref{eq:2SP-SAA} is to solve the \textit{extensive form}; a deterministic program formed by introducing scenario-specific copies of the second-stage variables. The size of the extensive form grows linearly with $M$, necessitating algorithms tailored to specific problem classes. Unfortunately, it is widely accepted that optimally solving \eqref{eq:2SP-SAA} is intractable due to complicating problem components such as non-linear objectives and integrality requirements \citep{patel2022neur2sp, wu2022learning}. Hence, the distribution is often represented by a small set of scenarios consisting of a positive integer $K$ number of scenarios, referred to as surrogate scenarios $\bs{\zeta}_{1 \hdots K} \coloneq (\bs{\zeta}_{1}, \hdots, \bs{\zeta}_{K}) \in \Omega^K$, with the following assumption\footnote{The surrogate scenarios $\bs{\zeta}_{1 \hdots K}$ can be interpreted as a matrix in $\mathbb{R}^{p \times K}$ with $p$ rows and $K$ columns.}.
\begin{assumption}\label{a2}
    An instance of \eqref{eq:2SP-SAA}, defined on $K$ scenarios, can be efficiently solved for sufficiently small values of $K$.
\end{assumption}
Ideally, $\bs{\zeta}_{1 \hdots K}$ are selected to ensure that the solution to \eqref{eq:2SP-SAA} on the reduced set is nearly optimal with respect to the \eqref{eq:Stage I} objective under $\mathbb{P}$. This process is known as \textit{scenario generation}. Typically, scenario generation is done by generating a discretized representation of $\mathbb{P}$ or by sampling a large number of scenarios $\widehat{\Omega}$ from $\mathbb{P}$ and performing \textit{scenario reduction} to reduce its size from $M$ to $K$ with $K << M$. Scenario generation encompasses both approaches. 

In practice, decision-makers are often tasked with making decisions repeatedly in different settings. Often, decision-makers find themselves in a contextual setting where they are faced with (i) a decision-making problem endowed with uncertainty $\bs{\omega}$ that impacts the problem's important elements and (ii) the opportunity to exploit side-information that is correlated with $\bs{\omega}$ at decision-making time. This work focuses on decision-making via 2SPs in a contextual setting. The side information, also referred to as contextual information $\boldsymbol{x} \in \mathcal{X} \subseteq \mathbb{R}^d$, follows a joint distribution with $\bs{\omega}$ denoted by $\mathbb{P}_{(\boldsymbol{x}, \bs{\omega})}$. The marginal distribution of $\boldsymbol{x}$ and the conditional distribution of $\bs{\omega}$ given $\boldsymbol{x}$ are denoted by $\mathbb{P}_{\boldsymbol{x}}$ and $\mathbb{P}_{\bs{\omega}|\boldsymbol{x}}$ respectively. A realization of the context $\hat{\boldsymbol{x}}$ is observed before the first-stage decision, and as such, the decision-maker wishes to solve \eqref{eq:Stage I} with expected costs evaluated over $\mathbb{P}_{\bs{\omega}|\boldsymbol{x} = \hat{\boldsymbol{x}}}$. We consider the following example to motivate the contextual setting. 

  \noindent \textbf{A Motivating Example:}  \citet{higle_sen_1996} introduce CEP1, a two-stage machine capacity expansion problem for a flexible manufacturing facility. The decision-maker must plan the expansion of production capacity for $N$ parts on $n_{\text{machines}}$ machines. In its baseline state, machine $j \in [n_{\text{machines}}]  \coloneq [1,2,\hdots, n_{\text{machines}}]$ is available for $h_j$ hours per week, with additional hours available at a cost of $c_j$ per hour. Machines have usage limits and require maintenance based on usage. Part type $i \in [m]$ can be produced on machine $j$ at a rate of $a_{ij}$ with cost $g_{ij}$ per hour. The demand $\omega_{i}$ for product type $i$ is uncertain. First, the decision-maker determines machine capacities (first stage) and then creates a production plan to minimize costs after demand $\omega_i$ is realized (second stage). Unmet demand is subcontracted at a premium. In the original problem, capacity decisions are made weekly, but if product demand fluctuates more frequently, adjusting capacity more often (e.g., minutely) may be valuable. This situation could arise if the product demands come from a manufacturing line with variable product demand. Here, the context $\bs{x}$ includes data from the assembly line, such as past demand data. Capacity decisions must be made quickly while accounting for uncertain demand, its relationship to the known context, and available recourse.
    
Thus, it is desirable to obtain an approach for efficiently solving 2SPs repeatedly in different contexts $\bs{x}$, such that the resulting solution performs well according to the 2SP objective, evaluated according to $\mathbb{P}_{\bs{\omega} | \bs{x}}$. However, decision-makers typically do not have access to the conditional distribution and instead rely only on historical data, motivating the following assumption.
\begin{assumption}\label{a3}
 The decision-maker does not have access to {\small $\mathbb{P}_{\bs{\omega}|\boldsymbol{x} = \hat{\boldsymbol{x}}}$} and only has access to an independent and identically distributed (iid) sample of $n$ observations from {\small $\mathbb{P}_{(\boldsymbol{x}, \bs{\omega})}$}, denoted by {\small $S = \{(\boldsymbol{x}^{(i)}, \bs{\omega}^{(i)})\}_{i=1}^{n} \overset{iid}{\sim} \mathbb{P}_{(\boldsymbol{x}, \bs{\omega})}$}. 
\end{assumption}

\begin{figure}
    \FIGURE
%    {\includegraphics{figure-filename.pdf}} %Callout External Image
%	Include LaTeX figures
{\includegraphics[width=\textwidth]{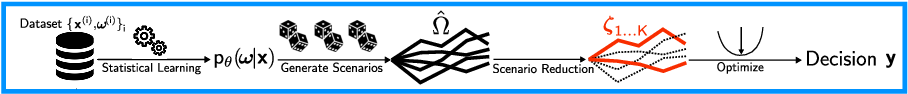}}
{Optimization with Scenario Generation \label{fig:fig1}}
{A standard contextual 2SP approach: (i) estimates the conditional distribution from the data, (ii) generates scenarios via sampling from the estimated distribution, (iii) reduces generated scenarios, and (iv) optimizes based on the reduced set of scenarios }
\end{figure}

It is sensible for the decision-maker to use $S$ to form an estimate of $\mathbb{P}_{\bs{\omega}|\boldsymbol{x}}$, denoted by $p_{\theta}(\bs{\omega} | \boldsymbol{x})$ and parameterized by $\theta$. Figure \ref{fig:fig1} outlines a traditional approach to contextual 2SP based on the abovementioned considerations. First, data is used to estimate $p_{\theta}(\bs{\omega} | \boldsymbol{x})$. A realization of the context $\hat{\boldsymbol{x}}$ is observed. However, solving \eqref{eq:Stage I} with the estimated conditional distribution is not analytically tractable; thus, a common approach is to employ SAA and generate samples via $p_{\theta}(\bs{\omega} | \boldsymbol{x} = \hat{\boldsymbol{x}})$. The resulting instance of \eqref{eq:2SP-SAA} is often unwieldy, and the scenarios are reduced to the set of surrogate scenarios $\bs{\zeta}_{1 \hdots K}$ using a scenario reduction approach. The resulting instance of \eqref{eq:2SP-SAA} is finally solved on the reduced set of scenarios, yielding a first-stage solution. The combination of sampling from the estimated conditional distribution and the subsequent scenario reduction constitutes the standard approach for generating $K$ scenarios in a contextual setting. Figure~\ref{fig:fig1} visually outlines this approach.

A downside of the approach in Figure \ref{fig:fig1} is that when the decision-maker observes a context realization $\hat{\boldsymbol{x}}$, the scenarios must be sampled from $p_{\theta}(\bs{\omega}| \boldsymbol{x})$ and subsequently reduced. Generally, scenario reduction algorithms have computational costs that scale at least linearly in the number of sampled scenarios. The justification of this claim is contained in the overview of scenario generation methods.

\begin{figure}
    \FIGURE
%    {\includegraphics{figure-filename.pdf}} %Callout External Image
%	Include LaTeX figures
{\includegraphics[width=0.9\textwidth]{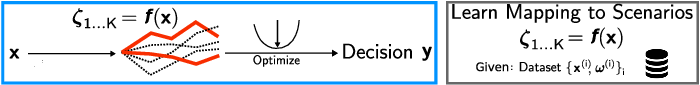}}
{Contextual Scenario Generation and Optimization \label{fig:fig2}}
{Contextual Scenario Generation and Optimization: (i) aims to learn a mapping $\bs{f}$ from data such that $\bs{f}: \bs{x} \mapsto \bs{\zeta}_{1 \hdots K}$, (ii) the predicted scenarios  $\bs{f}$ are close in distributional distance to $\mathbb{P}_{\bs{\omega} | \bs{x}}$, and (iii) solving \eqref{eq:Stage I} on $\bs{f}$ yields high-quality solutions according to the \eqref{eq:Stage I} objective evaluated using $\mathbb{P}_{\bs{\omega}|\bs{x}}$ }
\end{figure}

In time-sensitive applications where solutions to \eqref{eq:Stage I} are desired under different contexts, a mapping $\boldsymbol{f}: \boldsymbol{x} \mapsto \bs{\zeta}_{1 \hdots K}$, that is cheaply evaluated, is desirable so that one obtains first-stage decisions in a short amount of time. This process of mapping context to surrogate scenarios and ultimately to a decision is referred to as \textit{Contextual Scenario Generation (CSG)} and forms the base of the methodology proposed in this work. The mapping $\boldsymbol{f}$ is called the task-mapping. Ultimately, the success of CSG rests on the ability to select $\boldsymbol{f}$ such that solving \eqref{eq:2SP-SAA} on the resulting scenarios yields high-quality decisions. The ultimate goal of this work is to develop a general methodology for constructing a high-quality $\boldsymbol{f}$ for general 2SPs, given that we have access to historical data. The proposed CSG approach is shown visually in Figure \ref{fig:fig2}. 

\subsection*{Contributions}
The main contributions of this work are as follows:

\begin{itemize}
    \setlength\itemsep{-0.0em}
    \item While others consider contextual 2SPs, this work presents the first framework for solving 2SPs in a contextual setting that proceeds by learning a mapping to a set of surrogate scenarios from the context available to the decision-maker. We propose a distributional approach that produces surrogate scenarios that mimic the true conditional distribution, followed by a problem-driven approach that considers the cost of the first-stage solution obtained by solving \eqref{eq:2SP-SAA} on the surrogate scenarios.
     \item The proposed methodology makes minimal assumptions about the structure of the 2SP and only requires Assumptions \ref{a1}--\ref{a3}. 
     \item We provide finite-sample performance bounds that guarantee, with high probability, the out-of-sample performance of the proposed approaches. For the problem-driven approach we focus on the class of linear 2SPs restricted to right-hand-side uncertainty, complete recourse, and dual-feasibility for the recourse problem. 
    \item Computational experiments are performed using four application problems to demonstrate the proposed methodology's ability to produce high-quality solutions in various application domains with differing problem structures. We illustrate the proposed approach using the newsvendor problem, CEP1, portfolio optimization, and a multidimensional newsvendor problem with customer-directed substitutions. 
\end{itemize}

\subsection*{Structure of the paper}
First, Section~\ref{section:Literature} provides background on scenario generation in stochastic programming, learning-based approaches and their use in stochastic programming, and an overview of stochastic optimization in contextual environments. Section~\ref{section:methodology} introduces the proposed distributional and problem-driven approaches. Section~\ref{section:guarantees} presents the performance guarantees associated with the proposed methodologies. Section~\ref{section:experiments} presents the experimental setup and results of the computational experiments. Lastly, Section~\ref{section:conclusion} summarizes the results, contributions, and concludes. All the proofs and omitted results are relegated to an Online Supplement (OS).

\section{Related Works}\label{section:Literature}

\subsection{Scenario Generation}

Scenario generation techniques can be categorized into distributional and problem-based approaches. The distributional approach generates surrogate scenarios to mimic some features of the underlying distribution. This aim is typically achieved by minimizing a distributional distance \citep{dupavcova2003scenario} or by matching moments \citep{hoyland2003heuristic} between the scenarios and the underlying distribution. In contrast, problem-driven scenario generation methods account for problem-specific structures so that the resulting scenarios yield high-quality solutions when evaluated using the underlying distribution. Some problem-driven approaches include: clustering scenarios based on objective values of candidate solutions \citep{bertsimas2023optimization}, identifying irrelevant scenarios based on the objective \citep{fairbrother2022tail}, and approximating the recourse function of a pool of candidate solutions \citep{narum2024problem}. In our setting, the scenario generation methods mentioned assume access to the conditional distribution or a sufficiently accurate scenario representation with enough samples $M$. One can estimate the conditional distribution as in the conditional-density-estimation-then-optimize approach; however, the generic problem-driven approaches depend at least linearly on $M$. In contrast, besides optimizing implementation error, the proposed CSG approach generates surrogate scenarios via a neural network forward pass, eliminating dependence on $M$ at decision time.

\subsection{Learning-based Approaches for Stochastic Programs}

Machine learning has been applied to solve stochastic programs by predicting solutions or costs. In particular, approaches that leverage machine learning to predict the recourse cost are most relevant to this work. For example, \citet{patel2022neur2sp}, \citet{lee2023value} and \citet{bae2023deep} predict expected recourse in generic 2SPs, leveraging these predictions in the solution of 2SPs. However, these methods do not directly address scenario generation. Some works apply machine learning to scenario generation. For instance, \citet{bengio2020learning} use regression to predict a single representative scenario for 2SP such that the scenario minimizes implementation error, but their method relies on heuristically generated datasets and doesn't handle multiple scenarios. \citet{wu2022learning} employ semi-supervised learning and conditional variational auto-encoders to generate scenario embeddings that align with the optimal expected cost of 2SP based on a subset of solved instances. Unlike their approach, which doesn't explicitly minimize implementation error, this work does not require the ability to solve large-scale 2SPs and instead assumes the ability to solve 2SPs on up to $K$ scenarios.

\subsection{Contextual Optimization for Stochastic Programming}
% Learning-based approaches for stochastic programming have experienced success in recent years. Specifically, neural network approaches have been leveraged to solve challenging stochastic programs in both one-off and contextual settings. 

%% Standard contextual approaches conditional-density-estimation-then-optimize

This section highlights works focused on contextual stochastic optimization that are most related to the proposed methodology. We refer the reader to \citet{sadana2024survey} for a more general survey.  Furthermore, we leverage the categorization of contextual stochastic optimization introduced by \citet{estes2023smart}. The first general approach, referred to as \textit{conditional-density-estimation-then-optimize}, estimates $\mathbb{P}_{\bs{\omega}|\bs{x}}$ then solves \eqref{eq:Stage I} with the estimated conditional distribution. For example, in the case of 2SPs, \citet{ban2019dynamic} use residuals from the trained regression models to estimate conditional distributions, whereas \citet{bertsimas2020predictive} reweights samples to approximate $\mathbb{P}_{\bs{\omega}|\bs{x}}$. However, as pointed out in the introduction, this typically yields a difficult problem, motivating the generation of a manageable number of scenarios. 

The second approach referred to as \textit{direct-solution-prediction}, aims to directly estimate a mapping from the context $\bs{x}$ to decisions $\bs{y}$ such that the decisions are of high quality, yielding a policy optimization problem. For example, in the case of 2SPs, \citet{yilmaz2023deep} propose a multi-agent actor-critic approach to solving contextual two-stage knapsack problems where agents predict solutions for the first and second stages, respectively. Although the solution-prediction approach is effective in tailored settings, it is difficult to handle general integrality constraints and other complicating features \citep{patel2022neur2sp,zharmagambetov2023landscape}. 

The \textit{predict-then-optimize} approach generates a point prediction of uncertainty from the context and solves \eqref{eq:Stage I} using a corresponding singleton distribution. Naive versions use standard loss functions (e.g., least-squares), while \textit{smart predict-then-optimize} selects the predictor based on decision performance. For linear programs (LPs), \citet{elmachtoub2021smart} minimize decision regret with a surrogate loss function. \citet{estes2023smart} apply similar methods to linear 2SPs. Other approaches update models using gradients of downstream performance obtained by differentiating optimality conditions and the implicit function theorem \citep{agrawal2019differentiable, donti2017task}. Still, these methods can struggle with uninformative gradients. 
To address this issue, one line of work focuses on generating smoothed proxy gradients. For example, in the case of combinatorial problems: \citet{berthet2020learning} use randomly perturbed optimizers to evaluate gradients, \citet{poganvcic2019differentiation} apply an informed perturbation to the solver input and calling the solver one additional time, and \citet{sahoo2023backpropagationcombinatorialalgorithmsidentity} propose simply backpropagating the incoming gradient from the downstream loss (i.e., they set the solver gradient to negative identity). 
Another approach by \citet{shah2022decision} considers generic problem classes for which access to a solver-oracle is available and constructs a locally convex decision-focused approximate loss for each available observation of uncertainty.

Similar to \citet{shah2022decision}, \citet{zharmagambetov2023landscape} also propose modelling a surrogate loss function, although they define the loss over the domain rather than for each available datapoint. They point out, in comparison to local approaches, that defining the loss over the domain rather than per sample allows for generalization to unseen instances. Furthermore, they point out that a global approximation enables more efficient training as training local models for each sample can be expensive. However, in our work, we observe that selecting a prediction model to optimize a surrogate loss can often yield predictions that lie outside the region where the model accurately captures the true loss, resulting in what we call \textit{loss error maximization}. 

Moreover, these approaches ultimately reduce uncertainty to point predictions, which implicitly assume near-perfect knowledge and can lead to poorly hedged decisions against risk \citep{king2012modeling}. In contrast, \citet{homem2024forecasting} show that for linear 2SPs with fixed recourse and costs, mild conditions guarantee the existence of a single scenario—possibly outside the true support—such that solving the 2SP using only that scenario still yields an optimal solution. Motivated by this result, they use contextual features to predict the optimal scenario. They fit this context-to-scenario map using Nelder-Mead for linear mappings and for piecewise-linear mappings induced by decision-tree partitions. Viewed this way, applying \citet{zharmagambetov2023landscape}'s landscape-surrogate approach to the same problem class provides a natural extension: it targets the same underlying bilevel objective (selecting the scenario prediction model so as to maximize downstream decision quality), but replaces restricted parametric forms with more expressive function classes (e.g., neural networks) that can be trained efficiently via gradient-based optimization. 

Relatedly, within the conditional-density-estimation-then-optimize paradigm, \citet{grigas2021integrated} address the same uninformative gradient issue in a contextual stochastic optimization setting where a learned distribution is fed into the downstream problem. They study convex stochastic programs in which the conditional distribution is represented by context-dependent probability weights over a fixed pre-specified scenario support. Their goal is to learn a mapping from context to these probabilities so that solving the induced stochastic program yields high-quality decisions. Because the optimal solution map is generally non-differentiable with respect to the probabilities, they construct a smooth approximation that permits backpropagation, enabling gradient-based learning of the context-to-probability mapping.

\begin{table}[h]
\TABLE
{Comparison of related approaches along four dimensions. \label{tab:related_work_comparison}}
{
\setlength{\tabcolsep}{4pt}
\renewcommand{\arraystretch}{1.25}
\begin{tabularx}{\textwidth}
{
  >{\raggedright\arraybackslash}p{0.17\textwidth}
  *{4}{>{\raggedright\arraybackslash}X}
}
\toprule
Dimension
& This paper
& \citet{grigas2021integrated}
& \citet{zharmagambetov2023landscape}
& \citet{homem2024forecasting} \\
\midrule
\textbf{Fixed support}
& No; learns $K$ scenarios mapping $\bs{f}$ in $\Omega^K$.
& Yes; learns mapping from $\mathcal{X}$ to probabilities on a fixed scenario support.
& No; learns a single scenario mapping.
& No; learns a single scenario mapping. \\
\textbf{Gradient-based methodology}
& Yes; uses approximation of problem-driven loss to enable gradient-based training.
& Yes; uses approximations of optimal solution map to enable gradient-based learning.
& Yes; uses approximation of problem-driven loss to enable gradient-based training.
& No; uses Nelder-Mead heuristic to solve bilevel formulation for two classes of scenario mappings. \\
\textbf{Guarantees}
& Yes, finite-sample generalization bounds.
& Yes, finite-sample generalization bounds and asymptotic guarantees.
& No performance guarantees provided.
& Yes, one-scenario optimality and asymptotic optimality among scenario mappings. \\
\textbf{Problem class}
& Contextual 2SP; broad (incl.\ mixed-integer) provided 2SP on $K$ scenarios is amenable to repeated solves.
& Convex contextual stochastic optimization provided that the full-problem is amenable to repeated solves.
& Generic; assumed access to a parametric solver that is amenable to repeated solves.
& 2SP with fixed recourse/fixed costs provided the instance defined on a single scenario is amenable to repeated solves.\\
\bottomrule
\end{tabularx}
}
{}
\end{table}

The approaches of \citet{zharmagambetov2023landscape, grigas2021integrated} and \citet{homem2024forecasting} are the most closely related to this work. Table \ref{tab:related_work_comparison} summarizes the key differences along four dimensions: (i) \textbf{fixed support}:  whether the distribution is defined over a fixed, predetermined support, (ii) \textbf{gradient-based methodology}: whether the method relies on gradient-based optimization, enabling flexible context-to-distribution mappings, (iii) \textbf{guarantees}: the type of theoretical guarantees available, and (iv) \textbf{problem class}: the types of stochastic optimization problems to which the approach applies. This work is inspired by \citet{grigas2021integrated}'s, but we do not assume a predetermined fixed support. In addition, we target the same problem-class flexibility as \citet{zharmagambetov2023landscape}'s approach, tailored to 2SPs and effectively develop a gradient-based extension of \citet{homem2024forecasting}'s work to $K$ scenarios.

The proposed CSG approach aims to address the above-mentioned concerns in the general setting of a large class of 2SPs. It circumvents computational concerns by generating a small subset of scenarios and subsequently solving \eqref{eq:2SP-SAA}. This avoids feasibility issues commonly encountered in the direct-solution-prediction approach. Furthermore, the approach inherits all the problem-class generality of \citet{zharmagambetov2023landscape}'s approach while addressing the issue of loss error maximization. The proposed approach makes no assumptions regarding fixing the support of the uncertainty, while still enjoying finite-sample performance guarantees. Lastly, CSG produces solutions that hedge against uncertainty by considering the surrogate scenarios as opposed to a single prediction. We demonstrate both theoretically via performance bounds and empirically in experiments that the additional flexibility of allowing for $K$ scenarios improves decision-making quality.

\section{Proposed CSG Methodology}\label{section:methodology}
This section introduces the relevant tools from scenario generation in Section \ref{subsection:dsg}, followed by the proposed distributional contextual scenario generation approach in Section \ref{subsection:DCSG}. Analogously, Section \ref{subsection:prelim_PSG} introduces a bilevel approach to scenario generation and discusses its properties, followed by its contextual extension to our setting in Section \ref{subsection:PCSG} and the associated solution approach in Section \ref{subsection:solution_approach}.

\subsection{Preliminaries: Distributional Scenario Generation}\label{subsection:dsg}

\citet{romisch2003stability} argues that integral probability metrics (IPMs) are a sensible choice of distances for 2SPs since IPMs such as the Fortet-Mourier and Wasserstein distances lead to stability bounds for 2SPs and form the basis of several existing scenario generation approaches. Given a class of real-valued bounded measurable functions $\mathcal{G}$ associated with the sample space $\Omega$ along with two distributions $\mathbb{P}_{\bs{\omega}}$ and $\mathbb{P}_{\bs{\xi}}$ on $\Omega$, IPMs take the form:
\begin{equation*}\label{eq:ipm}
d_{\mathcal{G}}(\mathbb{P}_{\bs{\omega}}, \mathbb{P}_{\bs{\xi}}) = \sup_{g \in \mathcal{G}} \left|  \mathbb{E}_{\bs{\omega} \sim \mathbb{P}_{\bs{\omega}}} [g(\bs{\omega})] - \mathbb{E}_{\bs{\xi} \sim \mathbb{P}_{\bs{\xi}}} [g(\bs{\xi})] \right|.
\end{equation*}

Different choices of $\mathcal{G}$ yield different distances. For example, when $\mathcal{G}$ is the set of 1-Lipschitz functions $\mathcal{G}_{\text{Lip}}$, $d_{\mathcal{G}_{\text{Lip}}}(\cdot, \cdot)$ corresponds to the 1-Wasserstein distance. Another metric of interest, maximum mean discrepancy (MMD), corresponds to the class $\mathcal{G}_{\text{MMD}} = \{g \in \mathcal{H}: \| g\|_{\mathcal{H}} \leq 1 \} $ where $\mathcal{H}$ denotes a Reproducing Kernel Hilbert Space (RKHS) with associated positive semi-definite kernel $k: \Omega \times \Omega \rightarrow \mathbb{R}$ and norm $\norm{.}_{\mathcal{H}}$. The kernel mean embedding of $\mathbb{P}_{\bs{\omega}}$ in $\mathcal{H}$ is given by 
$\mu_{\mathbb{P}_{\bs{\omega}}} \coloneq \mathbb{E}_{\bs{\omega} \sim \mathbb{P}_{\bs{\omega}}}[k(\cdot,\bs{\omega})]$, and is guaranteed to be an element of $\mathcal{H}$ if $\mathbb{E}_{\bs{\omega} \sim \mathbb{P}_{\bs{\omega}}}[\sqrt{k(\bs{\omega}, \bs{\omega})}] < \infty$\footnote{This technical assumption is not overly restrictive since it holds for continuous kernels on
compact domains or continuous bounded kernels.} \citep{JMLR:v13:gretton12a}. \citet{JMLR:v13:gretton12a} showed that the squared MMD is given by
\begin{align*}
   d^2_{\mathcal{G}_{\text{MMD}}}(\mathbb{P}_{\bs{\omega}}, \mathbb{P}_{\bs{\xi}})& = \norm{\mu_{\mathbb{P}_{\bs{\omega}}} - \mu_{\mathbb{P}_{\bs{\xi}}}}^2_{\mathcal{H}} \\ & =  \mathbb{E}_{(\bs{\omega}, \bs{\omega}') \sim \mathbb{P}_{\bs{\omega}}}[k(\bs{\omega}, \bs{\omega}')] + \mathbb{E}_{(\bs{\xi}, \bs{\xi}') \sim \mathbb{P}_{\bs{\xi}}}[k(\bs{\xi}, \bs{\xi}')] - 2\mathbb{E}_{\bs{\omega} \sim \mathbb{P}_{\bs{\omega}}, \bs{\xi} \sim \mathbb{P}_{\bs{\xi}}}[k(\bs{\omega}, \bs{\xi})].
\end{align*}
A kernel $k$ is characteristic if $\mu : \mathbb{P} \mapsto \mu_{\mathbb{P}}$ is injective. This matters because characteristic kernels guarantee that $d^2_{\mathcal{G}_{\text{MMD}}}(\mathbb{P}_{\bs{\omega}}, \mathbb{P}_{\bs{\xi}}) = 0 \iff \mathbb{P}_{\bs{\omega}} = \mathbb{P}_{\bs{\xi}}$. For more details regarding RKHS and embeddings of probability distributions, see \citep{JMLR:v13:gretton12a, muandet2017kernel}.

For a given context $\bs{x}$, the distributional scenario generation (DSG) problem solves
\begin{equation}\label{eq:DSG}
\min_{\bs{\zeta}_1,\hdots, \bs{\zeta}_K} d_{\mathcal{G}}(\mathbb{P}_{\bs{\zeta}_{1 \hdots K}}, \mathbb{P}_{\bs{\omega} | \bs{x}}), \tag{DSG}
\end{equation}
where $\mathbb{P}_{\bs{\zeta}_{1 \hdots K}} = \frac{1}{K} \sum_{k=1}^K \delta_{\bs{\zeta}_k}$ is the empirical measure associated with the scenarios $\bs{\zeta}_{1 \hdots K}$ and $\delta_{\bs{\zeta}}$ is the Dirac delta function centered at $\bs{\zeta}$. However, as pointed out in the introduction, in practice one only has historical data $S$ that is used to form an estimate $p_{\theta}(\bs{\omega} | \bs{x})$ of $\mathbb{P}_{\bs{\omega} | \bs{x}}$. Subsequently, a sufficiently large number of samples are sampled from $p_{\theta}(\bs{\omega} | \bs{x})$ and \eqref{eq:DSG} is solved with $\mathbb{P}_{\bs{\omega} | \bs{x}}$ replaced with the empirical measure supported on the samples. The next section presents the proposed contextual extension of \eqref{eq:DSG}.

\subsection{Distributional Contextual Scenario Generation}\label{subsection:DCSG}

Given a contextual realization $\hat{\bs{x}}$, the decision-maker wishes to use $\mathbb{P}_{\bs{\omega} | \bs{x} = \hat{\bs{x}}}$ to evaluate the expectation in \eqref{eq:Stage I}. We propose selecting $\bs{f}: \bs{x} \mapsto \bs{\zeta}_{1 \hdots K}$ from a vector-valued function class $\mathcal{F}$ such that the distributional distance between the empirical distribution supported on $\bs{f}(\bs{x})$ and $\mathbb{P}_{\bs{\omega}|\bs{x}}$ is minimized in expectation over $\mathbb{P}_{\bs{x}}$. We refer to $\mathcal{F}$ as the hypothesis set for the task-mapping. This approach is referred to as distributional contextual scenario generation (DCSG):
\begin{equation}\label{eq:DCSG}
\min_{\bs{f} \in \mathcal{F}} \mathcal{L}_{\text{dist}}(\bs{f}) := \mathbb{E}_{\bs{x} \sim \mathbb{P}_{\bs{x}}} \left[ d \left( \nu_{\bs{f}(\bs{x})}, \mathbb{P}_{\bs{\omega} |\bs{x}} \right) \right], \tag{DCSG}
\end{equation}
where $\nu_{\bs{f}(\bs{x})} = \frac{1}{K} \sum_{k=1}^K \delta_{\bs{f}_k(\bs{x})}$ is the empirical measure associated with the scenarios $\bs{f}(\bs{x})$ and $d(\cdot,\cdot )$ is a measure of distance between the distributions. It is not clear what class of distances should be used in the contextual setting. Several authors have considered this question in the context of generative modelling using both Wasserstein and MMD distances. A more thorough comparison of the different approaches to comparing conditional distributions is provided in Section~\ref{appendix:comparing}. \citet{huang2022evaluating}'s work is similar to ours as they consider the same \eqref{eq:DCSG} objective but in the setting of generative models. In this work, we set $d = d^2_{\mathcal{G}_{\text{MMD}}}$ due to the desirable properties it affords $\mathcal{L}_{\text{dist}}$ in the proposed contextual setting, which we discuss next. 

First, \citet{huang2022evaluating} provide the following theorem, showing that $\mathcal{L}_{\text{dist}}(\bs{f})$ is a metric. 
\begin{theorem}[Theorem 4, from \citep{huang2022evaluating}]\label{prop:metric}
 If (i) $k(\cdot,\cdot )$ is characteristic and measurable, (ii) $\mathbb{E}_{\bs{\omega} \sim \mathbb{P}_{\bs{\omega}|\bs{x}}}[k(\bs{\omega}, \bs{\omega})] < \infty$ and (iii) $\mathbb{E}_{\bs{\omega} \sim \nu_{\bs{f}(\bs{x})}}[k(\bs{\omega}, \bs{\omega})] < \infty$ for all $\bs{x} \in \mathcal{X}$, then $\mathcal{L}_{\text{dist}}(\bs{f}) \geq 0$ and $\mathcal{L}_{\text{dist}}(\bs{f}) = 0 \iff\nu_{\bs{f}(\bs{x})} = \mathbb{P}_{\bs{\omega}|\bs{x}}$ almost everywhere according to $\mathbb{P}_{\bs{x}}$. 
\end{theorem}

Next, we note that $\mathcal{L}_{\text{dist}}(\bs{f}) $ can be optimized over $\bs{f}$ using only joint samples. $\mathcal{L}_{\text{dist}}(\bs{f}) $ can be written as
\begin{equation*}
  \begin{aligned}
  \mathcal{L}_{\text{dist}}(\bs{f}) =\ & \mathbb{E}_{\bs{x} \sim \mathbb{P}_{\bs{x}}} {\Bigg [} \mathbb{E}_{(\bs{\omega}, \bs{\omega}') \sim \mathbb{P}_{\bs{\omega}|\bs{x}}}{\Big[}k(\bs{\omega}, \bs{\omega}'){\Big]} + \frac{1}{K^2} \sum_{i=1}^K \sum_{i' = 1}^K k{\big (} \bs{f}_i(\bs{x}), \bs{f}_{i'}(\bs{x}) {\big )}  \\ 
   \quad & -\frac{2}{K} \sum_{i=1}^K \mathbb{E}_{\bs{\omega} \sim \mathbb{P}_{\bs{\omega}|\bs{x}}}{\Big[} k {\big (} \bs{\omega}, \bs{f}_i(\bs{x}) {\big )}{\Big]} 
  {\Bigg]}  \\
  =\ & C + \frac{1}{K^2} \sum_{i=1}^K \sum_{i' = 1}^K \mathbb{E}_{\bs{x} \sim \mathbb{P}_{\bs{x}}} {\Big[} k {\big (} \bs{f}_i(\bs{x}), \bs{f}_{i'}(\bs{x}) {\big )} {\Big]} \\
  \quad & -\frac{2}{K} \sum_{i=1}^K \mathbb{E}_{\bs{x} \sim \mathbb{P}_{\bs{x}}} \mathbb{E}_{\bs{\omega} \sim \mathbb{P}_{\bs{\omega}|\bs{x}}} {\Big[} k {\big (} \bs{\omega}, \bs{f}_i(\bs{x}) {\big )}{\Big]} \\
   = \ & C + \mathbb{E}_{(\bs{x}, \bs{\omega}) \sim \mathbb{P}_{\bs{x},\bs{\omega}}} {\Bigg [} \underbrace{\frac{1}{K^2} \sum_{i=1}^K \sum_{i' = 1}^K k{\big (} \bs{f}_i(\bs{x}), \bs{f}_{i'}(\bs{x}) {\big )} - \frac{2}{K} \sum_{i=1}^K  k {\big (} \bs{\omega}, \bs{f}_i(\bs{x}) {\big )}}_{\coloneq \ell_{\text{MMD}}(\bs{f}(\bs{x}), \bs{\omega})} {\Bigg ],}
  \end{aligned}
  \end{equation*}  
where $C = \mathbb{E}_{\bs{x} \sim \mathbb{P}_{\bs{x}}} \mathbb{E}_{(\bs{\omega}, \bs{\omega}') \sim \mathbb{P}_{\bs{\omega}|\bs{x}}}{\big[}k(\bs{\omega}, \bs{\omega}'){\big]}$ is a constant with respect to $\bs{f}$. \citet{huang2022evaluating} point out that evaluating $C$ requires $\bs{\omega}$ and $\bs{\omega}'$ to be drawn in a conditionally independent manner for each $\bs{x}$, which is not equivalent to globally sampling $(\bs{x}, \bs{\omega}, \bs{\omega}')$, since the latter is not necessarily conditionally independent. Let $\mathcal{L}_{\text{MMD}}(\bs{f}) \coloneq \mathbb{E}_{(\bs{x}, \bs{\omega})} [\ell_{\text{MMD}}(\bs{f}(\bs{x}), \bs{\omega})]$ denote the part of $\mathcal{L}_{\text{dist}}(\bs{f})$ that depends on $\bs{f}$. 

Since $C$ does not depend on $\bs{f}$, optimizing $\mathcal{L}_{\text{MMD}}(\cdot)$ is equivalent to optimizing $\mathcal{L}_{\text{dist}}(\cdot)$ over $\bs{f}$. Let $\hat{\mathcal{L}}_{\text{MMD}}(\bs{f}) \coloneq \frac{1}{n} \sum_{i=1}^n \ell_{\text{MMD}}(\bs{f}(\bs{x}^{(i)}), \bs{\omega}^{(i)})$ denote the sample estimate. In practice, one selects $\bs{f}$ via empirical loss minimization
$
\min_{\bs{f} \in \mathcal{F}} \ \hat{\mathcal{L}}_{\text{MMD}}(\bs{f})
$.
Although $\hat{\mathcal{L}}_{\text{MMD}}(\cdot)$ is non-convex in general, optimization can still be performed via gradient-based methods. In the naive implementation, evaluating the batch gradient over $S$ requires $K^2 + nK$ kernel evaluations, making optimization computationally tractable for small $K$.

In this work, we consider kernel $k_E(\bs{\omega}, \bs{\omega}') = \frac{1}{2} \left(\norm{\bs{\omega}}_2 + \norm{\bs{\omega}'}_2 - \norm{\bs{\omega} - \bs{\omega}'}_2 \right)$ corresponding to the energy-distance \citep{sejdinovic2013equivalence}. \citet{sejourne2023unbalanced} point out that optimization over MMD distances induced by parameterized kernels (e.g., the Gaussian kernel) is sensitive to their parameters, making the parameter-free $k_E$ an attractive choice. Furthermore, one can observe that $d^2_{\text{MMD}}$ with kernel $k_E$ is equivalent to selecting the negative Euclidean kernel $\tilde{k}_E(\bs{\omega}, \bs{\omega}') = -\norm{\bs{\omega} - \bs{\omega}'}_2$. In this case, $d^2_{\text{MMD}}$ respects the underlying geometry via scale equivariance \citep{szekely2012uniqueness}. That is, for $d^2_{\mathcal{G}_{\text{MMD}}}(\mathbb{P}_{\bs{\omega}}, \mathbb{P}_{\bs{\eta}})$, scaling the sample spaces of $\mathbb{P}_{\bs{\omega}}$ and $\mathbb{P}_{\bs{\eta}}$ by $c \in \mathbb{R}$, scales $d^2_{\mathcal{G}_{\text{MMD}}}(\mathbb{P}_{\bs{\omega}}, \mathbb{P}_{\bs{\eta}})$ by $|c|$. Furthermore, it is interesting to note that, in the case of the energy kernel $\tilde{k}_E$, along with setting $K=1$, optimizing $\mathcal{L}_{\text{MMD}}(\bs{f})$ reduces to least-squares. 

The energy distance satisfies the assumptions in Theorem \ref{prop:metric}. Assumptions (ii) and (iii) in Theorem~\ref{prop:metric} reduce to the bounded moment conditions present in Proposition 1 by \citet{szekely2012uniqueness}. Proposition 1 by \citet{szekely2012uniqueness} implies that the kernel $\tilde{k}_E$ is characteristic. Therefore, under the bounded moment conditions by \citet{szekely2012uniqueness}, assumption (i) in Theorem \ref{prop:metric} is automatically satisfied, implying Theorem \ref{prop:metric} holds.

% In this work, we consider kernel $k_E(\bs{\omega}, \bs{\omega}') = \frac{1}{2} \left(\norm{\bs{\omega}}_2 + \norm{\bs{\omega}'}_2 - \norm{\bs{\omega} - \bs{\omega}'} \right)$ corresponding to the energy-distance \citep{sejdinovic2013equivalence}. \citet{sejourne2023unbalanced} point out that optimization over MMD distances induced by parameterized kernels (e.g., the Gaussian kernel) are sensitive to their parameters, making the parameter-free $k_E$ an attractive choice. Furthermore, one can observe that $d^2_{\text{MMD}}$ with kernel $k_E$ is equivalent to selecting the negative Euclidean kernel $\tilde{k}_E(\bs{\omega}, \bs{\omega}') = -\norm{\bs{\omega} - \bs{\omega}'}_2$. In this case, $d^2_{\text{MMD}}$ respects the underlying geometry via scale equivariance \citep{szekely2012uniqueness}. That is, for $d^2_{\mathcal{G}_{\text{MMD}}}(\mathbb{P}_{\bs{\omega}}, \mathbb{P}_{\bs{\eta}})$, scaling the sample spaces of $\mathbb{P}_{\bs{\omega}}$ and $\mathbb{P}_{\bs{\eta}}$ by $c \in \mathbb{R}$, scales $d^2_{\mathcal{G}_{\text{MMD}}}(\mathbb{P}_{\bs{\omega}}, \mathbb{P}_{\bs{\eta}})$ by $|c|$. Furthermore, it is interesting to note that, in the case of the energy kernel $\tilde{k}_E$, along with setting $K=1$, optimizing $\mathcal{L}_{\text{MMD}}(\bs{f})$ reduces to least-squares. 

\subsection{Preliminaries: Problem-driven Scenario Generation}\label{subsection:prelim_PSG}
This section discusses problem-driven scenario generation. We highlight two important features that guide our proposed contextual approach. Firstly, we consider computational issues regarding problem-driven scenario generation, followed by concerns caused by non-unique solutions to \eqref{eq:2SP-SAA} defined on the surrogate scenarios. 

\subsubsection*{Bi-level Problem-driven Scenario Generation}
One can attempt to naively employ a bilevel approach to scenario generation based on selecting the scenarios so that the solution of \eqref{eq:2SP-SAA}, based on said scenarios, minimizes the \eqref{eq:Stage I} objective. For a given context $\bs{x}$, this results in the following bilevel problem:
\begin{align}
\min_{\bs{\zeta}_1, ... ,\bs{\zeta}_K} \quad & L_{\text{task}}(\bs{\zeta}_{1 \hdots K}) \coloneq  h( \bs{y}(\bs{\zeta}_{1 \hdots K}))+ \mathbb{E}_{\bs{\omega} \sim \mathbb{P}_{\bs{\omega}| \bs{x}}}\left[q \left(\bs{z}(\bs{y}(\bs{\zeta}_{1 \hdots K}), \bs{\omega}), \bs{\omega} \right) \right] \quad \tag{BI-SAA} \label{eq:bisaa} \\ \textit{s.t.}  \quad & \bs{z}(\bs{y}(\bs{\zeta}_{1 \hdots K}), \bs{\omega}) \in \underset{z \in \mathcal{Z}(\bs{y}(\bs{\zeta}_{1 \hdots K}), \bs{\omega})}{\text{argmin}} q(\bs{z}, \bs{\omega}) \quad \forall\  \bs{\omega} \in \Omega &   \tag{SP} \\
\quad & \bs{y}(\bs{\zeta}_{1 \hdots K}) \in \underset{\bs{y}}{\text{proj}}\  \underset{\bs{y}, \bs{z}_1, ... \bs{z}_K}{\text{argmin}}\  h( \bs{y}) + \frac{1}{K}\sum_{i = 1}^K q(\bs{z}_i, \bs{\zeta}_i) &  \tag{$\bs{\zeta}$-SAA}\\
\quad & \hspace{11em}\bs{y} \in \mathcal{Y}, \  \bs{z}_i \in \mathcal{Z}(\bs{y}, \bs{\zeta}_i) \ \forall i\  \in \{1,\hdots ,K\}, \notag
\end{align}  
where $\text{proj}_{\bs{y}}$ selects the $\bs{y}$ component of the \hbox{($\bs{\zeta}$-SAA)} optimal solution set. Some notable aspects of \eqref{eq:bisaa} exist. First, there is no need for feasibility constraints in the upper-level problem since $\bs{y}(\bs{\zeta}_{1 \hdots K})$ is feasible by construction. The lower level problem \hbox{($\bs{\zeta}$-SAA)} with solution $\bs{y}(\bs{\zeta}_{1 \hdots K})$ is the surrogate SAA problem defined on $K$ scenarios. \eqref{eq:bisaa} implicitly assumes \hbox{($\bs{\zeta}$-SAA)} admits a unique solution for any $\bs{\zeta}_{1 \hdots K} \in \Omega^K$. Similar to \eqref{eq:DSG}, one can sample a sufficiently large number of scenarios from $p_{\theta}(\bs{\omega} | \bs{x})$ and attempt to solve \eqref{eq:bisaa} with $\mathbb{P}_{\bs{\omega} | \bs{x}}$ replaced with the empirical measure supported on the samples.

\subsubsection*{Challenges with Gradients} \label{section:gradient_problems}
\begin{figure}[H]
     \centering
 \begin{subfigure}[b]{0.49\textwidth}
         \centering
        \includegraphics[width=0.9\textwidth]{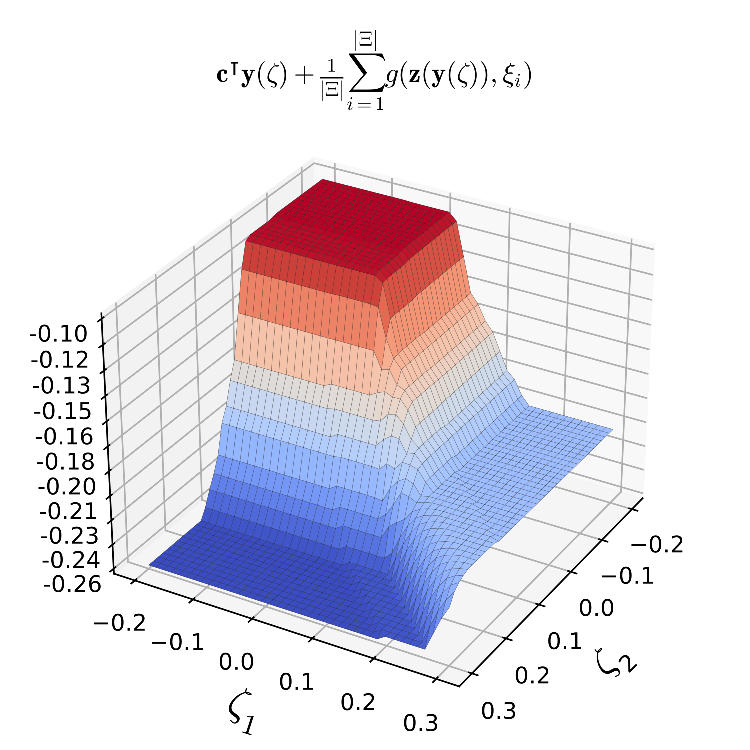}
         \caption{The loss $L_{\text{task}}(\bs{\zeta})$ plotted against two components of the surrogate scenario $\bs{\zeta}$. }
         \label{fig:non_diff}
     \end{subfigure}
     \hfill
   \begin{subfigure}[b]{0.49\textwidth}
         \centering
        \includegraphics[width=0.9\textwidth]{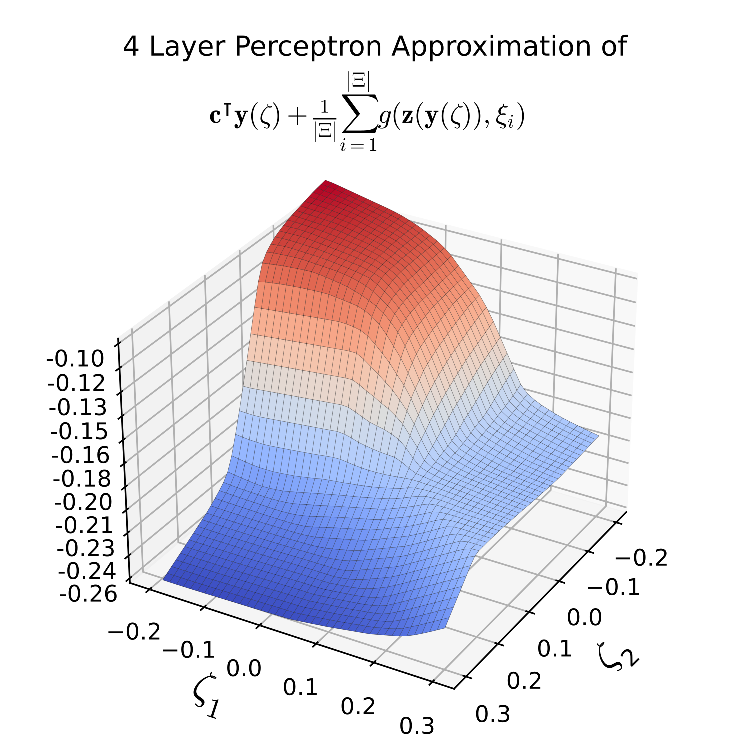}
         \caption{Approximating the loss surface with a neural network has a smoothing effect.}
         \label{fig:diff}
     \end{subfigure}
     \caption{The true loss vs an approximated loss}
     \label{fig:differentiability}
 \end{figure}

Consider the case where $\mathbb{P}_{\bs{\omega}|\bs{x}}$ is replaced with $M$ scenarios. A heuristic approach to solve \eqref{eq:bisaa} is to use gradients of the upper-level objective with respect to the upper-level decision variables $\bs{\zeta}_{1 \hdots K}$. The scenarios can be initialized and updated via gradient descent,
$
\bs{\zeta} \leftarrow \bs{\zeta} - \alpha \ \partial L_{\text{task}} / \partial \bs{\zeta}
$, where the chain rule yields
$$\frac{\partial{L_{\text{task}}}}{\partial{\bs{\zeta}}} = \frac{\partial{L_{\text{task}}}}{\partial{\bs{y}}} \frac{\partial{\bs{y}}}{\partial{\bs{\zeta}}} + \frac{\partial{L_{\text{task}}}}{\partial{\bs{z}}} \frac{\partial{\bs{z}}}{{\partial{\bs{y}}}} \frac{\partial{\bs{y}}}{\partial{\bs{\zeta}}},$$
where $\bs{z} \in \mathbb{R}^{M s_2}$ refers to the $M$ second-stage solutions stacked in a vector. The gradient of the upper-level cost relies on computing the gradients of the surrogate problem solution map with respect to the surrogate scenarios $\frac{\partial{\bs{y}}}{\partial{\bs{\zeta}}}$ and the gradients of the subproblem solution map with respect to the first-stage solution $\frac{\partial{\bs{z}}}{{\partial{\bs{y}}}}$. These gradients or their proxies are typically computed via various methods depending on the problem structure. For example, in the case of convexity, \citet{agrawal2019differentiable}'s \textsf{cvxpy layers} can be used in backpropagation frameworks. Unfortunately, the gradients (if they exist) are commonly sparse and uninformative (even in non-combinatorial problems) \citep{grigas2021integrated,zharmagambetov2023landscape}. Figure \ref{fig:differentiability} demonstrates this behavior for a convex two-stage portfolio selection problem parameterized by a single scenario $\bs{\zeta}$. For instance, Figure \ref{fig:non_diff} plots $L_{\text{task}}(\bs{\zeta})$ against two components of $\zeta$ and shows large flat regions. Methods that rely on differentiating the KKT conditions (e.g. \citep{donti2017task}) for this problem class would likely struggle due to the flat regions. 

Figure \ref{fig:diff} shows a neural network approximation to the true loss surface. It is clear that the approximation has more informative gradients and can approximately optimize over $\bs{\zeta}$ more effectively. This empirical observation motivates using neural architectures to smooth the loss over a particular set of surrogate scenarios and subsequently back-propagate smoothed, problem-driven proxy gradients to the task-mapping $\bs{f}$. \citet{grigas2021integrated} observed similar issues arising from the optimal solution map and proposed approximating it with neural networks. Unlike in \citep{grigas2021integrated}, we do not approximate a high-dimensional solution mapping; instead, we approximate only a scalar-valued loss.

\subsubsection*{Non-uniqueness of Optimal Solutions}\label{section:nonunique_optimal_solutions}

As written, \eqref{eq:bisaa} is not well-defined if \hbox{($\bs{\zeta}$-SAA)} has multiple optimal solutions. Several approaches can address the issue of multiple optimal solutions.  If \hbox{($\bs{\zeta}$-SAA)} is a convex program with a compact feasible set, one can add a small quadratic regularization term $\epsilon \left(\norm{\bs{y}}_2^2 + \sum_{i=1}^K \norm{\bs{z}_i}_2^2 \right)$ to the \hbox{($\bs{\zeta}$-SAA)} objective to induce strict convexity and hence uniqueness. The unique solution obtained by solving the regularized problem is guaranteed to have objective value within $\epsilon D^2$ of the optimal objective, where $D$ is the diameter of the compact feasible set \citep{wilder_melding}. If \hbox{($\bs{\zeta}$-SAA)} has a more complex problem structure (e.g., mixed-integer constraints), then the regularization trick described above may not be sufficient for uniqueness. We resort to the standard conventions in the bilevel optimization literature \citep{sinha2017review}. It is unclear which \hbox{($\bs{\zeta}$-SAA)} solution is implemented in the upper-level objective. In the classic optimistic (pessimistic) setting, nature selects the optimal solution to \hbox{($\bs{\zeta}$-SAA)} that minimizes (maximizes) the upper-level objective. 

We let $\mathcal{Y}^*(\bs{\zeta}_{1 \hdots K})$ denote the set of $\bs{y}$ that form a part of an optimal solution to \hbox{($\bs{\zeta}$-SAA)}, given $\bs{\zeta}_{1 \hdots K}$ i.e. 
\begin{align}
  \mathcal{Y}^*(&\bs{\zeta}_{1 \hdots K}) \coloneq \notag \\
  &\Bigg{\{}\bs{y} \in \mathcal{Y} : h( \bs{y}) + \frac{1}{K}\sum_{i = 1}^K  q(\bs{z}_i, \bs{\zeta}_i) \leq v^*(\bs{\zeta}_{1 \hdots K}), \nonumber 
 \bs{z}_i \in \mathcal{Z}(\bs{y}, \bs{\zeta}_i), \ i \in \{1,\hdots ,K\} \Bigg{\}} , \nonumber
\end{align}
where $v^*(\bs{\zeta}_{1 \hdots K})$ is the optimal objective value of \hbox{($\bs{\zeta}$-SAA)}.
The optimistic version of \eqref{eq:bisaa} can be written as
\begin{equation}\label{eq:optbisaa}
 \min_{\bs{\zeta}_1, ... ,\bs{\zeta}_K}\  \min_{\bs{y} \in \mathcal{Y}^*(\bs{\zeta}_{1 \hdots K})} \quad   h( \bs{y})+ \mathbb{E}_{\bs{\omega} \sim \mathbb{P}_{\bs{\omega}| \bs{x}}}\left[Q \left(\bs{y}(\bs{\zeta}_{1 \hdots K}), \bs{\omega} \right) \right].  \tag{Opt-BI}
\end{equation}

The optimistic setting assumes that any of the best possible solutions among the optimal solution set $\mathcal{Y}^*(\bs{\zeta}_{1 \hdots K})$ according to $L_{\text{task}}(\bs{\zeta}_{1 \hdots K})$ are selected. One could also consider the pessimistic version; given by replacing $\min_{\bs{y} \in \mathcal{Y}^*(\bs{\zeta}_{1 \hdots K})} $ with $\max_{\bs{y} \in \mathcal{Y}^*(\bs{\zeta}_{1 \hdots K})}$ in \eqref{eq:optbisaa}.

\subsection{Problem-driven Contextual Scenario Generation}\label{subsection:PCSG}

We aim to select $\bs{f} \in \mathcal{F}$ such that, given $\bs{x}$, the predicted surrogate scenarios $\bs{f}(\bs{x})$ produce the highest quality set of optimal \hbox{($\bs{\zeta}$-SAA)} solutions $\mathcal{Y}^*(\bs{f}(\bs{x}))$. We measure the quality of the optimal solution set $\mathcal{Y}^*(\bs{f}(\bs{x}))$ by the best possible two-stage performance among solutions in $\mathcal{Y}^*(\bs{f}(\bs{x}))$. The goal of the problem-driven approach is to select $\bs{f}$ such that $\bs{f}$ is in \eqref{eq:optbisaa}'s solution set $\textit{(a.s)}$ with respect to $\mathbb{P}_{\bs{x}}$. This corresponds to 
\begin{equation*}
  \begin{aligned}
\bs{f}(\bs{x}) \in \underset{\bs{\zeta}_1, ... ,\bs{\zeta}_K}{\text{argmin}} \ \min_{\bs{y} \in \mathcal{Y}^*(\bs{\zeta}_{1 \hdots K})} \quad  & h( \bs{y}) + \mathbb{E}_{\bs{\omega} \sim \mathbb{P}_{\bs{\omega} | \bs{x}}} \left[ Q \left(\bs{y}, \bs{\omega} \right) \right] 
\end{aligned}
\end{equation*}
holding (\textit{a.s.}) with respect to $\mathbb{P}_{\bs{x}}$. Theorem 14.60 due to \citet{rockafellar2009variational}, implies this is equivalent to $\bs{f}$ being an optimal solution to the following problem:
\begin{equation} \label{eq:pCSG-opt1}
\min_{\bs{f} \in \mathcal{F}} \quad  \mathbb{E}_{\bs{x} \sim \mathbb{P}_{\bs{x}}} \left[ \min_{\bs{y} \in \mathcal{Y}^*(\bs{f}(\bs{x}))} h(\bs{y}) + \mathbb{E}_{\bs{\omega} \sim \mathbb{P}_{\bs{\omega} | \bs{x}}} \left[ Q \left(\bs{y}, \bs{\omega} \right) \right] \right]. \tag{Opt-PCSG}
\end{equation} 
The following inequality holds:
\begin{align} 
\min_{\bs{y} \in \mathcal{Y}^*(\bs{f}(\bs{x}))} h(\bs{y}) + \mathbb{E}_{\bs{\omega} \sim \mathbb{P}_{\bs{\omega} | \bs{x}}} \left[ Q \left(\bs{y}, \bs{\omega} \right) \right]
= \ & \min_{\bs{y} \in \mathcal{Y}^*(\bs{f}(\bs{x}))} \mathbb{E}_{\bs{\omega} \sim \mathbb{P}_{\bs{\omega} | \bs{x}}} \left[ h(\bs{y}) + Q \left(\bs{y}, \bs{\omega} \right) \right] \nonumber \\ \nonumber
\geq \ &  \mathbb{E}_{\bs{\omega} \sim \mathbb{P}_{\bs{\omega} | \bs{x}}} \left[ \min_{\bs{y} \in \mathcal{Y}^*(\bs{f}(\bs{x}))} h(\bs{y}) + Q \left(\bs{y}, \bs{\omega} \right) \right],
\end{align}
where the inequality follows by relaxing nonanticipativity among the optimal solution set $\mathcal{Y}^*(\bs{f}(\bs{x}))$. Thus, the following optimistic relaxation of \eqref{eq:pCSG-opt1} follows:
\begin{equation}\label{eq:pCSG-opt2}
    \min_{\bs{f} \in \mathcal{F}} \quad   \mathbb{E}_{(\bs{x}, \bs{\omega}) \sim \mathbb{P}_{\bs{x}, \bs{\omega}}} \left[
    \min_{\bs{y} \in \mathcal{Y}^*(\bs{f}(\bs{x}))} h(\bs{y}) + Q \left(\bs{y}, \bs{\omega} \right) \right]. \tag{Opt-PCSG'}
\end{equation} 
The advantage of \eqref{eq:pCSG-opt2} over \eqref{eq:pCSG-opt1} is that it does not require sampling access to $\mathbb{P}_{\bs{\omega}|\bs{x}}$. Consequently, \eqref{eq:pCSG-opt2} suggests the following loss function:
$$\ell_{\text{opt}}(\bs{\zeta}_{1 \hdots K}, \bs{\omega}) \coloneq \min_{\bs{y} \in \mathcal{Y}^*(\bs{\zeta}_{1 \hdots K})} h(\bs{y}) + Q \left(\bs{y}, \bs{\omega} \right).$$ 
Evaluating $\ell_{\text{opt}}(\bs{\zeta}_{1 \hdots K}, \bs{\omega})$ requires solving ($\bs{\zeta}$-SAA), obtaining the optimal value $v^*(\bs{\zeta}_{1 \hdots K})$, then given $\bs{\omega}$, finding the solution among $\mathcal{Y}^*(\bs{f}(\bs{x}))$ that minimizes $h(\bs{y}) + Q \left(\bs{y}, \bs{\omega} \right)$. The ability to efficiently evaluate $\ell_{\text{opt}}$ relies on i) the ability to efficiently solve \hbox{($\bs{\zeta}$-SAA)} on $K$ scenarios (Assumption \ref{a2}), and ii) the ability to efficiently minimize the 2SP objective defined by a single scenario $\bs{\omega}$ over the optimal solution set $\mathcal{Y}^*(\bs{f}(\bs{x}))$. Suppose \hbox{($\bs{\zeta}$-SAA)} obtains optimal value $v^*(\bs{\zeta}_{1 \hdots K})$, then evaluating $\ell_{\text{opt}}(\bs{\zeta}_{1 \hdots K}, \bs{\omega})$ is equivalent to
  \begin{align}\label{eq:opt_search}
  \ell_{\text{opt}}(\bs{\zeta}_{1 \hdots K}, \bs{\omega} ) = \min_{\bs{y}, \bs{z}, \bs{z}_1, ... \bs{z}_K} \quad & h( \bs{y})+q(\bs{z}, \bs{\omega}) \quad  \tag{Opt-Search} \\ \textit{s.t.}  
 \quad & h( \bs{y}) + \frac{1}{K}\sum_{i = 1}^K  q(\bs{z}_i, \bs{\zeta}_i) \leq v^*(\bs{\zeta}_{1 \hdots K}) & \label{eq:optimality} \\
 \quad & \bs{y} \in \mathcal{Y},\ \bs{z} \in \mathcal{Z}(\bs{y}, \bs{\omega}),\  \bs{z}_i \in \mathcal{Z}(\bs{y}, \bs{\zeta}_i), \ \forall i\  \in \{1,\hdots ,K\}.\nonumber
  \end{align} 
\eqref{eq:opt_search} has the same constraints as \hbox{($\bs{\zeta}$-SAA)}, along with constraint \eqref{eq:optimality} that ensures the optimality of $(\bs{y}, \bs{z}_1, ... \bs{z}_K)$ with respect to \hbox{($\bs{\zeta}$-SAA)}. An additional decision variable $\bs{z}$ is introduced to model the recourse in response to scenario $\bs{\omega}$. In addition to LPs and convex programs with multiple solutions, the optimistic approach is generally amenable to mixed-integer programs (MIP) with convex relaxations since, in this case, the resulting instance of \eqref{eq:opt_search} is a convex MIP with one additional constraint \eqref{eq:optimality} and decision variable $\bs{z}$.

In the case where ($\bs{\zeta}$-SAA) exhibits unique solutions, $\ell_{\text{opt}}$ reduces to the simplified loss function $h(\bs{y}(\bs{\zeta}_{1 \hdots K})) + Q(\bs{y}(\bs{\zeta}_{1 \hdots K}), \bs{\omega})$ whose evaluation only requires solving \hbox{($\bs{\zeta}$-SAA)} (Assumption \ref{a2}) and a single subproblem (SP) (Assumption \ref{a1}) to obtain $\bs{y}(\bs{\zeta}_{1 \hdots K})$ and $\bs{z}(\bs{y}(\bs{\zeta}_{1 \hdots K}), \bs{\omega})$, respectively. Furthermore, if solving \eqref{eq:opt_search} proves too computationally burdensome, then one can easily compute $h(\bs{y}(\bs{\zeta}_{1 \hdots K})) + Q(\bs{y}(\bs{\zeta}_{1 \hdots K}), \bs{\omega})$ as a heuristic loss (although learning may not be well-defined in this case). This work uses $\ell_{\text{opt}}$ for all problems under consideration since they lack uniqueness guarantees and are amenable to direct computation by solving \eqref{eq:opt_search}. Similar to the optimistic approach, one can consider a pessimistic variant. Following the same line of reasoning suggests $\max_{\bs{y} \in \mathcal{Y}^*(\bs{\zeta}_{1 \hdots K})} h(\bs{y}) + Q \left(\bs{y}, \bs{\omega} \right),$ as a loss function. However, evaluating this loss is more difficult than evaluating $\ell_{\text{opt}}$. For instance, in the case of a two-stage LP, evaluating the pessimistic loss is equivalent to maximizing a convex function, which is NP-hard \citep{zwart}. As discussed in the next Section, our proposed methodology hinges on repeatedly computing $\ell_{\text{opt}}$ to construct its approximation. Due to this added complexity, this work exclusively focuses on the optimistic case, and we leave exploring the merits of assuming the worst-case among the optimal solution set as future work.

\subsection{Solution Methodology}\label{subsection:solution_approach}
Thus far, two loss functions have been introduced. In the case of  $\ell_{\text{MMD}}$, empirical risk minimization over a parametric class $\bs{f}_{\phi}, \ \phi \in \varPhi$ via gradient-based methods is straightforward. We will typically take $\varPhi$ as a class of neural networks and refer to $\bs{f}$ as a Task-Net. More details are provided in the experimental section. As discussed in Section \ref{section:gradient_problems}, $\ell_{\text{opt}}$ will in general exhibit sparse gradients with respect to the output of $\bs{f}_{\phi}$, making back-propagation of gradients ineffective. To address this, a neural network architecture (Loss-Net) is proposed to approximate the problem-based loss $\ell_{\text{opt}}$.

\begin{figure}
  \FIGURE
   {\includegraphics[width=\textwidth]{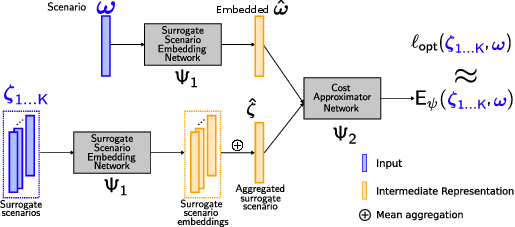}}
    {Loss-Net Architecture
    \label{fig:lossnet}}
    {}
  \end{figure}

The network is a mapping $E_{\psi}$ from $\Omega^{K+1} \rightarrow \mathbb{R}$ with the aim that $E_{\psi}(\bs{\zeta}_{1 \hdots K}, \bs{\omega}) \approx \ell_{\text{opt}}(\bs{\zeta}_{1 \hdots K}, \bs{\omega})$. The network's parameters are denoted by $\psi$, which are assumed to lie in a set $\varPsi$. The architecture's design is motivated by permutation invariant neural architectures \citep{zaheer2017deep,tabaghi2024universal}. In particular, \citet{patel2022neur2sp}'s Neur2SP architecture is notable because it also applies permutation-invariant neural architectures to 2SPs. The network embeds each surrogate scenario $\bs{\zeta}_k,\  k \in [K]$ into a latent space $\Omega' \subseteq \mathbb{R}^{p_{\text{latent}}}$ using an encoder $ \Psi_1 : \Omega \rightarrow \Omega'$. The embedded surrogate set is then represented as a single encoded scenario $\widehat{\bs{\zeta}}$ via mean aggregation i.e., $\widehat{\bs{\zeta}} = \frac{1}{K} \sum_{k=1}^K  \Psi_1(\bs{\zeta}_{k})$. This ensures that the predictions from the network are invariant to the ordering of the input set $\bs{\zeta}_{1 \hdots K}$. The input scenario is also embedded via $\Psi_1$: $\hat{\bs{\omega}} =  \Psi_1(\bs{\omega})$. The embedding of the surrogate scenarios $\widehat{\bs{\zeta}}$ and embedded input scenario $\hat{\bs{\omega}}$ are fed into a separate network $ \Psi_2:\Omega' \times \Omega' \rightarrow \mathbb{R}$ that outputs the approximation of $\ell_{\text{opt}}(\bs{\zeta}_{1 \hdots K}, \bs{\omega})$ given by $$E_{\psi}(\bs{\zeta}_{1 \hdots K}, \bs{\omega}) = \Psi_2\left(\widehat{\bs{\zeta}},\  \hat{\bs{\omega}}\right) = \Psi_2\left(\frac{1}{K} \sum_{k = 1}^K  \Psi_1(\bs{\zeta}_k),\  \Psi_1(\bs{\omega}) \right).$$ The networks, $\Psi_1$ and $\Psi_2$, are taken to be fully connected feedforward neural networks with Rectified Linear Unit (ReLU) activations of appropriate input and output dimensions, each with hyperparameters such as the number of hidden layers and activations. The architecture is shown in Figure~\ref{fig:lossnet}.

It would be ideal for $\ell_{\text{opt}}$ to be representable via a structure similar to $E_{\psi}$. Indeed, the following proposition shows that $\ell_{\text{opt}}$ can be represented as a composition of functions like $E_{\psi}$ when the scenario embedding space $\Omega'$ is of large dimensionality.

\begin{restatable}{proposition}{representation}\label{prop:representation}
% \begin{prop}\label{prop:representation}
  For a fixed positive integer $K$, the (continuous) task-based loss $\ell_{\text{opt}} : \Omega^K \times \Omega \rightarrow \mathbb{R}$, where $\Omega \subseteq \mathbb{R}^p$ is compact, satisfies
  \begin{equation*}\label{eqn:representer}
    \ell_{\text{opt}}(\bs{\zeta}_{1 \hdots K}, \bs{\omega}) = \rho\left(\frac{1}{K}\sum_{k=1}^K \bs{\tau}(\bs{\zeta}_k), \bs{\tau}(\bs{\omega})\right) \quad \forall \bs{\zeta}_{1 \hdots K} \in \Omega^K, \bs{\omega} \in \Omega,
    \end{equation*}
  with  continuous $\bs{\tau} : \mathbb{R}^p \rightarrow \mathbb{R}^{{K+p \choose p} - 1}$, (continuous) $\rho : \mathbb{R}^{{K+p \choose p} - 1} \times \mathbb{R}^{{K+p \choose p} - 1} \rightarrow \text{codom}(\rho)$, and $\text{codom}(\ell_{\text{opt}}) \subseteq \text{codom}(\rho)$.
 %\end{prop}
\end{restatable}
The proof of Proposition \ref{prop:representation} is a direct application of Proposition 1 by \citet{tabaghi2024universal}. The proof introduces additional notation and is provided in Section~\ref{appendix:prop2}. Given the representation provided by Proposition~\ref{prop:representation}, approximating $\ell_{\text{opt}}$ by $E_{\psi}$ amounts to replacing $\rho$ and $\bs{\tau}$ with fully connected feedforward ReLU networks $\Psi_2$ and $\Psi_1$, respectively. Proposition~\ref{prop:representation} produces a representation such that the continuity of $\ell_{\text{opt}}$ determines the continuity of $\rho$. At the same time, $\bs{\tau}$ is continuous, independent of the continuity of $\ell_{\text{opt}}$. In either case, by using neural architectures, we leverage a continuous approximation to $\ell_{\text{opt}}$, irrespective of whether $\ell_{\text{opt}}$ is continuous, by simply replacing $\rho$ with $\Psi_2$. Furthermore, the representation in Proposition~\ref{prop:representation} has $\Omega' \subseteq \mathbb{R}^{{K+p \choose p} - 1}$, i.e. leverages a high-dimensional embedding of the scenarios.

\begin{figure}
  \FIGURE
  {\includegraphics[width=\textwidth]{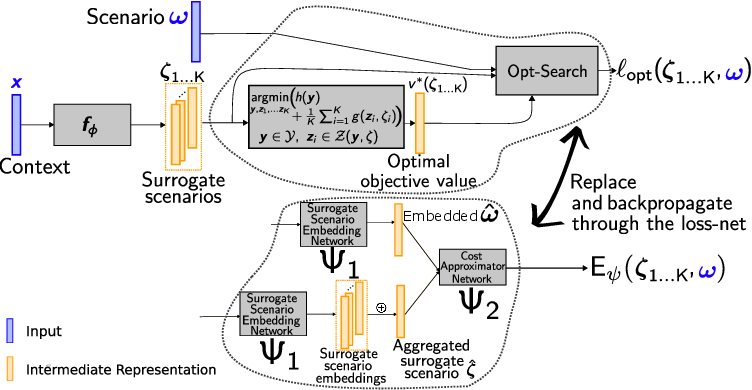}}
  {Replacing $\ell_{\text{opt}}$ with $E_{\psi}$ \label{fig:overall approach}}
  {Replace the true downstream loss with the Loss-Net approximated loss to avoid sparse gradients}
\end{figure}

The loss network acts as a replacement for $\ell_{\text{opt}}(\bs{\zeta}_{1 \hdots K}, \bs{\omega})$. The basic idea underpinning the proposed approach is to replace $\ell_{\text{opt}}(\bs{\zeta}_{1 \hdots K}, \bs{\omega})$ with $E_{\psi}(\bs{\zeta}_{1 \hdots K}, \bs{\omega})$ and use it to learn $\bs{f}_{\phi}$. This idea is visually depicted in Figure~\ref{fig:overall approach}. Given a trained $E_{\psi}$, one can use it to infer $\bs{\zeta}_{1 \hdots K}$ via gradient-based minimization. The remainder of this section describes training procedures to learn $E_{\psi}$ and $\bs{f}_{\phi}$. 

\subsubsection*{Static Approach} 
We do not seek to construct $E_{\psi}$ such that it uniformly approximates $\ell_{\text{opt}}$ over $\Omega^{K+1}$. Instead, we wish to approximate $\ell_{\text{opt}}$ over a distribution $\mathbb{P}_{\bs{\zeta}_{1 \hdots K}, \bs{\omega}}$ that is reflective of the surrogate scenarios $\bs{\zeta}_{1 \hdots K}$ and uncertainty that $E_{\psi}$ will encounter in practice. Once trained, $E_{\psi}$ is used to guide the learning for $\bs{f}_{\phi}$. We refer to this approach as the static approach.

A simple approach is to train $E_{\psi}$ over surrogate scenarios produced by $\bs{\zeta}_{1 \hdots K} = \bs{\hat{f}}_{\text{MMD}}(\bs{x})$, where $\bs{\hat{f}}_{\text{MMD}}$ denotes the Task-Net obtained via minimization of $\hat{\mathcal{L}}_{\text{MMD}}$ over the iid sample $S$. Given $ \bs{\hat{f}}_{\text{MMD}}$, the sample for training $E_{\psi}$ is generated by forming the surrogate scenarios $\bs{\zeta}_{1 \hdots K}^{(i)} = \bs{\hat{f}}_{\text{MMD}}(\bs{x}^{(i)})$, and finally evaluating the task-based loss $\ell_{\text{opt}}^{(i)} = \ell_{\text{opt}}(\bs{\zeta}_{1 \hdots K}^{(i)}, \bs{\omega}^{(i)})$, yielding a dataset \hbox{$S' = \{(\bs{x}^{(i)}, \bs{\zeta}_{1 \hdots K}^{(i)}, \bs{\omega}^{(i)} ), \ell_{\text{opt}}^{(i)} \}_{i=1}^{n}$}. $E_{\psi}$ is then trained by gradient-based methods to minimize the prediction error over $S'$:
\begin{equation*}
\min_{\psi \in \varPsi} \quad \frac{1}{n} \sum_{i = 1}^{n} \left(E_{\psi}(\bs{\zeta}_{1 \hdots K}^{(i)}, \bs{\omega}^{(i)})  - \ell_{\text{opt}}^{(i)} \right)^2,
\end{equation*}
where it is assumed that hyperparameters associated with the learning procedure, such as optimizer learning rates, batch size, and $l_2$ regularization on the weights, can be set such that the resulting Loss-Net achieves generalization. The resulting loss approximation network, obtained from training over $S'$, is denoted by $\hat{E}$.

Although it is tempting to select the Task-Net to minimize the approximate task loss over a sample, we observed that optimizing the approximate task loss tended to yield surrogate scenario predictions $\bs{f}_{\phi}(\bs{x})$ that are far from the distribution of the samples used to train the loss network. I.e., directly minimizing the approximate task loss tends to yield a Task-Net whose predictions maximize the error of the approximate task loss (\textit{loss error maximization}). To mitigate this, we proposed regularizing the approximate task loss by the MMD loss to select $\bs{f}_{\phi}$:
\begin{equation}\label{eq:static_PCSG}
 \min_{\phi \in \varPhi} \quad \frac{1}{n} \sum_{i=1}^{n} \hat{E}(\bs{f}_{\phi}(\bs{x}^{(i)}), \bs{\omega}^{(i)}) + \lambda \ell_{\text{MMD}}(\bs{f}_{\phi}(\bs{x}^{(i)}), \bs{\omega}^{(i)}), \tag{Static-PCSG}
\end{equation}
where $\lambda \geq 0$ is a regularization penalty. The MMD regularization ensures that the surrogate scenario predictions do not stray too far from the input distribution used to train the Loss-Net. Furthermore, the proposed problem is amenable to empirical risk minimization via the sample $S$, making training straightforward. Once $\hat{E}_{\psi}$ is trained, one can input any number of surrogate scenarios into $\Psi_1$, and a numerical result will still be produced. Although we do not explore it here, this can potentially reduce training time by using an $\hat{E}_{\psi}$ that is trained via cheaper evaluations of $\ell_{\text{opt}}$ using $K' < K$ surrogate scenarios. The static approach for selecting $\bs{f}_{\phi}$ is shown in Algorithm \ref{alg:static}. We denote the Task-Net obtained by training on $S'$ by $\bs{\hat{f}}_{\text{PCSG}}$.

\begin{algorithm}[H] 
  \small
  \caption{Training and Using $E_{\psi}$ for Static PCSG}
  \begin{algorithmic}[1]
  \Require Loss function $\ell_{\text{opt}}$, regularization parameter $\lambda$, 
  Network architecture and training parameters (e.g. batch size ($B$), number of epochs ($E$), No. hidden units), Sample $S = \{(\bs{x}^{(i)}, \bs{\omega}^{(i)})\}_{i=1}^n$, Trained MMD network $\bs{\hat{f}}_{\text{MMD}}$
    \State Initialize networks $E_{\psi}$ and $\bs{f}_{\phi}$ and hyperparameters: e.g. batch size ($B$)
  \For{each $(\bs{x}^{(i)}, \bs{\omega}^{(i)}) \in S$}
      \State Evaluate $\bs{\zeta}_{1 \hdots K}^{(i)} = \bs{\hat{f}}_{\text{MMD}}(\bs{x}^{(i)})$,  $\ell_{\text{opt}}^{(i)} = \ell_{\text{opt}}(\bs{\zeta}_{1 \hdots K}^{(i)}, \bs{\omega}^{(i)})$
      % \State Evaluate task-based loss $\ell_{\text{opt}}^{(i)} = \ell_{\text{opt}}(\bs{\zeta}_{1 \hdots K}^{(i)}, \bs{\omega}^{(i)})$
  \EndFor
  \State Form dataset $S' = \{(\bs{\zeta}_{1 \hdots K}^{(i)}, \bs{\omega}^{(i)}, \ell_{\text{opt}}^{(i)})\}_{i=1}^{n}$
  \State Train $\hat{E}$ over $S'$: $\min_{\psi} \  \frac{1}{|S'|} \sum_{\bs{\zeta}_{1 \hdots K}, \bs{\omega}, \ell_{\text{opt}} \in S'} (E_{\psi}(\bs{\zeta}_{1 \hdots K}, \bs{\omega}) - \ell_{\text{opt}} )^2$, yielding $\hat{E}$
  \State Solve \eqref{eq:static_PCSG} over $S$ using $\hat{E}$
  \State \textbf{Output:} Trained Task-Net $\bs{\hat{f}}_{\text{PCSG}}$ and Loss-Net $\hat{E}$ 
  \end{algorithmic}\label{alg:static}
  \end{algorithm}

\subsubsection*{Dynamic Approach}

After running Algorithm \ref{alg:static}, there is a trained Task-Net $\bs{\hat{f}}_{\text{PCSG}}$ and Loss-Net $\hat{E}$. By construction, $\hat{E}$ approximates $\ell_{\text{opt}}$ over the distribution of inputs $(\bs{\hat{f}}_{\text{MMD}}(\bs{x}), \bs{\omega})$, where $(\bs{x}, \bs{\omega}) \sim \mathbb{P}_{\bs{x}, \bs{\omega}}$. However, there is no guarantee that $\hat{E}$ approximates $\ell_{\text{opt}}$ over the analogous distribution induced by $\bs{\hat{f}}_{\text{PCSG}}$. Thus, additional samples for the Loss-Net are generated using $\bs{\hat{f}}_{\text{PCSG}}$, by forming the surrogate scenarios $\bs{\zeta}_{1 \hdots K}^{(i)} = \bs{\hat{f}}_{\text{PCSG}}(\bs{x}^{(i)})$, and evaluating the task-based loss $\ell_{\text{opt}}^{(i)} = \ell_{\text{opt}}(\bs{\zeta}_{1 \hdots K}^{(i)}, \bs{\omega}^{(i)})$. The resulting samples, in addition to the sample generated by $\bs{\hat{f}}_{\text{MMD}}$ can then be used to train $\hat{E}$ further. Once $\hat{E}$ is updated, re-solving \eqref{eq:static_PCSG} updates the Task-Net. The repetition of this process for $T$ iterations is referred to as the dynamic approach and is summarized in Algorithm \ref{algo:dynamic-training}.

\begin{algorithm}[H]
  \caption{Dynamic Training Algorithm for $E_{\psi}$ and $\bs{f}_{\phi}$}
  \label{algo:dynamic-training}
  \small 
  \begin{algorithmic}[1]
    \Require Same inputs as Algorithm \ref{alg:static}, iteration limit $T$
    \State Perform steps 1-5 in Algorithm \ref{alg:static}, yielding $S', \hat{E}$ and $\bs{\hat{f}}_{\text{PCSG}}$
    % \If{\texttt{MMD Pre-Train?}}
    %   \State $(\bs{f}_{\phi}, E_{\psi})$ = Algorithm \ref{alg:static} \Comment{Static PCSG}
    % \EndIf
    \For{$t= 1$ ... $T$}
      \For{each $(\bs{x}^{(i)}, \bs{\omega}^{(i)}) \in S$}
      \State Evaluate $\bs{\zeta}_{1 \hdots K}^{(i)} = \bs{\hat{f}}_{\text{PCSG}}(\bs{x}^{(i)})$,  $\ell_{\text{opt}}^{(i)} = \ell_{\text{opt}}(\bs{\zeta}_{1 \hdots K}^{(i)}, \bs{\omega}^{(i)})$ 
      \State $S' \gets S' + \{(\bs{\zeta}_{1 \hdots K}^{(i)}, \bs{\omega}^{(i)}, \ell_{\text{opt}}^{(i)})\}_{i=1}^{n}$ \Comment{Replay buffer: keep past $T' < T$ iterations}
      \EndFor 
  \State Update $\hat{E}$ over $S'$: $\min_{\psi} \  \frac{1}{|S'|} \sum_{\bs{\zeta}_{1 \hdots K}, \bs{\omega}, \ell_{\text{opt}} \in S'} (E_{\psi}(\bs{\zeta}_{1 \hdots K}, \bs{\omega}) - \ell_{\text{opt}} )^2$
    \State Fix $\hat{E}$ and update $\bs{f}_{\phi}$ via \eqref{eq:static_PCSG} over $S$ using $\hat{E}$
    \EndFor
    \State \textbf{Output:} Trained Task-Net $\bs{\hat{f}}_{\text{PCSG}}$ and Loss-Net $\hat{E}$ 
  \end{algorithmic}
\end{algorithm}

\begin{figure}
  \FIGURE
  { \includegraphics[width=\textwidth]{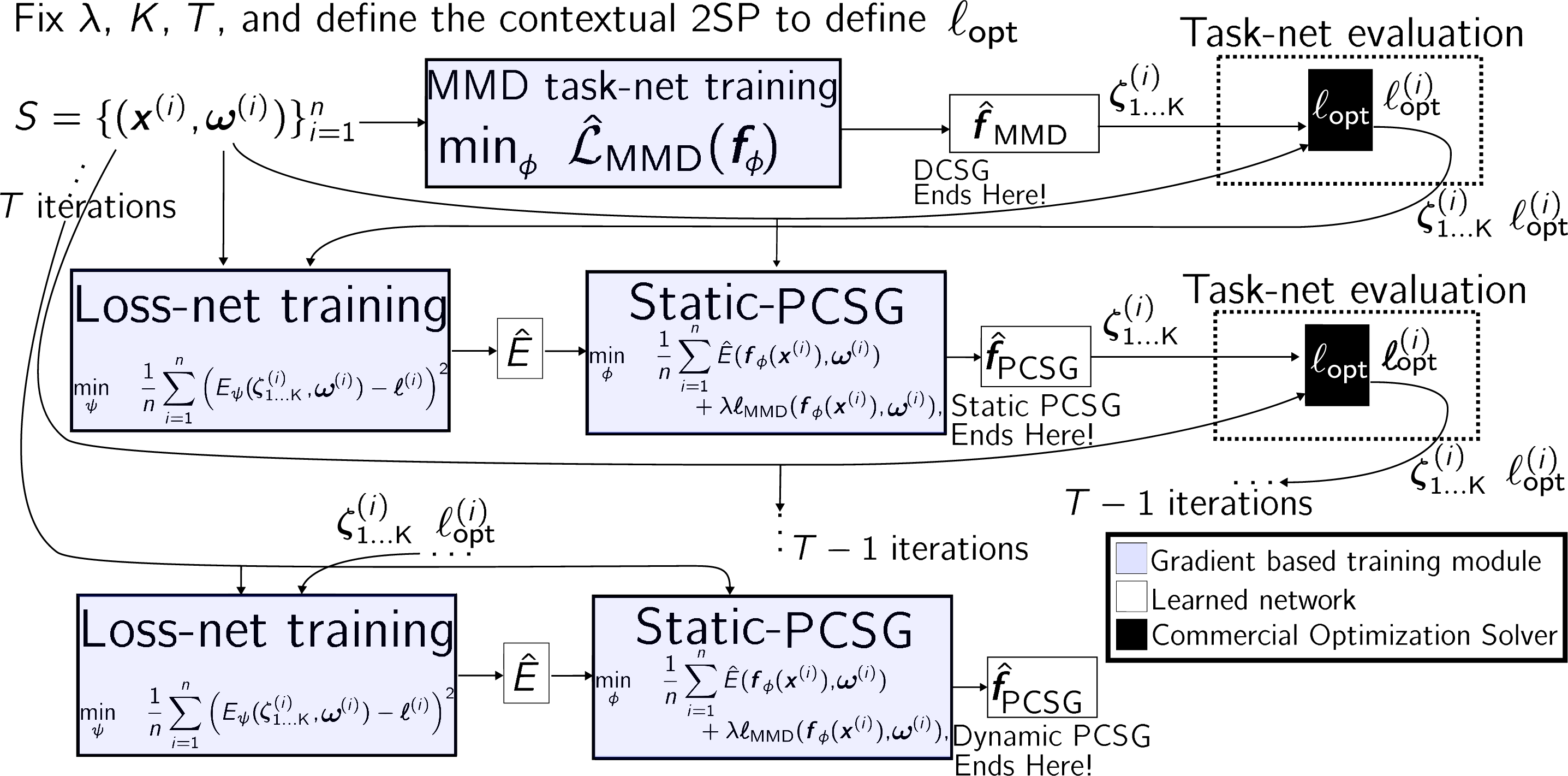} }
  {Visual representation of the DCSG and PCSG training procedures \label{fig:overall training}}
  {}
\end{figure}

Algorithm \ref{algo:dynamic-training} is similar to the approach proposed by \citet{zharmagambetov2023landscape}. In the case of scenario generation, we observe that \citet{zharmagambetov2023landscape}'s approach fails without MMD regularization. Thus, the implicit use of MMD to initialize the samples for $E_{\psi}$ and regularize the Task-Net $\bs{f}_{\phi}$ constitutes a critical difference between this work and theirs. Similar to \citet{zharmagambetov2023landscape}, the proposed methods do not maintain the entire history of examples for training $E_{\psi}$ at every iteration and instead, a \textit{replay buffer} \citep{lin1992self} is used, keeping, at most, the last $T = 3$ iterations of data. Furthermore, Algorithm \ref{algo:dynamic-training} builds upon Algorithm \ref{alg:static}, which relies upon empirical minimization of $\mathcal{L}_{\text{MMD}}(\bs{f}_{\phi})$. Figure \ref{fig:overall training} visually displays the entire training process and the relationships between DCSG, Static PCSG, and Dynamic PCSG. 

The experimental results will evaluate the performance of these three approaches relative to each other. In terms of training, minimizing $\hat{\mathcal{L}}_{\text{MMD}}(\bs{f}_{\phi})$ has the advantage that no optimization problems are required to evaluate the loss. However, the MMD approach does not account for the problem structure in \eqref{eq:2SP-SAA}. The static and dynamic approaches have the advantage of considering the downstream decision loss associated with a particular choice of $\bs{f}$. The static approach has the benefit over the dynamic approach that the training of $E_{\psi}$, and thus the repeated solution of optimization problems to evaluate $\ell_{\text{opt}}$, need only be performed once. While the dynamic approach offers an advantage over the static approach by aiming to construct better approximations of $\ell_{\text{opt}}$ near the chosen $\bs{f}$, it also has a drawback in that computations of $\ell_{\text{opt}}$ must be performed at each iteration to update the Loss-Net.
\section{Performance Guarantees}\label{section:guarantees}
In this section, we present finite-sample performance guarantees for our proposed methodologies. Sections~\ref{subsection:theory_dsg} and \ref{subsection:theory_PCSG} provide finite-sample generalization bounds based on Rademacher complexities for both $\ell_{\text{MMD}}$ and $\ell_{\text{opt}}$ when learning over classes of $\mathcal{F}$ with bounded Rademacher complexity. These results provide theoretical justification for learning $\bs{f}$ using our proposed methodology. 

A natural question is whether we can guarantee performance in $\ell_{\text{opt}}$ when we do not optimize $\ell_{\text{opt}}$ directly, but instead optimize its surrogate $E_{\psi}$. Section \ref{subsection:theory_PCSG} also extends the generalization results for $\ell_{\text{opt}}$ to this setting and relates the resulting bound to our methodology.

Our finite-sample guarantees mirror the generalization bounds obtained by \citet{grigas2021integrated}. The primary difference is that we rely on boundedness and Lipschitz properties of our loss functions, whereas they assume these properties for optimal solution maps. Unlike \citet{grigas2021integrated}, who do not consider a problem-agnostic approach, we relate \ref{eq:DCSG} to our problem-driven approaches.

\subsection{Statistical Learning Preliminaries}
We rely on single and multivariate Rademacher complexities to provide finite-sample generalization bounds. The Rademacher complexity measures the richness of a hypothesis class with respect to a probability distribution and is defined as follows.

\begin{samepage}
\begin{definition}[Rademacher Averages]
  Let $\mathcal{G}$ be a family of functions mapping $\Xi \subseteq \mathbb{R}^{\hat{p}}$ to $\mathbb{R}$, for some $\hat{p} \in \mathbb{N}$ and $S = \{\bs{\xi}^{(i)}\}_{i=1}^n$ denote a fixed sample of size $n$, with $\bs{\xi}^{(i)} \in \Xi, \ i \in [n]$. Then, the empirical Rademacher complexity of $\mathcal{G}$ with respect to $S$ is defined as
  $$
  \hat{R}_{n}(\mathcal{G}, S) 
  \coloneq \mathbb{E}_{\bs{\sigma}} \left[ \sup_{g \in \mathcal{G}} \frac{1}{n}\sum_{i=1}^n \sigma_i g(\bs{\xi}^{(i)}) \right],
  $$
  where $\sigma_1, \hdots \sigma_n$ denote independent uniform random variables taking values in $\{-1, +1\}$. 

  Furthermore, one can define the empirical vector Rademacher complexity by considering classes of vector-valued functions $\mathcal{F} \subseteq \{\bs{f}:\Xi \rightarrow \mathbb{R}^K\}$. The empirical Rademacher complexity of the class of vector-valued functions $\mathcal{F}$ with respect to $S$ is then defined as
  $$
  \hat{R}_{n}(\mathcal{F}, S) 
  \coloneq \mathbb{E}_{\bs{\sigma}} \left[ \sup_{\bs{f} \in \mathcal{F}} \frac{1}{n} \sum_{i=1}^n \sum_{j=1}^K \sigma_{ij} f_j(\bs{\xi}^{(i)}) 
  \right],$$
  where $\sigma_{ij}$ are also independent uniform variables taking values in $\{-1, +1\}$ \hbox{$\forall \ i \in [n]$}, $j\in [K]$, and $f_j(\bs{\xi}^{(i)})$ is the $j$th component of $\bs{f}(\bs{\xi}^{(i)})$. We refer to $\bs{\sigma}$ as the Rademacher variables.

  For both the scalar and vector-valued cases, the Rademacher complexity is given by
  $$
  R_{n}(\mathcal{G}) 
  \coloneq \mathbb{E}_{S} \left[\hat{R}_{n}(\mathcal{G}, S)\right].
  $$
\end{definition}
\end{samepage}

The Rademacher averages allow one to uniformly bound, over a function class, expectations in terms of empirical averages. We present the following theorem that is useful in providing performance bounds on $\ell_{\text{opt}}$ and $\ell_{\text{MMD}}$.

\begin{restatable}{theorem}{multivariaterademacher}
  \label{theorem:multivariate}
Let $\ell : \Omega^{K+1} \rightarrow \mathbb{R}$ be bounded $\sup_{\bs{\zeta}_{1 \hdots K}, \bs{\omega} \in \Omega^{K+1}}|\ell(\bs{\zeta}_{1 \hdots K}, \bs{\omega})| \leq M$, such that for any fixed $\bs{\omega}$, \hbox{$\ell_{\bs{\omega}}: \bs{\zeta}_{1 \hdots K} \mapsto  \ell(\bs{\zeta}_{1 \hdots K}, \bs{\omega})$} is $L$-Lipschitz. Also, let $\mathcal{F}$ be the vector-valued class of functions $\bs{f} : \mathcal{X} \rightarrow \Omega^K$. Then for any $\delta > 0$, with probability at least $1-\delta$ over the draw of an iid sample $S = \{(\bs{x}^{(i)}, \bs{\omega}^{(i)})\}_{i=1}^n$ of size $n$, the following holds for all $\bs{f} \in \mathcal{F}$:
\begin{equation*}
  \mathbb{E}_{(\bs{x}, \bs{\omega}) \sim \mathbb{P}_{\bs{x}, \bs{\omega}}}[\ell(\bs{f}(\bs{x}), \bs{\omega})] \leq \frac{1}{n} \sum_{i=1}^{n} \ell(\bs{f}(\bs{x}^{(i)}), \bs{\omega}^{(i)}) + 2\sqrt{2} L R_n(\mathcal{F}) + M\sqrt{\frac{\log(1/\delta)}{2n}}
  \end{equation*}
  % and
  % \begin{equation*}
  %   \mathbb{E}_{(\bs{x}, \bs{\omega}) \sim \mathbb{P}_{\bs{x}, \bs{\omega}}}[\ell(\bs{f}(\bs{x}), \bs{\omega})] \leq \frac{1}{n} \sum_{i=1}^{n} \ell(\bs{f}(\bs{x}^{(i)}), \bs{\omega}^{(i)}) + 2\sqrt{2} L \hat{R}_n(\mathcal{F}, S) + 3 M\sqrt{\frac{\log(2/\delta)}{2n}}
  %   \end{equation*}
  where $R_n(\mathcal{F})$ denotes the Rademacher average over the vector-valued function class $\mathcal{F}$, and we consider the elements of $\Omega^K$ as vectors in $\mathbb{R}^{pK}$.
\end{restatable}

The proof of Theorem \ref{theorem:multivariate} is provided in Section \ref{section:appendix_learning_bounds}. The proof relies on standard techniques for bounding Rademacher averages of linear hypothesis classes \citep{mohri2018foundations,bartlett2002rademacher,maurer2016vector}. From Theorem \ref{theorem:multivariate} we conclude that, if (i) $\hat{R}_{n}(\mathcal{F}, S)$ asymptotically vanishes as $n \rightarrow \infty$, (ii) $\ell$ is bounded, and (iii) $\ell_{\bs{\omega}}$ is Lipschitz $\forall \ \bs{\omega} \in \Omega$, then the task-mapping achieves generalization on unseen data.

To demonstrate an example of an $\mathcal{F}$, in Section \ref{appendix:taskmapping} we provide a hypothesis class consisting of $K$ linear transformations of $\bs{x}$, composed with a non-linear function. Based on this choice of hypothesis class, $O(\sqrt{1/n})$ upper bounds on $\hat{R}_{n}(\mathcal{F}, S)$ are presented. Although in practice we do not parametrize $\mathcal{F}$ as described above and instead use neural networks, we confirm that there are hypothesis classes with bounded complexity for our scenario prediction problem. As such, we verify boundedness + Lipschitzness for $\ell_{\text{MMD}}$ and $\ell_{\text{opt}}$ next.

\subsection{Finite-Sample Guarantees for DCSG} \label{subsection:theory_dsg}

We proceed with the following lemma that provides an upper bound on the Lipschitz constant and the absolute value of $\ell_{\text{MMD}, \bs{\omega}}: \bs{\zeta}_{1 \hdots K} \mapsto \ell_{\text{MMD}}(\bs{\zeta}_{1 \hdots K}, \bs{\omega})$. Let $L_k$ denote a bound on the Lipschitz constant of the positive semi-definite kernel function.
\begin{restatable}{lemma}{kernellipschitz}\label{lemma:kernel_lipschitz}
  For any fixed $\bs{\omega} \in \Omega$, $\ell_{\text{MMD}, \bs{\omega}}$ has Lipschitz constant upper bounded by $L_{\text{MMD}} = 4L_k/\sqrt{K}$. Furthermore, if $k(\bs{\omega}, \bs{\omega}') \leq U \ \forall \ \bs{\omega}, \bs{\omega}' \in \Omega$, then $\ell_{\text{MMD}}$ is bounded by $2U$.
\end{restatable}
The proof of Lemma \ref{lemma:kernel_lipschitz} is provided in Section \ref{subsection:mmdloss}. Thus, according to \hbox{Theorem \ref{theorem:multivariate}}, if $k$ is bounded and Lipschitz, then one can provide finite-sample generalization bounds for learning $\bs{f}$ using $\ell_{\text{MMD}}$.

Recall that the energy kernel $k_E$ is equivalent to the negative Euclidean kernel $\tilde{k}_E(\bs{\omega}, \bs{\omega}') = -\norm{\bs{\omega} - \bs{\omega}'}_2$. Setting $k = \tilde{k}_E$ yields a Lipschitz kernel, because $-\norm{\cdot - \bs{\omega}'}_2$ is $1$-Lipschitz by the reverse-triangle inequality: Let $g_{\bs{\omega}'}: \Omega \rightarrow \mathbb{R}_+$ be such that $g_{\bs{\omega}'}(\bs{\omega}) = \norm{\bs{\omega} - \bs{\omega}'}_2$. Then, for any $\bs{\omega} \in \Omega $ and $\bs{\zeta} \in \Omega$:
\begin{align*}
\left| g_{\bs{\omega}'}(\bs{\omega}) - g_{\bs{\omega}'}(\bs{\zeta}) \right| & = \left| \norm{\bs{\zeta} - \bs{\omega}'}_2  -\norm{\bs{\omega} - \bs{\omega}'}_2 \right| \\
& \leq \norm{ \bs{\zeta} - \bs{\omega}' - \left(\bs{\omega} - \bs{\omega}'\right)}_2  = \norm{\bs{\zeta} - \bs{\omega}}_2,
\end{align*}
implying $-\norm{\cdot - \bs{\omega}'}_2$ is $1$-Lipschitz. Furthermore, the steps in the proof of Lemma \ref{lemma:kernel_lipschitz} hold equivalently if the bounded kernel condition is replaced with $\norm{\bs{\omega} - \bs{\omega}'}_2 \leq U$. 

\subsection{Finite-Sample Guarantees for PCSG}\label{subsection:theory_PCSG}

To illustrate finite-sample guarantees for our proposed problem-driven approach, we focus on the class of linear, fixed-recourse 2SPs in which uncertainty appears only on the right-hand side of the second-stage constraints. This setup is common in the stochastic programming literature; see, e.g., \citet{gose2016sequential} and \citet{dentcheva2000differential}. In this setting, the first-stage cost and second-stage cost are given by $h(\bs{y})=\bs{h}^\top \bs{y}$ and $q(\bs{\omega},\bs{z})=\bs{q}^\top \bs{z}$, respectively. The first- and second-stage feasible regions are given by
\[
\mathcal{Y} = \Big\{\bs{y}\in\mathbb{R}^{s_1}_+:\ \bs{A} \bs{y} \geq \bs{b},\  \bs{y} \leq \bs{e} \Big\} \neq \emptyset, \ \text{and} \ \mathcal{Z}(\bs{y},\bs{\omega}) = \Big\{\bs{z}\in\mathbb{R}^{s_2}_+:\ \bs{C}\bs{y} + \bs{D}\bs{z}=\bs{\omega} \Big\},
\]
where $\bs{A} \in \mathbb{R}^{p^{\prime} \times s_1}$ for some $p^{\prime} \in \mathbb{N}$, $\bs{e} > \bs{0}$, $\bs{C} \in \mathbb{R}^{p \times s_1}$, and $\bs{D} \in \mathbb{R}^{p \times s_2}$. ($\bs{\zeta}$-SAA) is then given by
\begin{align}
v^*(\bs{\zeta}_{1 \hdots K})
=\min_{\bs{y}\ge \bs{0},\,\bs{z}_1,\dots,\bs{z}_K\ge \bs{0}}\quad
& \bs{h}^\top \bs{y} + \frac{1}{K}\sum_{i=1}^K \bs{q}^\top \bs{z}_i   \tag{$\bs{\zeta}$-SAA-RHS} \label{eq:zeta_saa_rhs_only}\\
\textit{s.t.}\quad
& \bs{A}\bs{y} \ge \bs{b}, \ \bs{y} \leq \bs{e} \notag \\
& \bs{C}\bs{y} + \bs{D}\bs{z}_i = \bs{\zeta}_i,\quad \forall i \in \{1,\dots,K\}. \notag
\end{align}
Similarly, $\ell_{\text{opt}}$ can be evaluated by solving the following instance of \eqref{eq:opt_search}:
\begin{align}
\ell_{\text{opt}}(\bs{\zeta}_{1 \hdots K},\bs{\omega})
=\min_{\bs{y}\ge \bs{0},\,\bs{z}\ge \bs{0},\,\bs{z}_1,\dots,\bs{z}_K\ge \bs{0}} \quad
& \bs{h}^\top \bs{y} + \bs{q}^\top \bs{z} \tag{\text{Opt-Search-RHS}} \label{eq:opt_search_rhs_only}\\
\textit{s.t.}\quad
& \bs{A}\bs{y} \ge \bs{b}, \ \bs{y} \leq \bs{e} \notag\\
& \bs{C}\bs{y} + \bs{D}\bs{z} = \bs{\omega}, \notag \\
& \bs{C}\bs{y} + \bs{D}\bs{z}_i = \bs{\zeta}_i,\quad \forall i \in \{1,\dots,K\}, \notag \\
& \bs{h}^\top \bs{y} + \frac{1}{K}\sum_{i=1}^K \bs{q}^\top \bs{z}_i \le v^*(\bs{\zeta}_{1 \hdots K}). \notag
\end{align}
Furthermore, we make the following assumptions:
\begin{assumption}\label{ass:complete_recourse}
  It is assumed that complete recourse holds, i.e., there exists $\bs{z} \geq \bs{0}$ such that $\bs{D}\bs{z} = \bs{\eta}$ for all $\bs{\eta} \in \mathbb{R}^p$.
\end{assumption}

\begin{assumption}\label{ass:dual_feasibility}
  It is assumed that there exists $\bs{\pi}$ such that $\bs{D}^\top \bs{\pi} \le \bs{q}$.
\end{assumption}

Assumption \ref{ass:complete_recourse} ensures that both \eqref{eq:zeta_saa_rhs_only} and \eqref{eq:opt_search_rhs_only} are feasible for all choices of $\bs{\zeta}_{1 \hdots K}$ and $\bs{\omega}$. Assumption \ref{ass:dual_feasibility} ensures the dual programs associated with \eqref{eq:zeta_saa_rhs_only} and \eqref{eq:opt_search_rhs_only} contain feasible solutions. In Section~\ref{subsection:optloss} we show that Assumption \ref{ass:dual_feasibility}, along with the upper bounds $\bs{e}$ placed on $\bs{y}$, prevent dual infeasibility. \citet{birge1997introduction} point out that the two aforementioned assumptions are typically associated with well-defined models.

Next, we present a Lemma that provides an upper bound on the Lipschitz constant of $\ell_{\text{opt}, \bs{\omega}}: \bs{\zeta}_{1 \hdots K} \mapsto \ell_{\text{opt}}(\bs{\zeta}_{1 \hdots K}, \bs{\omega})$ and a uniform bound on $\ell_{\text{opt}}$.

\begin{restatable}{lemma}{optlipschitz}\label{lemma:opt_lipschitz}
  % Consider \eqref{eq:zeta_saa_rhs_only} and \eqref{eq:opt_search_rhs_only} under Assumptions \ref{ass:complete_recourse} and \ref{ass:dual_feasibility}. It follows that $v^*$ is Lipschitz continuous with respect to $\bs{\zeta}_{1 \hdots K}$ for some constant $L_v > 0$. Furthermore, $\ell_{\text{opt}}$ has Lipschitz constant upper bounded by $L_{\text{opt}} = L^{\prime} \sqrt{L_v^2 + 1}$ for some $L^{\prime} > 0$. Thus,
  % for any fixed $\bs{\omega} \in \Omega$, $\ell_{\text{opt}, \bs{\omega}}$ has Lipschitz constant upper bounded by $L_{\text{opt}}$ also. Lastly, if $\Omega$ is bounded such that $\sup_{\bs{\omega} \in \Omega} \|\bs{\omega}\|_2 \leq U_{\Omega}$, then there exists $U_{\text{opt}} > 0$ such that $\ell_{\text{opt}}$ is uniformly bounded by $U_{\text{opt}}$ for any $\bs{\zeta}_{1 \hdots K} \in \Omega^K,  \bs{\omega} \in \Omega$.
  Consider \eqref{eq:zeta_saa_rhs_only} and \eqref{eq:opt_search_rhs_only} under Assumptions \ref{ass:complete_recourse} and \ref{ass:dual_feasibility}. It follows that $v^*$ is Lipschitz continuous with respect to $\bs{\zeta}_{1 \hdots K}$ with constant $L_v>0$. Furthermore, $\ell_{\text{opt}}$ is Lipschitz continuous with constant
\[
L_{\text{opt}} := L'\sqrt{L_v^2+1},
\]
for some $L'>0$. Hence, for any fixed $\bs{\omega}\in\Omega$, the map $\ell_{\text{opt},\bs{\omega}}$ is also Lipschitz continuous with constant at most $L_{\text{opt}}$. Lastly, if $\Omega$ is bounded so that $\sup_{\bs{\omega}\in\Omega}\|\bs{\omega}\|_2 \le U_{\Omega}$, then there exists $U_{\text{opt}}>0$ such that
\[
|\ell_{\text{opt}}(\bs{\zeta}_{1 \hdots K},\bs{\omega})| \le U_{\text{opt}},
\qquad \forall\,\bs{\zeta}_{1 \hdots K} \in\Omega^K, \bs{\omega} \in \Omega.
\]

\end{restatable}

The proof of Lemma \ref{lemma:opt_lipschitz} is provided in Section \ref{subsection:optloss}. Thus, according to \hbox{Theorem \ref{theorem:multivariate}}, in the case of linear 2SPs with right-hand side uncertainty, complete recourse and dual-feasibility for the optimal response, and bounded uncertainty set $\Omega$, we obtain finite-sample performance guarantees uniformly over $\bs{f}$ using $\ell_{\text{opt}}$. Thus, for $\bs{f}$ obtained via any training methodology we can conclude, given enough samples $n$, the empirical estimate of our task loss $\ell_{\text{opt}}$ converges to its expected value over $\mathbb{P}_{\bs{x}, \bs{\omega}}.$

The performance guarantees obtained by Theorem \ref{theorem:multivariate} and  \hbox{Lemma \ref{lemma:opt_lipschitz}} apply uniformly across $\mathcal{F}$, nonetheless, they do not explicitly account for the fact that an approximation of $\ell_{\text{opt}}$ is selected from a function class $\mathcal{G}_{\varPsi} \coloneq \{E_{\psi}: \psi \in \varPsi \}$ to best approximate $\ell_{\text{opt}}$ over a distribution $\mathbb{P}_{\bs{\zeta}_{1 \hdots K}, \bs{\omega}}$ constructed from the $n$ observations $S$. 
\begin{samepage}
\begin{restatable}{lemma}{losserror}\label{lemma:losserror}
  Consider the assumptions in Lemma \ref{lemma:opt_lipschitz} so that $\ell_{\text{opt}}$ is bounded: $\sup_{\bs{\zeta}_{1 \hdots K}, \bs{\omega} \in \Omega^{K+1}} |\ell_{\text{opt}}(\bs{\zeta}_{1 \hdots K}, \bs{\omega} )| \leq U_{\text{opt}}$. Further, assume $\mathcal{G}_{\varPsi}$ is bounded: $\sup_{\bs{\zeta}_{1 \hdots K}, \bs{\omega} \in \Omega^{K+1}} \sup_{g \in \mathcal{G}_{\varPsi}} |g(\bs{\zeta}_{1 \hdots K}, \bs{\omega}) | \leq U_{\varPsi}$ Then for any $\mathbb{P}_{\bs{\zeta}_{1 \hdots K}, \bs{\omega}}$ and $\delta > 0$, with probability at least $1 - \delta$ over the draw of an iid sample $S'' = \{(\bs{\zeta}_{1 \hdots K}^{(i)}, \bs{\omega}^{(i)}) \}_{i=1}^m$ of size $m$, the following holds for all $g \in \mathcal{G}_{\varPsi}$:
  \begin{align*}
  &\mathbb{E}_{(\bs{\zeta}_{1 \hdots K}, \bs{\omega}) \sim \mathbb{P}_{\bs{\zeta}_{1 \hdots K}, \bs{\omega}}}
  \Big[ \big| \ell_{\text{opt}}(\bs{\zeta}_{1 \hdots K}, \bs{\omega}) - g(\bs{\zeta}_{1 \hdots K}, \bs{\omega}) \big| \Big] \\
  &\leq \frac{1}{m} \sum_{i=1}^{m}
  \big| \ell_{\text{opt}}(\bs{\zeta}^{(i)}_{1 \hdots K}, \bs{\omega}^{(i)}) - g(\bs{\zeta}^{(i)}_{1 \hdots K}, \bs{\omega}^{(i)}) \big|
  \;+\; 2 R_m(\mathcal{G}_{\varPsi}) + (U_{\text{opt}} + U_{\varPsi}) \sqrt{\frac{\log{1/\delta}}{2m}}.
  \end{align*}
\end{restatable}
\end{samepage}
\begin{proof}
The proof is immediate from Theorem 11.3 in \citet{mohri2018foundations} applied to the $\ell_1$ loss function. The theorem applies since $\big| \ell_{\text{opt}}(\bs{\zeta}_{1 \hdots K}, \bs{\omega}) - g(\bs{\zeta}_{1 \hdots K}, \bs{\omega}) \big| \leq U_{\text{opt}} + U_{\varPsi}, \forall \bs{\zeta}_{1 \hdots K}, \bs{\omega} \in \Omega^{K+1}$. 
\end{proof}
In Lemma \ref{lemma:losserror}, we work with a reference distribution $\mathbb{P}_{\bs{\zeta}_{1 \hdots K}, \bs{\omega}}$ on $\Omega^{K+1}$ that we use to generate samples for a regression method. The bound depends on $R_m(\mathcal{G}_{\varPsi})$ and the upper bounds on $\ell_{\text{opt}}$ and $\mathcal{G}_{\varPsi}$. For example, setting $\mathcal{G}_{\varPsi}$ to a kernel-based hypothesis with bounded RKHS-norm and bounded kernel controls $R_m(\mathcal{G}_{\varPsi})$ and $U_{\varPsi}$ (see Theorem 6 in \citet{mohri2018foundations} for more details). Similarly, the problem structure defined by \eqref{eq:zeta_saa_rhs_only} and \eqref{eq:opt_search_rhs_only} yields the bound $U_{\text{opt}}$.

Using Lemma \ref{lemma:losserror}, Lemma \ref{lemma:opt_lipschitz} and Theorem \ref{theorem:multivariate}, we obtain a finite-sample performance guarantee for our proposed methodology \eqref{eq:static_PCSG} that applies in the idealized case where the samples used to select $\bs{f} \in \mathcal{F}$ and $g \in \mathcal{G}_{\varPsi}$ are independent coming from $\mathbb{P}_{\bs{x}, \bs{\omega}}$ and $\mathbb{P}_{\bs{\zeta}_{1 \hdots K}, \bs{\omega}}$ respectively. We use the notation $\mathbb{P}_{\bs{f}(\bs{x})}$ to refer to the distribution (pushforward) on $\Omega^K$ induced by the marginal distribution $\mathbb{P}_{\bs{x}}$ for all $\bs{f} \in \mathcal{F}$. Similarly, we use $\mathbb{P}_{\bs{f}(\bs{x}), \bs{\omega}}$ to refer to the joint distribution on $\Omega^{K+1}$ induced by $\mathbb{P}_{\bs{x}, \bs{\omega}}$.

Essentially, we uniformly bound the expectation of $\ell_{\text{opt}}$ by the sum of five parts: (i) the sample average of the approximation of $\ell_{\text{opt}}$, (ii) the complexity of the hypothesis class $\mathcal{F}$ used for the task mapping, (iii) the complexity of the hypothesis class $\mathcal{G}$ used for the approximation of $\ell_{\text{opt}}$, (iv) the error associated with approximating $\ell_{\text{opt}}$ using $\mathcal{G}$, and (v) the 1-Wasserstein distance between $\mathbb{P}_{\bs{\zeta}_{1 \hdots K}, \bs{\omega}}$ and $\mathbb{P}_{\bs{f}(\bs{x}), \bs{\omega}}$.

\begin{samepage}
\begin{restatable}{theorem}{appxlossguarantee}\label{thm:appxlossguarantee}
  Consider \eqref{eq:zeta_saa_rhs_only} and \eqref{eq:opt_search_rhs_only} under
  Assumptions \ref{ass:complete_recourse} and \ref{ass:dual_feasibility}, and boundedness of $\Omega$.
  Furthermore, consider a function class $\mathcal{G}_{\varPsi}$ with bounded Rademacher complexity such that
  $g$ is bounded and Lipschitz for all $g \in \mathcal{G}_{\varPsi}$. Suppose the approximate task loss $\hat{E}$ is selected
  from $\mathcal{G}_{\varPsi}$. Then, with probability at least $1-\delta$ over the draw of the iid sample
  $S=\{(\bs{x}^{(i)},\bs{\omega}^{(i)})\}_{i=1}^n$ of size $n$, drawn from $\mathbb{P}_{\bs{x},\bs{\omega}}$,
  and over an independent sample
  $S''=\{(\bs{\zeta}^{(i)}_{1 \hdots K},\bs{\omega}^{(i)})\}_{i=1}^m$ of size $m$, drawn from
  $\mathbb{P}_{\bs{\zeta}_{1 \hdots K},\bs{\omega}}$, the following holds for all $\bs{f} \in \mathcal{F}$:
  \begin{align}
    \mathbb{E}_{(\bs{x},\bs{\omega}) \sim \mathbb{P}_{\bs{x},\bs{\omega}}}
    \big[\ell_{\text{opt}}(\bs{f}(\bs{x}), &\bs{\omega})\big]
    \leq
    \frac{1}{n}\sum_{i=1}^{n} \hat{E}\big(\bs{f}(\bs{x}^{(i)}),\bs{\omega}^{(i)}\big)
    + 2\sqrt{2}\,L_{\text{opt}}\,R_n(\mathcal{F}) \notag
    \\
    &+ U_{\text{opt}}\sqrt{\frac{\log(2/\delta)}{2n}}
    +
    \frac{1}{m}\sum_{i=1}^{m}
    \Big|\ell_{\text{opt}}\big(\bs{\zeta}^{(i)}_{1 \hdots K},\bs{\omega}^{(i)}\big)
          - \hat{E}\big(\bs{\zeta}^{(i)}_{1 \hdots K},\bs{\omega}^{(i)}\big)\Big| \notag
    \\
    &+
    2\Big(R_m(\mathcal{G}_{\varPsi}) + R_n(\mathcal{G}_{\varPsi})\Big) 
    + (L_{\text{opt}} + L_{\varPsi}) \ d_{\text{Lip}}(\mathbb{P}_{\bs{\zeta}_{1 \hdots K}, \bs{\omega}}, \mathbb{P}_{\bs{f}(\bs{x}), \bs{\omega}}) \notag
    \\
    &+ (U_{\text{opt}} + U_{\varPsi})
    \left(
      \sqrt{\frac{\log(4/\delta)}{2m}}
      + \sqrt{\frac{\log(4/\delta)}{2n}}
    \right), \label{eq:appx:bound}
  \end{align}
  where $L_{\text{opt}}$ and $U_{\text{opt}}$ are as in Lemma \ref{lemma:opt_lipschitz}, $U_{\varPsi}$ is as in Lemma \ref{lemma:losserror}, and $L_{\varPsi}$ is the assumed Lipschitz constant of $\mathcal{G}_{\varPsi}$.
\end{restatable}
\end{samepage}

The proof of Theorem \ref{thm:appxlossguarantee} is provided in Section \ref{subsection:optloss}. The last line of \eqref{eq:appx:bound} involves the distributional shift term
$d_{\text{Lip}}(\mathbb{P}_{\bs{\zeta}_{1 \hdots K},\bs{\omega}},\mathbb{P}_{\bs{f}(\bs{x}),\bs{\omega}})$, which measures how different the
Loss-Net training distribution on $(\bs{\zeta}_{1 \hdots K},\bs{\omega})$ is from the distribution induced by the learned task-map $\bs{f}$.
To heuristically control this shift, we replace the 1-Wasserstein distance with an integral probability metric over a less expressive function class, namely an MMD induced by a
carefully constructed kernel on $\Omega^{K+1}$. In Section~\ref{appx:joint_kernel_mmd}, we show that any base kernel
$k$ on $\Omega$ with feature map $\bs{\varphi}:\Omega\to\mathcal{H}$ induces a positive semi-definite kernel $\kappa$ on $\Omega^{K+1}$
via the feature map
$$
\bs{\varsigma}(\bs{\zeta}_{1 \hdots K},\bs{\omega})
:=\frac{1}{K}\sum_{i=1}^K \bs{\varphi}(\bs{\zeta}_i)\;-\;\bs{\varphi}(\bs{\omega})
\ \in\ \mathcal{H},
$$
and associated kernel: 
$
\kappa\big((\bs{\zeta}_{1 \hdots K},\bs{\omega}),(\bs{\zeta}'_{1 \hdots K},\bs{\omega}')\big)
:=\big\langle \bs{\varsigma}(\bs{\zeta}_{1 \hdots K},\bs{\omega}),\,\bs{\varsigma}(\bs{\zeta}'_{1 \hdots K},\bs{\omega}')\big\rangle_{\mathcal{H}}.
$
We show that the associated kernel $\kappa$ is permutation-invariant in the first $K$ coordinates and inherits positive semi-definiteness from $k$. Furthermore, we that see this choice of $\kappa$ has desirable properties under a particular reference distribution
$\mathbb{P}_{\bs{\zeta}_{1 \hdots K},\bs{\omega}}$ that we refer to as the \textit{oracle distribution} $\mathbb{P}^*_{\bs{\zeta}_{1 \hdots K},\bs{\omega}}$ on $\Omega^{K+1}$. The oracle distribution is defined by first sampling $\bs{x}$ from $\mathbb{P}_{\bs{x}}$, then sampling $\bs{\omega}$ from $\mathbb{P}_{\bs{\omega}\mid \bs{x}}$, and finally sampling $\bs{\zeta}_1,\ldots,\bs{\zeta}_K$ independently and identically from $\mathbb{P}_{\bs{\omega}\mid \bs{x}}$, yielding a joint distribution 
$\mathbb{P}^*_{\bs{\zeta}_{1 \hdots K},\bs{\omega}}$ on $\Omega^{K+1}$. In this setting, Section~\ref{appx:joint_kernel_mmd} shows that the
contextual MMD objective $\mathcal{L}_{\mathrm{MMD}}$ with kernel $k$ upper bounds the MMD distance between the induced joint distribution
$\mathbb{P}_{\bs{f}(\bs{x}),\bs{\omega}}$ and the oracle joint distribution $\mathbb{P}^*_{\bs{\zeta}_{1 \hdots K},\bs{\omega}}$ measured using
the induced kernel $\kappa$. 
% \[
% d^2_{\kappa}\!\Big(\mathbb{P}^*_{\bs{\zeta}_{1 \hdots K},\bs{\omega}},\,\mathbb{P}_{\bs{f}(\bs{x}),\bs{\omega}}\Big)
% \;\le\;
% \mathcal{L}_{\mathrm{MMD}}(\bs{f}).
% \]
Thus, minimizing $\mathcal{L}_{\mathrm{MMD}}$ over $f \in \mathcal{F}$ controls a MMD distance on joint-distributions defined on $\Omega^{K+1}$, where the reference distribution follows a reasonable structure. In light of Theorem~\ref{thm:appxlossguarantee}, the joint MMD under $\kappa$ provides a computationally tractable heuristic for controlling the 1-Wasserstein distribution shift term $d_{\text{Lip}}(\mathbb{P}_{\bs{\zeta}_{1 \hdots K},\bs{\omega}},\mathbb{P}_{\bs{f}(\bs{x}),\bs{\omega}})$ appearing in \eqref{eq:appx:bound}; while the MMD function class is typically not as rich as the Wasserstein class, the kernel-based IPM offers substantial computational advantages (as discussed in Section~\ref{appendix:comparing}).

Thus, the proposed task-driven methodologies can be viewed as heuristics for controlling the terms in Theorem~\ref{thm:appxlossguarantee}. First, by selecting $\bs{f}$ to optimize $\mathcal{L}_{\text{MMD}}$ we are attempting to approximate the oracle distribution $\mathbb{P}^*_{\bs{\zeta}_{1 \hdots K},\bs{\omega}}$. Given a learned mapping $\hat{\bs{f}}_{\text{MMD}}$, we then use it, along with the optimization solver to generate $\{\bs{\zeta}_{1 \hdots K}^{(i)}, \bs{\omega}^{(i)}, \ell_{\text{opt}}^{(i)} \}_{i=1}^n$; this can be interpreted as a way to decrease the error term associated with $\hat{E}$ that appears in \eqref{eq:appx:bound}. Finally, reselecting $\bs{f}$ with the aim of minimizing the average predicted task loss regularized by the MMD loss aims to control the first term and the 1-Wasserstein term in \eqref{eq:appx:bound}. 

Lastly, we note that the independence assumption used to derive Theorem~\ref{thm:appxlossguarantee} does not hold. In Algorithms \ref{alg:static} and \ref{algo:dynamic-training}, the samples $S$ are used to produce $\hat{\bs{f}}_{\text{MMD}}$, which is then used to generate the training data for $\hat{E}$. Thus, the samples used to train $\hat{E}$ are not independent from the samples used to select $\hat{\bs{f}}_{\text{MMD}}$. This is the case even if $S$ is split into two disjoint sets (one for $\bs{f}$ and one for $\hat{E}$), because the training data for $\hat{E}$ depends on $\hat{\bs{f}}_{\text{MMD}}$, which in turn depends on the samples $S$ used to select $\hat{\bs{f}}_{\text{MMD}}$. This dependence structure is not accounted for in Theorem~\ref{thm:appxlossguarantee}, and thus, the bound does not hold with high probability over the draw of $S$ alone but holds over the draw of $S$ and $S''$ sampled from $\mathbb{P}_{\bs{x}, \bs{\omega}}$ and $\mathbb{P}_{\bs{\zeta}_{1 \hdots K},\bs{\omega}}$, respectively. In our setting, we do not know or have any samples from $\mathbb{P}_{\bs{\zeta}_{1 \hdots K},\bs{\omega}}$ and so we resort selecting $f \in \mathcal{F}$ according to $\mathcal{L}_{\text{MMD}}$, which controls the MMD distance between $\mathbb{P}_{\bs{f}(\bs{x}),\bs{\omega}}$ and $\mathbb{P}^*_{\bs{\zeta}_{1 \hdots K},\bs{\omega}}$ under kernel $\kappa$. Optimizing  $\hat{\mathcal{L}}_{\text{MMD}}$ yields $\mathbb{P}_{\hat{\bs{f}}_{\text{MMD}}(\bs{x}),\bs{\omega}}$ that we then sample from using each $(\bs{x}, \bs{\omega}) \in S$. We leave the derivation of a finite-sample performance guarantee that accounts for this dependence structure as a direction for future work.

\section{Experiments}\label{section:experiments}

This section evaluates the proposed CSG methods across four contextual two-stage stochastic programming (2SP) settings. Details of the computing environment and software are provided in Section~\ref{appendix:experimental_tooling}. Our primary goal is to assess how the proposed approaches compare to the baseline $K=1$ case (i.e., \citep{zharmagambetov2023landscape} or a gradient-based analogue of \citep{homem2024forecasting}). We also report results from solving each 2SP on a single deterministic scenario given by $\mathbb{E}[\bs{\omega}\mid \bs{x}]$. This helps quantify the gain from using additional samples $S$ even in the least-squares setting (e.g., the MMD approach with $K=1$).

Across all experiments, we benchmark methods by comparing the decisions they produce to the idealized solution obtained by solving the 2SP under the true conditional distribution $\mathbb{P}_{\bs{\omega}\mid \bs{x}}$ for a given context $\bs{x}$; i.e., the best decision one could hope to achieve. Finally, we emphasize that our methods are complementary to conditional-density-estimate-then-optimize approaches. For instance, in Section~\ref{subsection:cvarexperiment} we consider a trader who constructs $\mathbb{P}_{\bs{\omega}\mid \bs{x}}$ using $n$ residuals from fitted predictive models, as in \citep{deng2022predictive} and \citep{ban2019dynamic}. The practical challenge is then to compress this distribution so that solving the 2SP on $K \ll n$ scenarios still yields high-quality decisions.

\subsection{Newsvendor}

\textbf{Problem Setting:} The classic newsvendor problem illustrates the proposed approach. Each day, a newsvendor purchases $y$ newspapers at unit cost $c$. The vendor is budget-constrained, so they can only purchase at most $u$ papers. After purchasing, the vendor sells as many papers as possible for a unit price of $q$. When the day is done, the vendor can return the papers at a salvage price $r < c$. Suppose the vendor observes contextual information $\bs{x}$ before purchasing the papers. For example, $\bs{x}$ could represent the day of the week, the weather, and sales from previous days. The demand for the papers $\omega$ is unknown at the time of purchase. Thus, the newsvendor wishes to solve:
    \begin{equation} \label{eq:newsvendor}
    \min_{y} cy + \mathbb{E}_{\omega \sim \mathbb{P}_{\omega|\bs{x}}} \left[ Q(y, \omega)\right] \quad  \textit{s.t.} \quad  y \in [0, u],  
    \end{equation}
where $
        Q(y, \omega) = \min_{z \geq 0, w \geq 0} \ - q z - r w  \quad   \textit{s.t.} \quad z \leq \omega, \  z + w \leq y$. The decision variables $z$ and $w$ represent the quantities of papers sold and salvaged, respectively. We consider the following problem parameter setting $(c, q, r, u) = (1.0, 1.05, 0.1, 60.0)$. 

\vspace{\littlespace}
\noindent \textbf{Experimental Setup:}  A synthetic environment is used for evaluation purposes. $\bs{x}$ is set to be distributed according to a two-dimensional normal distribution. A randomly initialized neural network $\bs{f}_{\text{Random}} : \mathcal{X} \rightarrow \mathbb{R}_+^{200}$ maps $\bs{x}$ to an empirical demand distribution $\mathbb{P}_{\omega | \bs{x}}$ supported on $200$ points. The resulting joint distribution has an expectation of $15.1$ and a standard deviation of $0.25$ units.

The decision-maker has a dataset of $n = 500$ samples of $\bs{x}$ and $\omega$. The \eqref{eq:DCSG}, \eqref{eq:static_PCSG}, and Dynamic-PCSG approaches are trained on this sample. In figures and tables, we label these approaches as $\textit{MMD}$, $\textit{Static}$ and $\textit{Dynamic}$. All approaches let $\bs{f}_{\phi}$ be parameterized as a fully connected feedforward neural network with ReLU activation, mapping $\bs{x}$ to $\mathbb{R}_{+}^K$. We select the architecture by considering the \eqref{eq:DCSG} problem. A random search over a single training sample determines hyperparameters such as the number of hidden layers, number of units, and optimizer parameters. It optimizes $\hat{\mathcal{L}}_{\text{MMD}}$ on a fixed holdout of $20\%$ of the training samples. The loss network $E_{\psi}$'s architecture is selected similarly. First, $\bs{\hat{f}}_{\text{MMD}}$, generates the dataset $\{\bs{\zeta}_{1 \hdots K}^{(i)}, \omega^{(i)}, \ell^{(i)}_{\text{opt}} \}_{i=1}^n$, then a similar holdout procedure selects the architecture. The loss network $E_{\psi}$ is selected to be the same in both the static and dynamic approaches. We take this approach to selecting the network architectures in all the experiments considered in this work.

After training, we evaluate \eqref{eq:DCSG}, \eqref{eq:static_PCSG}, and Dynamic-PCSG on $n_{\text{eval}}=100$ out-of-sample contexts $\{\bs{x}^{(i)}_{\text{val}}\}_{i=1}^{n_{\text{eval}}}$ and their associated conditional distributions $\{\mathbb{P}_{\bs{\omega}\mid \bs{x}^{(i)}_{\text{val}}}\}_{i=1}^{n_{\text{eval}}}$. For each context $\bs{x}^{(i)}_{\text{val}}$ and learned task map $\bs{f}$, we form surrogate scenarios $\bs{\zeta}^{(i)}_{1\hdots K}=\bs{f}(\bs{x}^{(i)}_{\text{val}})$ and solve $(\bs{\zeta}\text{-SAA})$ to obtain a first-stage decision. We then evaluate this decision under the \emph{true} conditional distribution by computing the 2SP objective with the expected recourse taken with respect to $\mathbb{P}_{\bs{\omega}\mid \bs{x}^{(i)}_{\text{val}}}$. In all the experiments, we also verified that evaluating solutions using \eqref{eq:opt_search} (with its expectation taken under $\mathbb{P}_{\bs{\omega}\mid \bs{x}^{(i)}_{\text{val}}}$) yields identical results, and therefore we report the 2SP-based evaluation throughout. Let $v_{\text{MMD}}^{(i)}$, $v_{\text{Static}}^{(i)}$, and $v_{\text{Dynamic}}^{(i)}$ denote the resulting objective values for the $i$th validation context under \eqref{eq:DCSG}, \eqref{eq:static_PCSG}, and Dynamic-PCSG, respectively. As an oracle benchmark, we also compute the optimal 2SP solution under $\mathbb{P}_{\bs{\omega}\mid \bs{x}^{(i)}_{\text{val}}}$ and denote its objective value by $v_{\text{2SP}}^{(i)}$.

We consider two additional benchmarks. The first considers the expected value solution, determined by considering the newsvendor problem with a single scenario given by $\mathbb{E}[\omega | \bs{x}]$. In practice, the conditional mean is not known, so this is a benchmark that cannot be implemented. Considering the expected value solution allows us to ascertain whether the proposed methodologies can unlock the value of the stochastic solution (VSS) \citep{birge1982value}. The second benchmark is a quantile regression approach. It is well known that the optimal solution to the newsvendor is
$$
y^* = 
\begin{cases} 
0 & \text{if } \frac{q-c}{q-r} < F_{\mathbb{P}}(0), \\
u & \text{if } \frac{q-c}{q-r} > F_{\mathbb{P}}(u), \\
F_{\mathbb{P}}^{-1}\left( \frac{q-c}{q-r} \right) & \text{otherwise}. 
\end{cases}\quad , 
$$
where  $F_{\mathbb{P}}$ is the cumulative density according to $\mathbb{P}$ and $F^{-1}_{\mathbb{P}}(\alpha)$ is the $\alpha$ quantile of $F$. This suggests a quantile estimation approach to obtain a linear model \hbox{$F_{\mathbb{P}_{\omega|\bs{x}}}^{-1}\left( \frac{q-c}{q-r} \right) = \bs{\beta}^{\intercal} \bs{x}$}, from the training sample. The quantile estimation approach is well-studied in newsvendor problems. \citet{liu2022newsvendor} point out that quantile estimation is equivalent to determining the optimal linear (in $\bs{x}$) purchasing policy that minimizes the opportunity cost over the sample. The quantile regression (QR) is performed using the $\textsf{stats-models}$ package \citep{seabold2010statsmodels} that implements QR as in \citep{koenker2001quantile}. Let $v_{\mathbb{E}[\omega | \bs{x}]}^{(i)}$ and $v_{\text{QR}}^{(i)}$ denote the task losses for the solutions obtained via the expected value and QR approaches respectively, for the $i$th validation sample.

In this experiment $K$ varies in $\{1,2,5\}$, along with the MMD regularization parameter $\lambda \in \{10^{-1}, 10^{0}, 10^{1}\}$ to determine how the proposed approaches compare against the benchmarks. The experiment consists of $n_{\text{trials}} = 20$ trials. In each trial, a training set of $n$ joint observations and the evaluation set of $n_{\text{eval}}$ contexts and conditional distributions are sampled. Each approach is trained and evaluated. For each trial and validation instance, the optimality gap for each method is calculated as
$$
\textsf{Gap}_{m} = (v_{m} - v_{\text{2SP}})/|v_{\text{2SP}}|, \quad m \in \{\mathbb{E}[\omega | \bs{x}], \text{MMD}, \text{Static}, \text{Dynamic}, \text{QR} \}.
$$

\begin{figure}
    \FIGURE
    {\includegraphics[width=\textwidth]{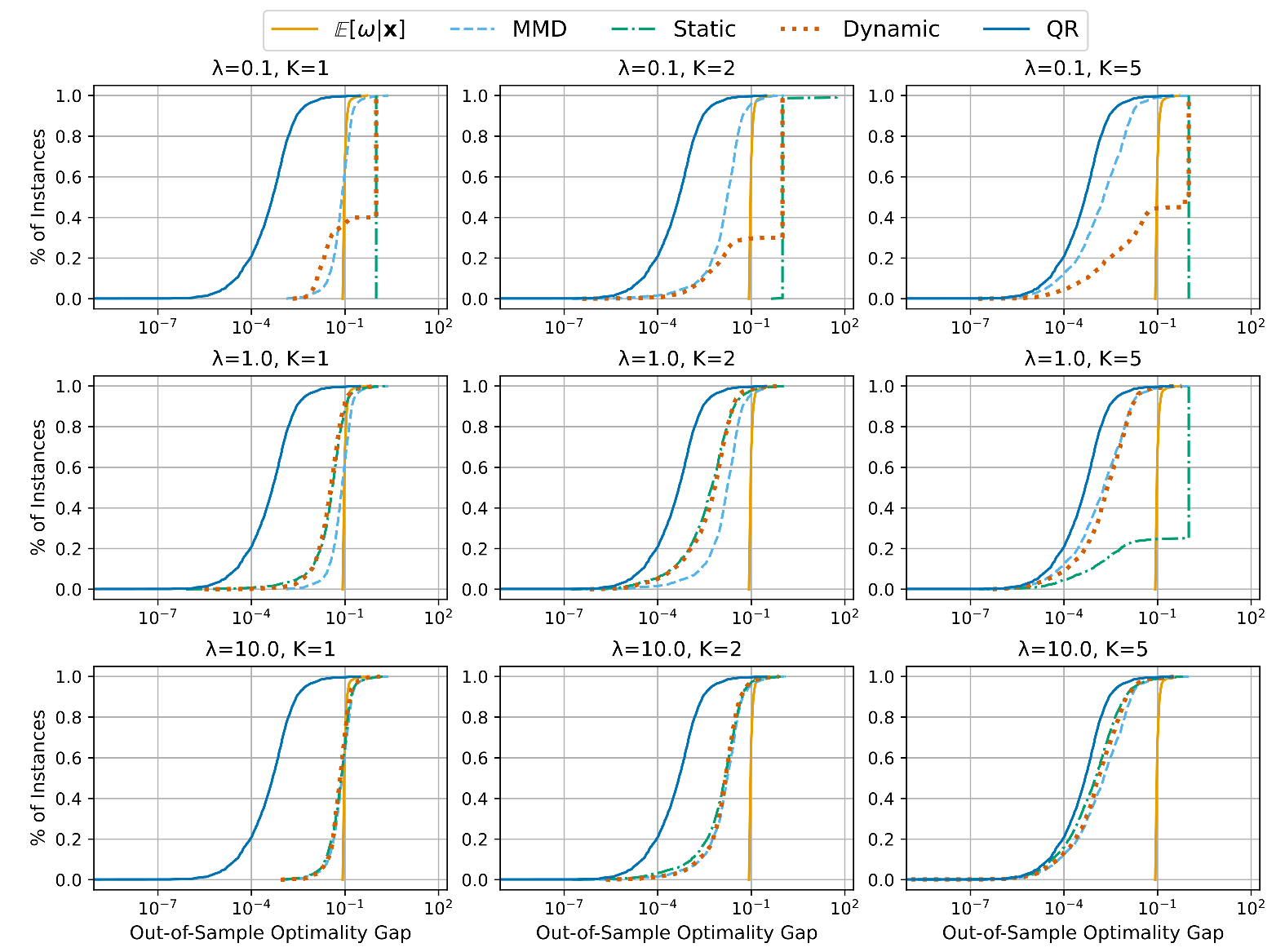}}
    {Out of sample optimality gaps for the contextual newsvendor problem.\label{fig:newsvendor}}
    {}
\end{figure}

\vspace{\littlespace}
\noindent \textbf{Results: } Figure \ref{fig:newsvendor} shows the cumulative distribution function (CDF) of the optimality gaps for each method by $K$ and $\lambda$ evaluated over the validation instances. If $(x,y)$ lies on the curve, the method achieves an optimality gap no more than $10^{x}$ on $100y\%$ of problems on the validation sets over all the trials. As expected, the quantile regression approach performs best as it explicitly uses the samples to predict the closed-form solution to the contextual newsvendor problem. The MMD approach captures more of the VSS as $K$ increases. However, the MMD approach has the downside that when $K = 1$, it reduces to conditional mean estimation via least-squares and, as such, is unlikely to offer any value over the expected value solution on average (since the expected value solution is the associated Bayes optimal predictor).

The task-based approaches yield poorly performing scenario mappings when \hbox{$\lambda = 0.1$}. In this setting, the dynamic approach slightly improves on the static approach. Furthermore, the MMD approach outperforms when the task-based approaches lack regularization. When $K \in \{1, 2\}$ and $\lambda = 1.0$, the task-based approaches outperform MMD significantly. Of the CSG approaches, the best-performing approaches when $K \in \{1, 2\}$ are the task-based approaches with $\lambda = 1.0$, and when $K = 5$, the best-performing approaches are the task-based approaches with $\lambda = 10$. However, when $K = 5$, the MMD approach is competitive with the task-based approaches. Furthermore, when $K > 1$ and $\lambda \geq 1.0$, the CSG approaches consistently exhibit value over the expected value solution. Lastly, task-based CSG outperforms the expected value solution when $K = 1$, $\lambda = 1.0$.

When $K=1$, one can also compare how well the proposed methods approximate the newsvendor-optimal scenario. In this case, the MMD approach reduces to a least-squares regression and aims to approximate $\mathbb{E}[\omega | \bs{x}]$. Conversely, it is observed that the task-based approaches, which aim to learn a scenario mapping that produces high-quality decisions downstream, yield scenarios that tend to be closer to the optimal quantile. This is in spite of them having no structural information about the newsvendor problem. Section \ref{appendix:newsvendor_scenarios} contains Figure~\ref{fig:newsvendor_scenarios} that compares the scenarios generated by each method.

Overall, task-based approaches offer more value than the MMD approach for smaller values of $K$. However, for larger $K$, the MMD approach can obtain higher-quality solutions. The importance of distributional regularization is also readily apparent, allowing PCSG methods to outperform DCSG significantly for small values of $K \leq 2$. Learning the surrogate scenario mapping for \hbox{($\bs{\zeta}$-SAA)} in a purely task-driven manner does not yield mappings that produce high-quality solutions. Although the QR approach exhibits the best performance, CSG methods still display desirable performance. Furthermore, CSG can be applied independently of the problem structure. The remaining experiments consider problem settings that are not endowed with a prescribed problem-driven approach to contextual optimization.

\subsection{Capacity Planning (CEP1)} 
% write out formulation
\textbf{Problem Setting:} Next, we test CSG via a contextual version of the CEP1 problem discussed in the introduction. In this application, the production manager wishes to decide the total hours of added capacity $\bs{y}_{\text{cap}} \in \mathbb{R}^{n_{\text{machines}}}_+$ and hours of operation $\bs{y}_{\text{op}} \in \mathbb{R}^{n_{\text{machines}}}_+$ to meet uncertain demand for the $N$ parts $\bs{\omega} \in \mathbb{R}_+^N$. The detailed formulation of the CEP1 problem is relegated to Section~\ref{appendix:cep1_formulation}. 

\vspace{\littlespace}
\noindent \textbf{Experimental Setup:}
Similar to the newsvendor, we consider a randomly initialized neural network to generate the conditional distribution $\mathbb{P}_{\bs{\omega} | \bs{x}}$ for context $\bs{x}$. The procedure defining the underlying $\mathbb{P}_{\bs{\omega} | \bs{x}}$ is described in Section~\ref{appendix:cep1_data_generation}. The experimental setup is nearly identical to that of the newsvendor. Table \ref{tab:CEP1_parameters} contains the CEP1 experiment parameters. 
% \begin{table}[h]
%     \TABLE
%     {Summary of CEP1 Experimental Parameters \label{tab:CEP1_parameters}}
%     {\begin{tabular}{ll}
%         \toprule
%         \textbf{Parameter} & \textbf{Value} \\
%         \midrule
%         Number of samples ($n$) & 500 \\
%         Number of out-of-sample instances ($n_{\text{eval}}$) & 200 \\
%         Number of surrogate scenarios ($K$) & $\{1, 2, 5\}$ \\
%         MMD regularization ($\lambda$, scaled by $\hat{\mathcal{L}}_{\text{MMD}}$) & $\hat{\mathcal{L}}_{\text{MMD}} * \tilde{\lambda}$ where $\tilde{\lambda} \in \{10^{-3+i}\}_{i=1}^4$ \\
%         Number of sampling repetitions ($n_{\text{trials}}$) & 20 \\
%         Number of Loss-Net updates in dynamic PCSG ($T$) & 4 \\
%         \bottomrule
%         \end{tabular}}{}
%     \end{table}

\begin{table}[h]
    \TABLE
    {Summary of CEP1 Experimental Parameters \label{tab:CEP1_parameters}}
    {\begin{tabular}{ll}
        \toprule
        \textbf{Parameter} & \textbf{Value} \\
        \midrule
        Number of samples ($n$) & 500 \\
        Number of out-of-sample instances ($n_{\text{eval}}$) & 200 \\
        Number of surrogate scenarios ($K$) & $\{1, 2, 5\}$ \\
        MMD regularization ($\lambda$, scaled by $\hat{\mathcal{L}}_{\text{MMD}}$) & $\hat{\mathcal{L}}_{\text{MMD}} * \tilde{\lambda}$ where $\tilde{\lambda} \in \{10^{-3+i}\}_{i=1}^4$ \\
        Number of sampling repetitions ($n_{\text{trials}}$) & 20 \\
        Number of Loss-Net updates in dynamic PCSG ($T$) & 4 \\
        \bottomrule
        \end{tabular}}{}
    \end{table}

The resulting 2SP is a two-stage LP; thus, there are no guarantees for uniqueness. Furthermore, the Task-Net $\bs{f}_{\phi}$ is parameterized as a fully connected feedforward neural network with ReLU activation, mapping $\bs{x}$ to $\mathbb{R}_{+}^{m \times K}$.  
%     \item the problem parameters in table 
%     \item describe the architecture 
%     \item CEP1 is LP, use $\ell_{\text{opt}}$
% \end{itemize}

\vspace{\littlespace}
\noindent \textbf{Results:} The CDF of the optimality gap for each method, by $K$ and $\lambda$ is relegated to Section~\ref{appendix:cep1_performance} for space considerations. Table \ref{tab:CEP1_gaps} shows the median optimality gaps over all the validation instances and trials. We observe in this setting that the VSS is not large since the expected value solution attains a median optimality gap of $1.69\%$. Unlike the newsvendor, the MMD approach does not yield solutions that are closer to optimal as $K$ increases, yielding some evidence of an overfitting effect for large $K$. Although the MMD approach still underperforms the expected value solutions when $K = 1$. In settings where $\lambda \geq 0.1$, the dynamic approach consistently improves on the static approach, highlighting the performance gains obtained by iteratively refining the Loss-Net via the dynamic approach. When $(K, \tilde{\lambda}) = (1, 0.01)$, the PCSG methods attain smaller optimality gaps than the other methods and can unlock the VSS with a single scenario. When $K \geq 2$, the PCSG methods perform similarly to MMD.

\begin{table}[h]
    \centering
        \begin{minipage}{0.49\textwidth}
        \centering
        \TABLE
        {(CEP1) Median optimality gaps on validation instances for various methods \label{tab:CEP1_gaps}}
        {
        \begin{tabular}{cccccc}
        \toprule
         & $K=1$ & $K=2$ & $K=5$ \\
        \midrule
        $\mathbb{E}[\bs{\omega}| \bs{x}]$ & 1.69\% & 1.69\% & 1.69\% \\
        MMD & 2.36\% & 0.43\% & 0.76\% \\
        \midrule
        \multicolumn{1}{c}{$\tilde{\lambda}$} & \multicolumn{3}{c}{\textbf{Static}} \\
        \midrule
        0.01 & 1.50\% & 5.04\% & 6.31\% \\
        0.1  & 2.07\% & 1.23\% & 1.20\% \\
        1    & 2.31\% & 0.44\% & 0.65\% \\
        10   & 2.36\% & 0.39\% & 0.65\% \\
        \midrule
        \multicolumn{1}{c}{$\tilde{\lambda}$} &  \multicolumn{3}{c}{\textbf{Dynamic}} \\
        \midrule
        0.01 & 1.18\% & 0.57\% & 0.78\% \\
        0.1  & 1.89\% & 0.70\% & 0.99\% \\
        1    & 2.23\% & 0.44\% & 0.61\% \\
        10   & 2.31\% & 0.39\% & 0.49\% \\
        \bottomrule
        \end{tabular}}{}
    \end{minipage}
     \hfill
        \begin{minipage}{0.49\textwidth}
        %\captionsetup{width=.75\textwidth}
        \centering
        \TABLE
        {(CEP1) Fraction of \\ instances with the lowest \\ out-of-sample 2SP cost. \label{tab:CEP1_competition}}
        {
        \begin{tabular}{cccccc}
        \toprule
         & $K=1$ & $K=2$ & $K=5$ \\
        \midrule
        $\mathbb{E}[\bs{\omega}| \bs{x}]$ & 11.6\% & 0.15\% & 0.10\% \\
        MMD & 0.0\% & 15.1\% & 11.7\% \\
        \midrule
        \multicolumn{1}{c}{$\tilde{\lambda}$} & \multicolumn{3}{c}{\textbf{Static}} \\
        \midrule
        0.01 & 29.8\% & 5.40\%  & 2.10\%  \\
        0.1  & 0.18\%  & 4.85\%  & 9.30\%  \\
        1    & 0.00\%  & 10.0\% & 13.7\% \\
        10   & 0.05\%  & 12.6\% & 14.4\% \\
        \midrule
        \multicolumn{1}{c}{$\tilde{\lambda}$} &  \multicolumn{3}{c}{\textbf{Dynamic}} \\
        \midrule
        0.01 & 55.7\% & 18.9\% & 15.5\% \\
        0.1  & 2.10\%  & 9.15\%  & 10.9\% \\
        1    & 0.40\%  & 10.1\% & 9.03\%  \\
        10   & 0.27\%  & 13.8\% & 13.3\% \\ 
        \bottomrule
        \end{tabular}}{}
    \end{minipage}
    \hfill
\end{table}

Table \ref{tab:CEP1_competition} compares the performance of the different methods across $K$ by showing the percentage of instances in which the particular method obtains the lowest out-of-sample 2SP cost (the columns sum to $100\%$). When $K=1$, the MMD (regression) approach never achieves the lowest 2SP objective. Knowing the true conditional expectation yields the best solution in only 11.6\% of instances. In contrast, the static and dynamic methods outperform when $\tilde{\lambda} = 0.01$, achieving the lowest cost in $29.8\%$ and $55.7\%$ of instances, respectively. When $K = 2$ $(5)$, the MMD approach becomes more competitive, obtaining the lowest 2SP cost in $15.1\%$ ($11.7\%$) instances. Still, when $K \in \{2, 5\}$, the dynamic approach with $\tilde{\lambda}$ obtains the best objective most often. Like the contextual newsvendor, task-based approaches offer more value than the MMD approach for smaller values of $K$. Although the MMD approach can obtain higher-quality solutions for larger $K$, with appropriate regularization, task-driven approaches tend to obtain the best results more often.

\subsection{Portfolio Optimization} \label{subsection:cvarexperiment}

\textbf{Problem Setting:} Next, we test CSG in a contextual portfolio optimization setting. In this application the trader wishes to select a portfolio $\bs{y} \in \mathbb{R}_+^{n_{\text{assets}}}$ to balance risk and expectation of the portfolio return $\bs{\omega}^{\intercal} \bs{y}$ where $\bs{\omega} \in \mathbb{R}^{n_{\text{assets}}}$ denotes the asset-returns. For execution reasons, the portfolio must consist of some subset $K_{\text{assets}} < n_{\text{assets}}$. The trader wishes to use contextual information to make better decisions quickly. In this setting, the trader takes $\bs{x} \in \mathbb{R}^{(n_{\text{assets}} + n_{\text{encoding}}) \times W}$ to be a panel of $W$ past return observations for each of the $n_{\text{assets}}$ assets, along with a representation of the 5-minute period within the 78 period trading session. The Conditional Value-at-Risk (CVaR) is the risk measure of choice. The reader is referred to \citet{krokhmal2002portfolio} for more details regarding portfolio selection based on CVaR. The detailed formulation of the CVaR problem is relegated to Section~\ref{appendix:cvar_formulation}.

\vspace{\littlespace}
\noindent \textbf{Experimental Setup:} We consider a data-driven approach to set up the contextual environment. The trader assumes that the returns satisfy the following relationship: $
\bs{\omega} = \Phi_{\bar{\bs{r}}}(\boldsymbol{x}) + \Phi_{\bs{\sigma}} (\boldsymbol{x}) \odot \boldsymbol{\epsilon}$, 
where $\Phi_{\bar{\bs{r}}} : \mathcal{X} \rightarrow \mathbb{R}^{n_{\text{assets}}}$ and $\Phi_{\bs{\sigma}} : \mathcal{X} \rightarrow \mathbb{R}_+^{n_{\text{assets}}}$ are models given to the trader by a statistical modeler to estimate the conditional mean $\mathbb{E}[\bs{\omega} |\bs{x}]$, and conditional standard deviation $\sqrt{\text{Var}[\bs{\omega} | \bs{x}]}$ along with $\boldsymbol{\epsilon}$ denoting the errors which obey a uniform distribution supported on $n^* = 1559$ points $\{\bs{\epsilon}^{(i)}\}_{i=1}^{n^*}$, derived from a dataset of returns. The procedure defining  $\mathbb{P}_{\bs{\omega} |\bs{x}}$ is fully described in Section~\ref{appendix:cvar_data_generation}.

For a given contextual observation, the trader faces an instance of CVaR optimization defined on $n^*$ scenarios. We aim to explore whether the CSG frameworks proposed in this work can enable the trader to directly generate a set of $K$ scenarios from the context, which can then be solved more quickly than solving the problem on all $n^*$ scenarios. We consider a training set and a validation set to evaluate our framework. The trained task-networks are evaluated on the validation set using $\mathbb{P}_{\bs{\omega} | \bs{x}}$ defined by $\Phi_{\bar{\bs{r}}}$, $\Phi_{\bs{\sigma}}$ and $\{\bs{\epsilon}^{(i)}\}_{i=1}^{n^*}$. This approach is sensible since our goal is to help the trader solve their CVaR optimization with $\mathbb{P}_{\bs{\omega} | \bs{x}}$ defined as described above. \hbox{Table \ref{tab:CVaR_parameters}} highlights the experiment parameters.

\begin{table}[h]
    \TABLE
    {Summary of CVaR Experimental Parameters \label{tab:CVaR_parameters}}
    {\begin{tabular}{ll}
        \toprule
        \textbf{Parameter} & \textbf{Value} \\
        \midrule
        Number of samples ($n$) & $\lfloor{0.8n^* +\frac{1}{2}}\rfloor = 1247$\\
        Number of out-of-sample instances ($n_{\text{eval}}$) & $n^* - \lfloor{0.8n^* +\frac{1}{2}}\rfloor = 312$ \\
        Number of surrogate scenarios ($K$) & $\{1, 3, 5, 10, 20, 40\}$ \\
        MMD regularization ($\lambda$) & $\{80*2^{i}\}_{i=0}^7$ \\
        %Number of Loss-Net updates in dynamic PCSG ($T$) & 4 \\
        \bottomrule
        \end{tabular}}{}
    \end{table}

The resulting 2SP is mixed-binary in the first stage and is linear in the second stage; thus, there are no uniqueness guarantees. Furthermore, The Task-Net $\bs{f}_{\phi}$ is parameterized by an LSTM network that outputs a \hbox{$n_{\text{assets}} \times K$} matrix.

\vspace{\littlespace}
\noindent \textbf{Results:} The CDF of the optimality gap for each method, by $K$ and $\lambda$ is relegated to Section~\ref{appendix:cvar_performance} for space considerations. Table \ref{tab:cvar_gaps} shows the median optimality gaps over the validation instances. Unlike CEP1, the VSS is typically large, with the expected value solution exhibiting a median gap of $241\%$. Like the other experiments, the MMD approach produces solutions closer to optimal as $K$ increases while underperforming the expected value solutions when $K=1$.

\begin{table}[h]
    \centering 
        \TABLE
        {(CVaR) Median optimality gaps on validation instances for various methods. \label{tab:cvar_gaps}}
        {
        \begin{tabular}{ccccccc}
        \toprule
         & $K=1$ & $K=3$ & $K=5$ & $K=10$ & $K=20$ & $K=40$ \\
        \midrule
        $\mathbb{E}[\bs{\omega}| \bs{x}]$ & 241\%      & 241\%      & 241\%      & 241\%      & 241\%      & 241\%      \\
        MMD & 483\%      & 166\%      & 122\%      & 55.1\%     & 18.7\%     & 10.6\%     \\
        \midrule
        \multicolumn{1}{c}{$\lambda$} & \multicolumn{6}{c}{\textbf{Static}} \\
        \midrule
        80    & 1.03e7\%   & 486\%      & 96.9\%     & 79.0\%     & 38.2\%     & 118\%      \\
        160   & 7.41e3\%   & 120\%      & 100\%      & 84.9\%     & 92.2\%     & 25.7\%     \\
        320   & 2.60e7\%   & 106\%      & 101\%      & 45.7\%     & 28.4\%     & 11.9\%     \\
        640   & 424\%      & 98.4\%     & 96.9\%     & 93.8\%     & 29.7\%     & 10.9\%     \\
        1280  & 456\%      & 115\%      & 77.1\%     & 88.5\%     & 15.8\%     & 8.67\%     \\
        2560  & 474\%      & 182\%      & 70.8\%     & 32.5\%     & 27.6\%     & 12.4\%     \\
        5120  & 486\%      & 152\%      & 94.1\%     & 45.0\%     & 19.5\%     & 9.27\%     \\
        10240 & 485\%      & 164\%      & 76.2\%     & 33.6\%     & 18.0\%     & 17.5\%     \\
        \midrule
        \multicolumn{1}{c}{$\lambda$} & \multicolumn{6}{c}{\textbf{Dynamic}} \\
        \midrule
        80    & 1.02e3\%   & 296\%      & 225\%      & 94.8\%     & 212\%      & 62.8\%     \\
        160   & 1.00e3\%   & 310\%      & 135\%      & 93.5\%     & 221\%      & 54.1\%     \\
        320   & 365\%      & 172\%      & 98.6\%     & 60.8\%     & 82.1\%     & 20.2\%     \\
        640   & 736\%      & 154\%      & 94.2\%     & 28.6\%     & 38.8\%     & 11.6\%     \\
        1280  & 433\%      & 108\%      & 91.2\%     & 51.8\%     & 24.2\%     & 18.6\%     \\
        2560  & 401\%      & 89.5\%     & 65.3\%     & 39.6\%     & 39.8\%     & 7.55\%     \\
        5120  & 439\%      & 160\%      & 66.2\%     & 28.1\%     & 14.4\%     & 10.3\%     \\
        10240 & 458\%      & 148\%      & 75.5\%     & 29.9\%     & 29.7\%     & 5.07\%     \\
        \bottomrule
        \end{tabular}}{}
    \end{table}
    \begin{table}[h]
         \centering 

        \TABLE
        {(CVaR) Fraction of instances with the lowest out-of-sample 2SP cost. \label{tab:cvar_competition}}
        {
        \begin{tabular}{ccccccc}
        \toprule
         & $K=1$ & $K=3$ & $K=5$ & $K=10$ & $K=20$ & $K=40$ \\
        \midrule
        $\mathbb{E}[\bs{\omega}| \bs{x}]$ & 74.4\%   & 0.00\%    & 0.00\%    & 0.00\%    & 0.00\%    & 0.00\%    \\
        MMD & 0.32\%   & 0.32\%    & 0.00\%    & 0.00\%    & 30.8\%    & 0.32\%    \\
        \midrule
        \multicolumn{1}{c}{$\lambda$} & \multicolumn{6}{c}{\textbf{Static}} \\
        \midrule
        80    & 0.00\%   & 0.00\%    & 12.8\%    & 6.09\%    & 0.00\%    & 0.00\%    \\
        160   & 0.00\%   & 8.65\%    & 8.65\%    & 4.81\%    & 0.00\%    & 0.00\%    \\
        320   & 0.00\%   & 13.1\%    & 4.81\%    & 0.00\%    & 0.00\%    & 4.81\%    \\
        640   & 2.24\%   & 14.7\%    & 26.0\%    & 0.00\%    & 10.3\%    & 29.5\%    \\
        1280  & 4.17\%   & 1.28\%    & 0.00\%    & 0.64\%    & 13.8\%    & 1.28\%    \\
        2560  & 0.64\%   & 0.00\%    & 2.56\%    & 0.32\%    & 6.73\%    & 0.00\%    \\
        5120  & 0.32\%   & 0.64\%    & 0.00\%    & 1.92\%    & 0.00\%    & 0.00\%    \\
        10240 & 0.32\%   & 0.96\%    & 0.00\%    & 20.2\%    & 0.32\%    & 0.00\%    \\
        \midrule
        \multicolumn{1}{c}{$\lambda$} & \multicolumn{6}{c}{\textbf{Dynamic}} \\
        \midrule
        80    & 0.00\%   & 0.00\%    & 0.00\%    & 0.00\%    & 0.00\%    & 0.00\%    \\
        160   & 0.00\%   & 0.00\%    & 0.00\%    & 0.00\%    & 0.00\%    & 0.00\%    \\
        320   & 10.9\%   & 0.00\%    & 3.85\%    & 0.32\%    & 0.00\%    & 5.13\%    \\
        640   & 0.00\%   & 6.73\%    & 12.5\%    & 19.9\%    & 6.73\%    & 3.21\%    \\
        1280  & 2.88\%   & 6.41\%    & 1.28\%    & 0.00\%    & 0.00\%    & 0.00\%    \\
        2560  & 2.56\%   & 44.2\%    & 13.5\%    & 0.00\%    & 0.00\%    & 3.53\%    \\
        5120  & 1.28\%   & 0.00\%    & 13.1\%    & 16.7\%    & 31.4\%    & 0.64\%    \\
        10240 & 0.00\%   & 2.88\%    & 0.96\%    & 29.2\%    & 0.00\%    & 51.6\%    \\
        \bottomrule
        \end{tabular}}{}
\end{table}

When $K = 1$, the CSG approaches do not outperform the expected value solution, although the static and dynamic methods tend to outperform the MMD approach for sufficiently large $\lambda$. The static and dynamic approaches also attain smaller optimality gaps than the MMD approach for $K \geq 3$ and sufficiently large $\lambda$. For example, when $K = 3$ and $\lambda \geq 160$, the static approach outperforms the MMD approach. Lastly, for any fixed $K$, the dynamic approach attains a smaller optimality gap than the static approach for some suitable value of $\lambda$. For instance, when $K = 10$, the dynamic approach attains a median gap of $28.1\%$ when $\lambda = 5120$, whereas the static approach's best gap is $32.5\%$ when $\lambda = 2560$. When $K = 40$, the static and dynamic approaches have median gaps less than $10.0\%$, whereas the MMD approach does not. Notably, the dynamic approach with $(K, \lambda)= (40, 10240)$ can achieve a median gap of $5.07\%$.

Table \ref{tab:cvar_competition} compares the performance of the different methods across $K$ by showing the percentage of instances in which the particular method obtains the lowest out-of-sample 2SP cost (the columns sum to $100\%$). When $K=1$, the expected value solution dominates; however, when \hbox{$K > 1$}, the expected value solution never obtains the lowest out-of-sample cost. Furthermore, the problem-driven approaches always outperform when \hbox{$K > 1$}, since for every $K$, there is some $\lambda$ such that the problem-driven approach wins most often. Unlike previous experiments, the task-based approaches offer more value than the MMD approach across all values of tested $K$. 

\subsection{Multidimensional Newsvendor with Substitution} 
% write out formulation, mip, non-linear

\noindent \textbf{Problem Setting:} Next, we test CSG in a contextual setting where the decision-maker is tasked with solving a multidimensional newsvendor problem. We consider the version of the problem introduced by \citet{narum2024problem} and provide a brief overview. The Multidimensional Newsvendor with Substitution (MNV) is a production planning problem with a structure similar to \eqref{eq:newsvendor}, where production decisions regarding $N$ products are made before uncertain demand $\bs{\omega} \in \mathbb{R}^N$ is realized. The primary distinction is that MNV considers multiple products and allows products to be substituted with each other. The original application by \citet{vaagen2011modelling} considers fashion apparel where retailers offer substitutable products. The detailed formulation of the MNV problem is relegated to Section~\ref{appendix:MNV_formulation}. 

% \item where this is this, that is that
% \end{itemize}

\vspace{\littlespace}
\noindent \textbf{Experimental Setup:}  Similar to newsvendor and CEP1, we consider randomly initialized neural networks to generate the conditional distribution $\mathbb{P}_{\bs{\omega} | \bs{x}}$ for context $\bs{x}$. The procedure defining the underlying $\mathbb{P}_{\bs{\omega} | \bs{x}}$ is described in Section~\ref{appendix:MNV_data}. The process defining $\mathbb{P}_{\bs{\omega} | \bs{x}}$ aims to capture key aspects of the product-demand distribution as pointed out by \citet{vaagen2011modelling}. Table \ref{tab:MNV_parameters} highlights the experiment parameters.  
\begin{table}[h]
    \TABLE
    {Summary of MNV Experimental Parameters  \label{tab:MNV_parameters}}
    {
    \begin{tabular}{ll}
    \toprule
    \textbf{Parameter} & \textbf{Value} \\
    \midrule
    Number of samples ($n$) & $300$\\
    Number of out-of-sample instances ($n_{\text{eval}}$) & $100$ \\
    Number of surrogate scenarios ($K$) & $\{1, 3, 5, 10\}$ \\
    MMD regularization ($\lambda$) & $\{5*4^{i}\}_{i=0}^6$ \\
    \bottomrule
    \end{tabular}}{}
\end{table}

The resulting MNV problem is a mixed-binary program (in both stages); thus, there are no uniqueness guarantees. Furthermore, the Task-Net $\bs{f}_{\phi}$ is parameterized as a fully connected feedforward neural network with ReLU activation, mapping $\bs{x}$ to $\mathbb{R}_{+}^{m \times K}$.

\vspace{\littlespace}
\noindent \textbf{Results:} The CDF of the optimality gap for each method, by $K$ and $\lambda$ is relegated to Section~\ref{appendix:MNV_performance} for space considerations. Table \ref{tab:mnv_gaps} shows the median optimality gaps over the validation instances. In the MNV setting, the VSS is small, with the expected value solution exhibiting a median gap of $2.70\%$. The MMD approach produces solutions closer to optimal as $K$ is increased. However, the MMD approach does not outperform for $K > 1$ and outperforms the expected value solution when $K = 10$.

When $K = 1$, the CSG approaches underperform the expected value solution, although the static approach outperforms the MMD approach for sufficiently large $\lambda \geq 80$. When $K = 3$, the static and dynamic approaches underperform the MMD approach and perform similarly when $\lambda = 20480$. In the case $K = 5$, the static approach ($\lambda \geq 320$) attains smaller median gaps than the expected value and MMD approach, whereas the dynamic approach does not. In particular, the static method when $(K, \lambda) \in \{(5, 5120), (5, 20480)\}$ attains smaller optimality gaps more frequently than the expected value and MMD approaches (see Section~\ref{appendix:MNV_performance}). Lastly, when $K = 10$, the static and dynamic approaches outperform the expected value approach, but only the static approach obtains similar performance 

\begin{table}[h]
    \begin{minipage}{0.49\linewidth}
        \TABLE
        {(MNV) Median optimality gaps on validation instances by method. \label{tab:mnv_gaps}}
        {
        \begin{tabular}{ccccc}
        \toprule
         & $K=1$ & $K=3$ & $K=5$ & $K=10$ \\
        \midrule
        $\mathbb{E}[\bs{\omega}| \bs{x}]$ & 2.70\%    & 2.70\%    & 2.70\%    & 2.70\%    \\
        MMD & 3.81\%    & 2.84\%    & 2.82\%    & 1.66\%    \\
        \midrule
        \multicolumn{1}{c}{$\lambda$} & \multicolumn{4}{c}{\textbf{Static}} \\
        \midrule
        5     & 102\%     & 37.5\%    & 18.4\%    & 13.4\%    \\
        20    & 4.61\%    & 9.78\%    & 5.40\%    & 2.60\%    \\
        80    & 2.80\%    & 4.86\%    & 3.44\%    & 1.93\%    \\
        320   & 3.32\%    & 3.17\%    & 2.61\%    & 1.75\%    \\
        1280  & 3.67\%    & 3.15\%    & 2.65\%    & 1.81\%    \\
        5120  & 3.73\%    & 2.89\%    & 2.07\%    & 1.73\%    \\
        20480 & 3.72\%    & 2.80\%    & 2.19\%    & 1.90\%    \\
        \midrule
        \multicolumn{1}{c}{$\lambda$} & \multicolumn{4}{c}{\textbf{Dynamic}} \\
        \midrule
        5     & 74.7\%    & 55.9\%    & 29.8\%    & 37.8\%    \\
        20    & 8.69\%    & 39.3\%    & 23.2\%    & 7.65\%    \\
        80    & 3.26\%    & 10.9\%    & 5.73\%    & 2.73\%    \\
        320   & 3.70\%    & 3.44\%    & 3.12\%    & 2.10\%    \\
        1280  & 4.18\%    & 3.76\%    & 2.79\%    & 2.16\%    \\
        5120  & 4.26\%    & 3.37\%    & 2.89\%    & 2.17\%    \\
        20480 & 4.26\%    & 2.83\%    & 2.85\%    & 2.01\%    \\
        \bottomrule
        \end{tabular}}{}
    \end{minipage}
    \hfill
    \begin{minipage}{0.49\linewidth}
        \TABLE
        {(MNV) Fraction of instances with the lowest out-of-sample 2SP cost. \label{tab:mnv_competition}}
        {
        \begin{tabular}{ccccc}
        \toprule
         & $K=1$ & $K=3$ & $K=5$ & $K=10$ \\
        \midrule
        $\mathbb{E}[\bs{\omega}| \bs{x}]$ & 20.8\%   & 13.8\%   & 9.33\%   & 10.0\%    \\
        MMD & 2.17\%   & 11.7\%   & 7.00\%   & 9.50\%    \\
        \midrule
        \multicolumn{1}{c}{$\lambda$} & \multicolumn{4}{c}{\textbf{Static}} \\
        \midrule
        5      & 2.17\%   & 3.67\%   & 0.67\%   & 1.50\%    \\
        20     & 12.5\%   & 0.50\%   & 8.00\%   & 9.50\%    \\
        80     & 11.5\%   & 5.50\%   & 6.50\%   & 7.67\%    \\
        320    & 2.00\%   & 7.83\%   & 6.17\%   & 9.50\%    \\
        1280   & 3.33\%   & 8.17\%   & 9.17\%   & 7.00\%    \\
        5120   & 2.00\%   & 6.83\%   & 10.0\%   & 5.50\%    \\
        20480  & 2.17\%   & 10.7\%   & 11.8\%   & 7.33\%    \\
        \midrule
        \multicolumn{1}{c}{$\lambda$} & \multicolumn{4}{c}{\textbf{Dynamic}} \\
        \midrule
        5      & 0.00\%   & 0.00\%   & 0.00\%   & 0.83\%    \\
        20     & 13.0\%   & 0.00\%   & 0.17\%   & 2.00\%    \\
        80     & 14.3\%   & 0.17\%   & 2.83\%   & 5.83\%    \\
        320    & 4.00\%   & 7.50\%   & 8.17\%   & 7.17\%    \\
        1280   & 4.50\%   & 7.50\%   & 7.33\%   & 5.67\%    \\
        5120   & 3.00\%   & 7.33\%   & 7.17\%   & 5.17\%    \\
        20480  & 2.50\%   & 8.83\%   & 5.17\%   & 5.83\%    \\
        \bottomrule
        \end{tabular}}{}
    \end{minipage}
\end{table}

Table \ref{tab:mnv_competition} compares the performance of the different methods across $K$ by showing the percentage of instances in which the particular method obtains the lowest out-of-sample 2SP cost (the columns sum to $100\%$). When $K=1$, the expected value solution dominates; however, the static and dynamic approaches with $\lambda \in \{20, 80\}$ are competitive. 
When $K > 1$, the expected value solution's relative performance deteriorates by roughly a factor of 2. Furthermore, in the cases $K > 1$, the wins are split more evenly among the MMD, static, and dynamic approaches, where the static and dynamic approaches are heavily regularized. 

Surprisingly, the dynamic approach fails to outperform the static approach (as one would expect that iteratively refining $E_{\psi}$ would yield better results). In the MNV example, $E_{\psi}$ was able to construct an accurate approximation on the dataset generated by the learned MMD Task-Net: $\bs{\zeta}_{1 \hdots K}^{(i)} = \bs{\hat{f}}_{\text{MMD}}(\bs{x}^{(i)})$, $\ell_{\text{opt}}^{(i)} = \ell_{\text{opt}}(\bs{\zeta}_{1 \hdots K}^{(i)}, \bs{\omega}^{(i)})$. However, the architecture for $E_{\psi}$, selected by random search with $20\%$ holdout, did not generalize to the dataset generated in subsequent iterations of the dynamic approach. When one faces this situation, a sensible approach is to consider task-nets $\bs{f}$ that produce $\tilde{\bs{\zeta}}_{1 \hdots K}^{(i)} = \bs{f}(\bs{x}^{(i)})$ that are close to $\bs{\zeta}_{1 \hdots K}^{(i)}$ with respect to MMD distance. Setting $\lambda$ to a large value aims to achieve this, and indeed, we observe in Tables \ref{tab:mnv_gaps} and \ref{tab:mnv_competition} that only heavily regularized PCSG approaches achieve outperformance relative to the MMD approach. 

\subsection{Interpretation across experiments}
\label{subsec:problem_driven_vs_agnostic}

The numerical results indicate that the relative performance of problem-driven versus problem-agnostic training is problem-dependent. We interpret this through (i) the value of the stochastic solution (VSS), which upper bounds the attainable benefit from modelling uncertainty, and (ii) the dominant error terms in Theorem~\ref{thm:appxlossguarantee}, which quantify the implicit costs of learning and optimizing a task surrogate.

VSS measures how much is gained by solving the 2SP under the full conditional distribution $\mathbb{P}_{\bs{\omega}\mid\bs{x}}$ rather than using the single  scenario $\mathbb{E}[\bs{\omega}\mid\bs{x}]$. Thus, when VSS is small, there is limited room to improve upon the expected value policy, so even modest learning errors can erase the potential gains of problem-driven methods.

Theorem~\ref{thm:appxlossguarantee} explains how such errors arise when optimizing a learned surrogate $\hat{E}$ (e.g., Loss-Net) in place of the true task loss $\ell_{\text{opt}}$. While problem-driven training can reduce the empirical surrogate loss, it may simultaneously incur (a) approximation error if $\hat{E}$ does not accurately represent $\ell_{\text{opt}}$, and (b) distribution shift if optimizing $\hat{E}$ moves $\mathbb{P}_{\bs{f}(\bs{x}),\bs{\omega}}$ away from the distribution on which $\hat{E}$ is reliable (in the case of \eqref{eq:static_PCSG}, approximated by $\mathbb{P}_{\bs{\hat{f}}_{\text{MMD}}(\bs{x}),\bs{\omega}}$ as a proxy for $\mathbb{P}^*_{\bs{\zeta}_{1:K},\bs{\omega}}$). When $\ell_{\text{opt}}$ is difficult to approximate, and the induced shift is large, i.e., $\bs{f}(\bs{x})$ places mass in regions of $\Omega^K$ poorly covered during surrogate training, then these costs can outweigh the VSS.

This perspective is consistent with our experiments. In the newsvendor and CVaR portfolio settings, VSS is large and problem-driven methods can substantially improve upon MMD. In CEP1 and MNV, the gains are less pronounced. Moreover, the difference between CEP1 and MNV aligns with the surrogate approximation term: for MNV, we observed surrogate errors on the order of the VSS, whereas for CEP1, the approximation errors were smaller.

Overall, Theorem~\ref{thm:appxlossguarantee} suggests that problem-driven training is most beneficial when (a) VSS is non-negligible and (b) $\hat{E}$ approximates $\ell_{\text{opt}}$ well on the relevant region of $\Omega^K$ while inducing only a small distribution shift. When VSS is small (as in MNV) or approximation/shift errors are large, the implicit costs of problem-driven training dominate, and a distribution-focused objective such as MMD can yield better and more stable performance.

\subsection{Timing Analysis}

This section explores the time investment trade-offs presented by the proposed methodologies. A time investment exists to train the Task-Net and Loss-Net used in the CSG approaches. However, once trained, the Task-Net can quickly generate solutions to 2SP by computing a forward pass through the Task-Net and solving ($\bs{\zeta}$-SAA) on the resulting $K$ surrogate scenarios. The following processes detail the different training and evaluation stages:
\begin{itemize} % Adjust indentation and spacing between items
    \item \textbf{MMD Training:} Time taken to train the Task-Net via DCSG.
    \item \textbf{Surrogate Solution Calculation:} Time taken to pass through the Task-Net and compute an optimal solution to ($\bs{\zeta}$-SAA) for all $n_{\text{eval}}$ contextual problems in the evaluation set.
    \item \textbf{2SP Solve:} Time taken to solve all $n_{\text{eval}}$ contextual two-stage stochastic programs in the evaluation set to optimality.
    \item \textbf{MMD Loss Evaluation:} Time taken to calculate the loss $\ell_{\text{opt}}$ over $S$.
    \item \textbf{Loss-Net Training:} Time taken to train the loss network.
    \item \textbf{Static Training:} Time taken to train the Task-Net via Static-PCSG.
    \item \textbf{Dynamic Training:} Time taken to train the Task-Net via Dynamic-PCSG.
\end{itemize}

The computation times for the processes above are shown in Tables \ref{tab:CVaR_times} and \ref{tab:MNV_times} for both the CVaR and MNV experiments, respectively, and are located in Section \ref{appendix:timing}. Tables \ref{tab:CVaR_times} and \ref{tab:MNV_times} indicate that CSG's trade-off between training time and quick surrogate solution evaluations can be worth it. For example, when $K = 40$, CSG can solve 312 CVaR instances in approximately $33$ seconds using the dynamic approach with $\lambda = 10240$. The resulting solutions have a median optimality gap of $5.07\%$. Similar statements also hold for MNV. Although the CSG methods require additional time to train $\bs{f}$, this time is incurred offline and is justified by repeatedly solving 2SPs in different contexts.

 \section{Conclusion}\label{section:conclusion}

This work introduces CSG, a framework for solving 2SPs in contextual settings. In time-sensitive applications, here-and-now decisions are required that account for the distribution of uncertainty conditioned on available context. It is desirable to efficiently obtain solutions to contextual 2SPs. To this end, CSG proposes learning a mapping from context to a set of $K$ scenarios, which can then be used to solve the deterministic equivalent of the 2SP defined on those scenarios. CSG leverages contextual information to learn task-mappings that produce high-quality decisions while reducing the computational burden of repeatedly solving large-scale stochastic programs.

Motivated by distributional approaches to scenario generation, we propose a distributional approach that leverages the MMD distance to create task mappings that remain close to the true conditional distribution across contexts. This method only requires joint samples, typically available in historical data, and does not need direct access to the true conditional distributions. Furthermore, we propose a problem-driven approach based on a bilevel scenario generation problem that addresses non-unique solutions and sparse gradients. Non-unique solutions are addressed by introducing an optimistic loss function that is easily computable for a large class of 2SPs. Furthermore, we mitigate the issue of sparse gradients by using a neural network to approximate the proposed problem-driven loss. We subsequently use the trained network, along with MMD regularization, to guide gradient-based search for the Task-Net. Lastly, we propose a dynamic approach that extends the static approach by iteratively refining the Loss-Net and Task-Net in an alternating fashion. We observe that the dynamic approach can improve the performance of the static approach if the Loss-Net can sufficiently approximate the problem-driven loss. Furthermore, we explore finite-sample generalization bounds for the proposed methodologies. 

A diverse array of contextual problem settings is employed to demonstrate the framework's effectiveness and illustrate its flexibility. We show that the proposed problem-driven methods rely critically on MMD regularization to produce high-quality Task-Nets. Furthermore, we observe that the proposed methodology can unlock the VSS and compute high-quality solutions in the considered problem settings. Additionally, the offline cost of training the models required for CSG is justified by the time savings in computing solutions to contextual 2SPs. 

Future work could explore the application and impacts of different distributional distance measures, the merits and drawbacks of the pessimistic approach described in Section~\ref{subsection:PCSG}, and extensions to distributionally robust optimization and multi-stage problems. Lastly, problem-specific extensions and structured architectures could reveal new insights in the field of contextual optimization.

\subsection*{Acknowledgements}
The authors received financial support from the joint International Doctoral Cluster between the KAIST and the University of Toronto, Department of Mechanical and Industrial Engineering.

\section*{Declarations}

\subsection*{Funding}
Financial support was provided by NSERC Grant 455963 for David Islip and Roy H. Kwon.

\subsection*{Conflict of interest}
The authors have no relevant financial or non-financial interests to disclose.

% \subsection*{Ethics approval and consent to participate}
% Not applicable.
% \subsection*{Consent for publication}
% Not applicable.
\subsection*{Data availability}
All data used in the experiments is either simulated (and accessible in the Online Supplement) or available via the Tiingo data service (\url{https://www.tiingo.com/}).

% \subsection*{Materials availability}
% Not applicable.
% \subsection*{Code availability}
% The code is available

\subsection*{Author contribution}

\begin{itemize}
  \item David Islip: Writing - review \& editing, Writing - original draft, Validation, Software, Methodology, Formal analysis, Conceptualization.
  \item Roy H. Kwon: Writing - review \& editing, Writing - original draft, Supervision, Resources, Project administration, Funding acquisition, Conceptualization.
  \item Sanghyeon Bae: Writing - review \& editing, Methodology, Formal analysis, Conceptualization
  \item Woo Chang Kim: Writing - review \& editing, Supervision, Resources, Project administration, Funding acquisition, Conceptualization.
\end{itemize}

% \begin{itemize}
% %\item Funding
% %\item Conflict of interest/Competing interests (check journal-specific guidelines for which heading to use)
% \item Ethics approval and consent to participate
% \item Consent for publication
% \item Materials availability
% \end{itemize}

% \noindent
% If any of the sections are not relevant to your manuscript, please include the heading and write `Not applicable' for that section. 

% \begin{appendices}

% Restart page numbering (only do this in the supplement PDF, not the main paper)

\bibliography{sn-bibliography}% common bib file

\clearpage
\setcounter{page}{1}
\pagenumbering{arabic}

% Reset counters
\setcounter{section}{0}
\setcounter{equation}{0}
\setcounter{figure}{0}
\setcounter{table}{0}

% Section numbering: OS1, OS2, ...
\renewcommand{\thesection}{OS-\arabic{section}}
\renewcommand{\thesubsection}{\thesection.\arabic{subsection}}
\renewcommand{\thesubsubsection}{\thesubsection.\arabic{subsubsection}}

% Eq/Fig/Table numbering: OS1.1, OS1.2, ...
\renewcommand{\theequation}{\thesection.\arabic{equation}}
\renewcommand{\thefigure}{\thesection.\arabic{figure}}
\renewcommand{\thetable}{\thesection.\arabic{table}}

% (Optional but recommended) reset eq/fig/table per section
% Requires: \usepackage{chngcntr}
\counterwithin{equation}{section}
\counterwithin{figure}{section}
\counterwithin{table}{section}

% (If using hyperref, avoid duplicate anchors)
\renewcommand{\theHsection}{\thesection}
\renewcommand{\theHsubsection}{\thesubsection}
\renewcommand{\theHequation}{\theequation}
\renewcommand{\theHfigure}{\thefigure}
\renewcommand{\theHtable}{\thetable}

\begin{center}
{\Large\bfseries Online Supplement (OS)}\\[6pt]
{\large\bfseries Contextual Scenario Generation for Two-Stage Stochastic Programming}\\[6pt]
{\normalsize Computational Optimization and Applications}\\[10pt]

David Islip$^{1}$, Roy H.\ Kwon$^{1}$, Sanghyeon Bae$^{2}$, Woo Chang Kim$^{2}$\\[8pt]

{\small
$^{1}$Department of Mechanical and Industrial Engineering, University of Toronto, Toronto, Ontario, Canada\\
$^{2}$Department of Industrial and Systems Engineering, Korea Advanced Institute of Science and Technology (KAIST), Daejeon, Republic of Korea\\[8pt]
\textit{Corresponding author:} Roy H.\ Kwon (\texttt{rkwon@mie.utoronto.ca})
}
\end{center}

\section{Comparing Conditional Distributions}\label{appendix:comparing}
In general, MMD distances differ from Wasserstein distances in several key respects. Computing Wasserstein metrics between two samples of size $n$ requires $O(n^3 \log(n))$ operations, whereas MMD requires $n^2$ operations \citepSM{pele2009fastSM}. Additionally, Wasserstein metrics have dimension-dependent sample complexity and require more samples in higher dimensional settings to bound the gap between the sampled distance and the true distance \citepSM{genevay2019sampleSM}. Wasserstein also suffers from biased sample gradients, unlike MMD \citepSM{bellemare2017cramerSM}. However, MMD induces a geometry on the space of probability measures that does not respect the distance metric on the underlying space. \citetSM{feydy2019interpolatingSM} use Sinkhorn divergence: an entropic regularization of the Wasserstein distance with bias removed. They show that the Sinkhorn divergence interpolates between Wasserstein and MMD of the Laplacian RKHS based on the amount of regularization. The interpolation property also applies to complexities, with Sinkhorn achieving $O(n^2)$ and $O(1/\sqrt{n})$ computational and sample complexity, respectively \citepSM{genevay2019sampleSM}. Still, computing Sinkhorn divergence takes approximately 20--50 times as long as computing MMD \citepSM{feydy2019interpolatingSM}. 

\citetSM{ren2016conditionalSM} and \citetSM{park2020measureSM} consider conditional versions of the MMD metric that only require joint sampling access, however, \citetSM{huang2022evaluatingSM} point out that optimizing these metrics via batched gradients requires matrix inversions that scale with the batch size. Some authors consider $d^2_{\mathcal{F}_\text{MMD}}(\mathbb{P}_{\boldsymbol{x}, \bs{\omega}}, \mathbb{P}_{\boldsymbol{x}, \bs{\xi}})$ to infer a relationship between $\mathbb{P}_{\bs{\omega}|\boldsymbol{x}}$ and $\mathbb{P}_{\bs{\xi}|\boldsymbol{x}}$. For instance, \citetSM{huang2022evaluatingSM} (Theorem 6) claim, under suitable assumptions (measurability, boundedness in expectation and characteristic kernels) that $d^2_{\mathcal{F}_\text{MMD}}(\mathbb{P}_{\boldsymbol{x}, \bs{\omega}}, \mathbb{P}_{\boldsymbol{x}, \bs{\xi}}) = 0 \iff \mathbb{P}_{\bs{\omega}|\boldsymbol{x}} = \mathbb{P}_{\bs{\xi}|\boldsymbol{x}}$ almost everywhere according to $\mathbb{P}_{\boldsymbol{x}}$. In the case of the energy distance, \citetSM{hagemann2023posteriorSM} bounds $\mathbb{E}_{\boldsymbol{x} \sim \mathbb{P}_{\boldsymbol{x}}} \left[ d^2_{\mathcal{F}_\text{MMD}}\left( \mathbb{P}_{\bs{\omega} |\boldsymbol{x}}, \mathbb{P}_{\bs{\xi} |\boldsymbol{x}} \right) \right]$ by $d^2_{\mathcal{F}_\text{MMD}}(\mathbb{P}_{\boldsymbol{x}, \bs{\omega}}, \mathbb{P}_{\boldsymbol{x}, \bs{\xi}})^{\frac{1}{4(1+p+d)}}$ under regularity and compactness assumptions. \citetSM{chemseddine2023diagonalSM} show the opposite relationship holds for the 1-Wasserstein metric and consequently introduce a Wasserstein distance of order $q$, conditioned on $\boldsymbol{x}$, denoted by $d_{q, \boldsymbol{x}}$ that satisfies $\mathbb{E}_{\boldsymbol{x} \sim \mathbb{P}_{\boldsymbol{x}}} \left[ d_q^q \left( \mathbb{P}_{\bs{\omega} |\boldsymbol{x}}, \mathbb{P}_{\bs{\xi} |\boldsymbol{x}} \right) \right] = d_{q, \boldsymbol{x}}^q(\mathbb{P}_{\boldsymbol{x}, \bs{\omega}}, \mathbb{P}_{\boldsymbol{x}, \bs{\xi}})$, where $d_q^q(.,.)$ denotes the $q$-Wasserstein distance raised to the power $q$. In practice, they implement their conditional distance by penalizing transport cost associated with $\boldsymbol{x}$ and use \citetSM{feydy2019interpolatingSM}'s debiased Sinkhorn approach to compute the distance.

\section{Proof of Proposition \ref{prop:representation}}\label{appendix:prop2}
This section provides proof of Proposition \ref{prop:representation}. The proof largely relies on Proposition 1 due to \hbox{\citetSM{tabaghi2024universalSM}}. First, we introduce some notation, followed by Proposition 1 by \citetSM{tabaghi2024universalSM}, and lastly, the proof of Proposition \ref{prop:representation}. 

\subsection*{Multiset Notation} A multiset is a pair $(U, \Gamma)$ where $U$ denotes a set of objects and $\Gamma$ maps $U$ to the non-negative integers, where $\Gamma(u)$ represents the cardinality of $u \in U$. Multisets are denoted by double brackets $\{\{ \}\}$. For example, the multiset $\{\{3, 3, 4\}\}$ has three elements with $\Gamma(3) = 2$ and $\Gamma(4) = 1$. For any domain $\Omega$, a multiset $U$ such that $U \subseteq \Omega$ means that the elements of $U$ are in $\Omega$. The cardinality of $U$, denoted by $|U|$ is the number of elements in $U$, counting repetitions, e.g. $|\{\{3, 3, 4\}\}| = 3$. For some $N$ in the positive integers $\mathbb{N}$ and domain $\Omega$, the following sets are introduced:
$$
\mathbb{U}_{\Omega, N} = \{ \text{multiset} \ U \subseteq \Omega: |U| = N \}\ 
\text{and}\ 
\mathbb{U}_{\Omega, [N]} = \{ \text{multiset} \  U \subseteq \Omega: |U| \in [N] \}.
$$
\subsection*{\citetSM{tabaghi2024universalSM}'s Proposition 1} \citetSM{tabaghi2024universalSM} prove the following proposition regarding functions of multisets, which we repeat here for the reader's convenience. 

\begin{proposition}[Proposition 1 from \protect\citepSM{tabaghi2024universalSM}]\label{prop:two_multisets}
    A (continuous) multiset function $f:\mathbb{U}_{\Omega,[N_1]} \times \mathbb{U}_{\Omega,[N_2]} \rightarrow \text{codom}(f)$, where $\Omega$ is a compact subset of $\ \mathbb{R}^{p}$, is (continuously) sum-decomposable via $\mathbb{R}^{{N+p \choose p}-1} \times \mathbb{R}^{{N+p \choose p}-1}$, that is, 
    \[
         f(U, U^{\prime}) = \tilde{\rho} \left( \sum_{\bs{u} \in U} \bs{\tau}(\bs{u}), \sum_{\bs{u}^{\prime} \in U^{\prime}} \bs{\tau}(\bs{u}^{\prime}) \right) \quad \forall U \in \mathbb{U}_{\Omega,[N_1]}, U^{\prime} \in \mathbb{U}_{\Omega,[N_2]}, 
    \]
    with continuous $\bs{\tau}: \mathbb{R}^{p} \rightarrow \mathbb{R}^{{N+p \choose p}-1}, N = \max \{N_1, N_2 \}$ and (continuous) $\tilde{\rho}: \mathbb{R}^{{N+p \choose p}-1} \times \mathbb{R}^{{N+p \choose p}-1} \rightarrow \text{codom}(\tilde{\rho})$, and $\text{codom}(f) \subset \text{codom}(\tilde{\rho})$.
\end{proposition}

\subsection*{Proof}
For the reader's convenience, we repeat Proposition \ref{prop:representation}, followed by its proof.
\representation*
\begin{proof}
        First, fix a value of $K \in \mathbb{N}$. We can identify $\ell_{\text{opt}}$ as a multiset function. Since the value $\ell_{\text{opt}}(\bs{\zeta}_{1 \hdots K^{\prime}}, \bs{\omega})$ does not depend on the order of the $K^{\prime} \in \mathbb{N}$ scenarios $\bs{\zeta}_{1 \hdots K^{\prime}}$, we can define the multisets $U_{\bs{\zeta}_{1 \hdots K^{\prime}}} = \{\{ \bs{\zeta}_k \}\}_{k = 1}^{K^{\prime}} \in \mathbb{U}_{\Omega, K^{\prime}}$ and $U_{\bs{\omega}} = \{\{ \bs{\omega} \}\} \in \mathbb{U}_{\Omega, 1}$ and define $\ell_{\text{opt}}^{\prime} : \mathbb{U}_{\Omega, [K^{\prime}]} \times \mathbb{U}_{\Omega, 1} \rightarrow \text{codom}(\ell_{\text{opt}})$ such that
        $$
        \ell_{\text{opt}}^{\prime} \left( U_{\bs{\zeta}_{1 \hdots K^{\prime}}}, U_{\bs{\omega}} \right) = \ell_{\text{opt}}(\bs{\zeta}_{1 \hdots K^{\prime}}, \bs{\omega}) \quad \forall \bs{\zeta}_{1 \hdots K^{\prime}} \in \Omega^{K^{\prime}}, \bs{\omega} \in \Omega, \ \forall K^{\prime} \in [K]. 
        $$
        The value of $\ell_{\text{opt}}^{\prime} \left( U_{\bs{\zeta}_{1 \hdots K^{\prime}}}, U_{\bs{\omega}} \right)$ is calculated by solving ($\bs{\zeta}$-SAA) on the multiset of scenarios $U_{\bs{\zeta}_{1 \hdots K^{\prime}}}$ and optimistically evaluating the resulting solutions on the singleton-multiset $U_{\bs{\omega}}$. Applying Proposition 1 from \citepSM{tabaghi2024universalSM}, to $\ell_{\text{opt}}^{\prime}$ implies 
        \begin{equation} \label{eqn:useful_equation}
            \ell_{\text{opt}}^{\prime} \left( U, U^{\prime} \right) = \tilde{\rho} \left( \sum_{\bs{u} \in U} \bs{\tau}(\bs{u}), \sum_{\bs{u}^{\prime} \in U^{\prime}} \bs{\tau}(\bs{u}^{\prime}) \right) \quad \forall U \in \mathbb{U}_{\Omega,[K]}, U^{\prime} \in \mathbb{U}_{\Omega,[1]}, 
        \end{equation}
        where $N = \max \{K, 1 \} = K$, continuous $\bs{\tau}: \mathbb{R}^{p} \rightarrow \mathbb{R}^{{K+p \choose p}-1}$, (continuous) $\tilde{\rho}: \mathbb{R}^{{K+p \choose p}-1} \times \mathbb{R}^{{K+p \choose p}-1} \rightarrow \text{codom}(\tilde{\rho})$, and $\text{codom}(\ell_{\text{opt}}^{\prime}) \subset \text{codom}(\tilde{\rho})$.
        Since \eqref{eqn:useful_equation} holds for all $U \in \mathbb{U}_{\Omega,[K]}$,  it also holds for all $U \in \mathbb{U}_{\Omega,K}$, implying 
        \begin{align*}
        \ell_{\text{opt}}(\bs{\zeta}_{1 \hdots K}, \bs{\omega})  &= \tilde{\rho} \left( \sum_{\bs{\zeta} \in U_{\bs{\zeta}_{1 \hdots K}}} \bs{\tau}(\bs{\zeta}), \sum_{\bs{\omega} \in U_{\bs{\omega}}} \bs{\tau}(\bs{\omega}) \right) \\
        &= \tilde{\rho} \left( \sum_{k = 1}^K \bs{\tau}(\bs{\zeta}_k),  \bs{\tau}(\bs{\omega}) \right) \quad \forall \bs{\zeta}_{1 \hdots K} \in \Omega^{K}, \bs{\omega} \in \Omega.
        \end{align*}
        Lastly, letting $\rho : \mathbb{R}^{{K+p \choose p}-1} \times \mathbb{R}^{{K+p \choose p}-1} \rightarrow \text{codom}(\tilde{\rho})$ be defined such that
        $$
        \rho(\widehat{\bs{\zeta}}, \hat{\bs{\omega}}) = \tilde{\rho}(K \widehat{\bs{\zeta}}, \hat{\bs{\omega}}) \quad \forall \widehat{\bs{\zeta}} \in \mathbb{R}^{{K+p \choose p}-1},  \hat{\bs{\omega}} \in \mathbb{R}^{{K+p \choose p}-1},
        $$
        yields the desired result. 
\end{proof}
\section{Experimental Information}

\subsection{Experimental Tools}\label{appendix:experimental_tooling}
The mathematical programming models are implemented in Python 3.7 and solved using Gurobi 11.0 on \textsf{Google Colab} equipped with an Intel(R) Xeon(R) CPU 2.20GHz, 51 gigabytes of random access memory (RAM) and a Tesla T4 GPU. All neural networks are trained using \textsf{skorch} \citepSM{skorchSM}. 

\subsection{Newsvendor Scenario Comparisons}\label{appendix:newsvendor_scenarios}
\begin{figure}[h]
    \FIGURE
    {\includegraphics[width=0.99\textwidth]{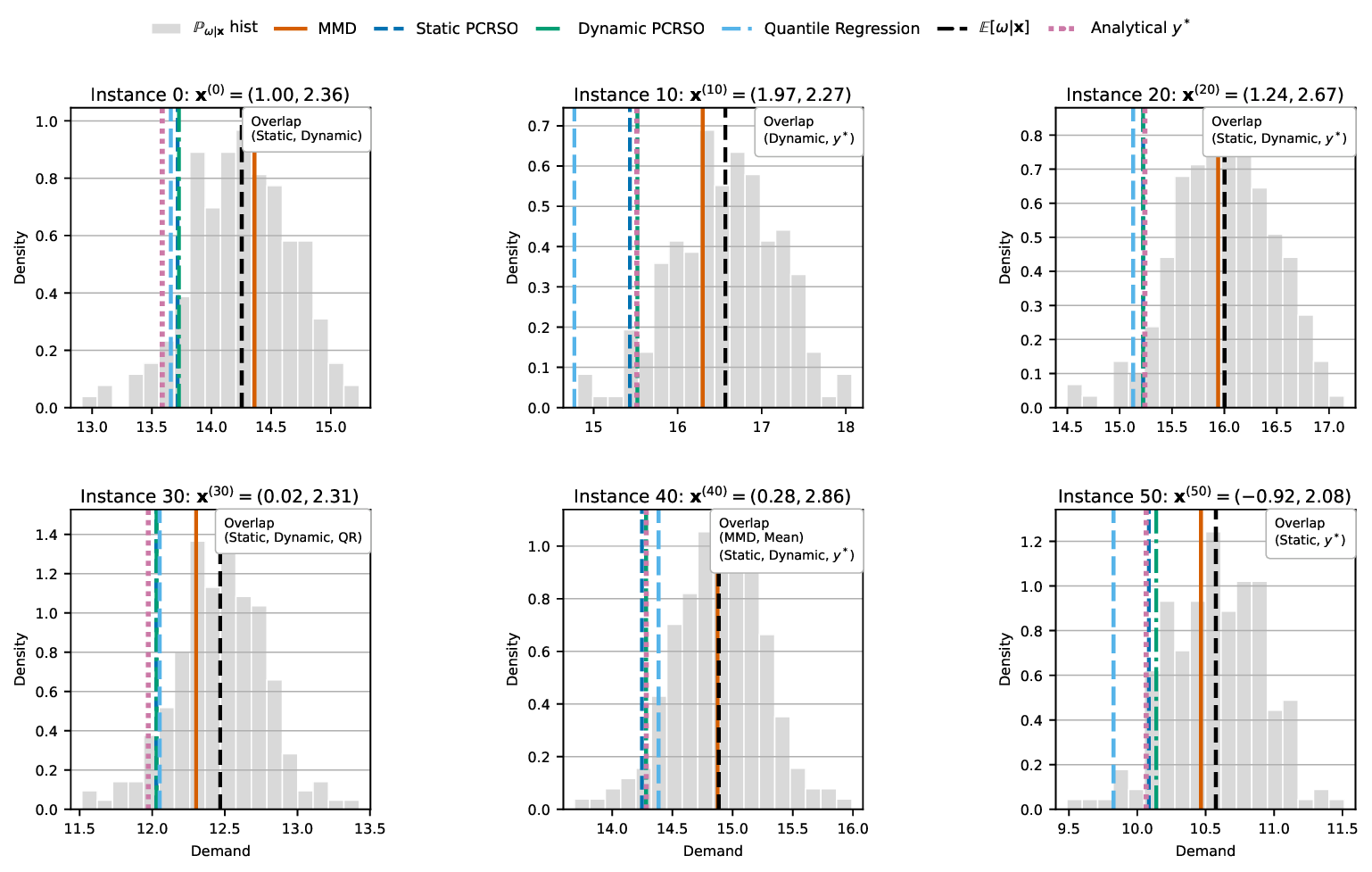}}
    {Sample Scenarios for the contextual newsvendor problem ($\lambda = 1.0$).\label{fig:newsvendor_scenarios}}
    {}
\end{figure}

Figure \ref{fig:newsvendor_scenarios} shows how the proposed methodologies differ in the scenario they generate for the newsvendor problem when $K=1$. Six different contexts are sampled from $\mathbb{P}_{\bs{x}}$, their conditional distributions $\mathbb{P}_{\bs{\omega}|\bs{x}}$ are plotted as histograms, and the scenarios generated by each method are overlaid. From Figure \ref{fig:newsvendor_scenarios}, we can see that the scenarios generated by the proposed problem-driven methods tend to approximate the optimal quantile $y^*$, whereas the proposed distributional method, which minimizes the MMD distance between the generated scenario and the true distribution, generates a scenario that is closer to $\mathbb{E}[\bs{\omega}|\bs{x}]$.

\subsection{CEP1 Experimental Supplement}\label{appendix:cep1_formulation_data}
\subsubsection{CEP1 Formulation} \label{appendix:cep1_formulation}
% Recall in CEP1, there are $N$ parts being produced on $n_{\text{machines}}$. Machine $j \in [n_{\text{machines}}]$ is available for $h_j$ hours per week, with the ability to purchase additional hours at an hourly cost of $c_j$. 
 Based on the CEP1 problem description, the decision-maker wishes to solve:
    \begin{align}\label{eq:cep1}
    \min_{\bs{y}}\quad   \bs{c}^{\intercal} \bs{y}_{\text{cap}} +  \mathbb{E}_{\bs{\omega} \sim \mathbb{P}_{\bs{\omega}|\bs{x}}}& \left[ Q(\bs{y}_{\text{op}} , \bs{\omega})\right]  \tag{CEP1}  \\ \textit{s.t.}  \quad  -y_{\text{cap}, j} + y_{\text{op},j} & \leq h_j  & \forall j \in [n_{\text{machines}}]  \notag \\
    \quad  \bs{t}^{\intercal} \bs{y}_{\text{op}} & \leq T \notag \\
    \quad  \bs{0} \leq \bs{y}_{\text{cap}},\ & \bs{0} \leq \bs{y}_{\text{op}} \leq \bs{u}, \notag
    \end{align}
where $\bs{y}_{\text{cap}}$ and $\bs{y}_{\text{op}}$ denote the hours of added capacity and total operation of the machines, respectively. $\bs{c} \in \mathbb{R}_+^{n_{\text{machines}}}$ denotes the capacity cost vector, $\bs{h} \in \mathbb{R}_+^{n_{\text{machines}}}$ denotes the baseline capacities of the machines, $\bs{t} \in \mathbb{R}_+^{n_{\text{machines}}}$ is the maintenance requirements incurred per hour for all the machines, $T$ is the total maintenance limit, and $\bs{u} \in \mathbb{R}_+^{n_{\text{machines}}}$ denotes the upper bounds on utilization for each machine. The recourse cost is given by
    \begin{align}\label{eq:cep1stage2}
    Q(\bs{y}_{\text{op}}, \bs{\omega}) = \min_{\bs{z}, \bs{s}} \quad   \sum_{i=1}^m \sum_{j=1}^{n_{\text{machines}}} g_{ij} z_{ij} &+ \sum_{i=1}^m p_i s_i  \tag{CEP1-Stage II} \\ \textit{s.t.}  \quad  \sum_{j=1}^{n_{\text{machines}}} a_{ij} z_{ij} +s_i & \geq \omega_i  & \forall i \in [m]  \notag \\
    \quad  \sum_{i=1}^m z_{ij} & \leq y_{\text{op}, j} & \forall j \in [n_{\text{machines}}] \notag \\
    \quad  z_{ij} & \geq 0 & \forall i \in [m], \ j \in [n_{\text{machines}}]  \notag \\
    \quad s_i & \geq 0 & \forall i \in [m], \notag
    \end{align} 
 where $z_{ij}$ and $s_i$ denote the hours producing part $i$ on machine $j$ and the production shortfall of part $i$ respectively. The uncertain demand for part $i$ is $\omega_i$. Machine $j$ produces part $i$ at a rate of $a_{ij} \geq 0$ with an hourly cost of $g_{ij} \geq 0$. Production shortfalls incur a penalty of $p_i \geq 0$. We use the same deterministic problem data as defined by \citepSM{higle_sen_1996SM} with $n_{\text{machines}} = 4$ and $m = 3$. The static data below is used to define the CEP1 instance considered in the manuscript. 
\[
n = 4, \quad m = 3, \quad T = 100, \quad
\bs{c} = (2.5, 3.75, 5.0, 3.0), \quad \bs{t} = (0.08, 0.04, 0.03, 0.01)
\]
\[
\bs{h} = (500, 500, 500, 500), \quad \bs{u} = (2000, 2000, 3000, 3000), \quad
\bs{p} = (400, 400, 400)
\]
\[
[a_{ij}] = \begin{bmatrix}
0.6 & 0.6 & 0.9 & 0.8 \\
0.1 & 0.9 & 0.6 & 0.8 \\
0.05 & 0.2 & 0.5 & 0.8
\end{bmatrix}, \quad 
[g_{ij}] = \begin{bmatrix}
2.6 & 3.4 & 3.4 & 2.5 \\
1.5 & 2.3 & 2.0 & 3.6 \\
4.0 & 3.8 & 3.5 & 3.2
\end{bmatrix}
\]
\subsubsection{CEP1 Stochastic Modelling} \label{appendix:cep1_data_generation}
\citetSM{higle_sen_1996SM} originally consider $\omega_i$, $ i \in \{1,2,3\}$ as an iid discrete random variable with $\omega_i \in \{0, 600, 1200, 1800, 2400, 3000\}$ being equally likely. We construct our contextual setting similarly. The contextual information $\bs{x}$ is distributed according to a four-dimensional normal distribution $N(\bs{0}, \bs{I})$. A neural network $\bs{h}_{\text{Random}} : \mathcal{X} \rightarrow \mathbb{R}^{6 \times m}$ (6 rows and $m$ columns) is randomly initialized. For part $i \in \{1, 2, 3\}$, the $i$th column of $\bs{h}_{\text{Random}}(\bs{x})$ corresponds to the $6$ points defining the support of $\omega_i$. Since each demand $\omega_i$ is supported on 6 points, there are $6^3 = 216$ possibilities defining the conditional demand distribution. Scaling is applied to the outputs from $\bs{h}_{\text{Random}}$ such that the demands expectation and standard deviation are $1500$ and $1024.7$, respectively. Any negative demand values are clipped to 0. This process defines $\mathbb{P}_{\bs{\omega} | \bs{x}}$. 

\subsubsection{CEP1 Optimality Gap CDFs}\label{appendix:cep1_performance}
\begin{figure}[h]
    \FIGURE
    {\includegraphics[width=0.93\textwidth]{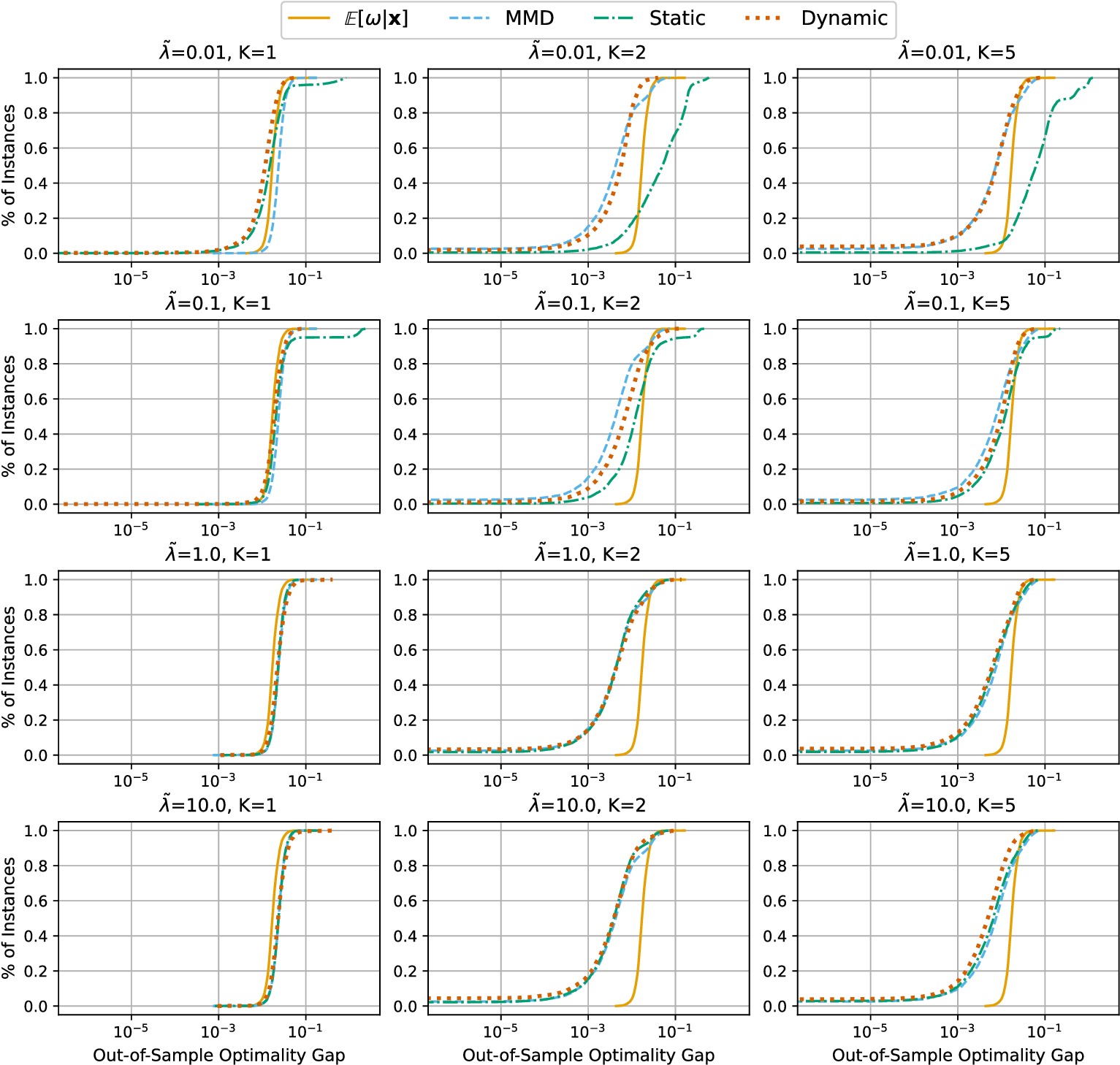}}
    {Out of sample optimality gaps for the contextual CEP1 problem.\label{fig:cep1_gaps}}
    {}
\end{figure}

\newpage 
\subsection{CVaR Experimental Supplement}
\subsubsection{CVaR Formulation} \label{appendix:cvar_formulation}
We consider a trader who wishes to balance the trade-off between risk and expected return. CVaR is used as the risk measure, and we employ the formulations introduced by \citetSM{krokhmal2002portfolioSM} to model the trader's problem. The trader wishes to solve:
 \begin{align}\label{eq:CVaR}
 \min_{\bs{y}, \gamma, \bs{t}}\quad   - \mathbb{E}_{\bs{\omega} \sim \mathbb{P}_{\bs{\omega}|\bs{x}}}[\bs{\omega}]^{\intercal} \bs{y}  
        &+  \lambda_{\text{risk}} \left( \gamma + \frac{1}{1-\alpha} \mathbb{E}_{\bs{\omega} \sim \mathbb{P}_{\bs{\omega}|\bs{x}}} \left[ Q(\bs{y} , \gamma,  \bs{\omega})\right] \right) \tag{CVaR} \\
 \textit{s.t.} \quad  &\bs{1}^{\intercal} \bs{y} = 1, \notag \\
        &y_i \leq t_i \quad \forall i \in [n_{\text{assets}}],  \notag \\
        &\bs{1}^{\intercal} \bs{t} \leq K_{\text{assets}}, \notag  \\
        &\bs{y} \in \mathbb{R}_+^{n_{\text{assets}}},\ \gamma \in \mathbb{R}, \ \bs{t} \in \{0,1\}^{n_{\text{assets}}}, \notag
 \end{align} 
 where $\bs{y}$, $\gamma$, and $\bs{t}$ represent the selected portfolio, a placeholder variable that at optimality equals the portfolio's value at risk (VaR) according to $\mathbb{P}_{\bs{\omega} | \bs{x}}$, and the binary indicator variables that model $t_i = 0 \implies y_i = 0, \ \forall i \in [n_{\text{assets}}]$, respectively. The variables $\bs{y}$, $\gamma$, and $\bs{t}$ constitute the first-stage decisions. The relative risk aversion and VaR confidence level are denoted by $\lambda_{\text{risk}} > 0$ and $\alpha \in (0,1)$, respectively. The $n_{\text{assets}}$ random returns for each asset are represented by $\bs{\omega}$. Lastly, the recourse cost is given by
 \begin{align}\label{eq:CVaRstage2} \tag{CVaR}
 Q(\bs{y}, \gamma,  \bs{\omega}) = \min_{z} \quad  z \notag  \\ \textit{s.t.}  \quad  z & \geq -\bs{\omega}^{\intercal} \bs{y} - \gamma  \notag \\
 \quad  z & \geq 0, \notag
 \end{align}
where $z$ is a variable that models the portfolio loss exceeding the value at risk. We consider the following instance parameters $\lambda_{\text{risk}} = 10$, $n_{\text{assets}} = 10$, $K_{\text{assets}} = 5$, and $\alpha = 0.9$.

\subsubsection{CVaR Stochastic Modelling} \label{appendix:cvar_data_generation}

A data-driven approach is considered to set up the contextual environment. The decision-maker limits the selection to tickers in the S\&P 500 as of December 29, 2023, and chooses the \(n_{\text{assets}}\) tickers with the smallest traded volumes as their asset universe. They desire the ability to obtain portfolios that are of high quality according to \eqref{eq:CVaR} every 5 minutes for hedging purposes. The context window $W$ is set to $78$ periods (a full trading session). Furthermore, as discussed in the text, the trader assumes that the returns satisfy the following relationship: 
\begin{equation}\label{eq:appendix_returns_law}
\bs{\omega} = \Phi_{\bar{\bs{r}}}(\boldsymbol{x}) + \Phi_{\bs{\sigma}} (\boldsymbol{x}) \odot \boldsymbol{\epsilon}, 
\end{equation} 
where $\Phi_{\bar{\bs{r}}} : \mathcal{X} \rightarrow \mathbb{R}^{n_{\text{assets}}}$ and $\Phi_{\bs{\sigma}} : \mathcal{X} \rightarrow \mathbb{R}_+^{n_{\text{assets}}}$ are models given to the trader by a statistical modeler to estimate the conditional mean $\mathbb{E}[\bs{\omega} |\bs{x}]$, and conditional standard deviation $\sqrt{\text{Var}[\bs{\omega} | \bs{x}]}$ along with $\boldsymbol{\epsilon}$ denoting the random errors. There is a large dataset $S^* = \{\bs{x}^{(i)}, \bs{\omega}^{(i)} \}_{i=1}^{n^*}$ of size $n^*$ used to produce $\Phi_{\bar{\bs{r}}}$ and $\Phi_{\bs{\sigma}} (\boldsymbol{x})$. Similar to \citetSM{deng2022predictiveSM} and \citetSM{ban2019dynamicSM}, the residuals obtained from the trained models define the conditional distribution. That is, $\bs{\epsilon}$ is assumed to be uniformly distributed over $\{\bs{\epsilon}^{(i)}\}_{i=1}^{n^*}$, where 
\begin{equation}\label{eq:residuals}
\bs{\epsilon}^{(i)} = \left(\bs{\omega}^{(i)} - \Phi_{\bar{\bs{r}}}(\boldsymbol{x}^{(i)}) \right) \odot \begin{pmatrix}
    1/\Phi_{\sigma, 1} (\boldsymbol{x}^{(i)}) \\
    1/\Phi_{\sigma, 2} (\boldsymbol{x}^{(i)}) \\
    \vdots \\
    1/\Phi_{\sigma, n_{\text{assets}}} (\boldsymbol{x}^{(i)})
    \end{pmatrix}, \quad \forall i \in [n^*].
\end{equation}
Thus, the trader aims to solve \eqref{eq:CVaR} with $\mathbb{P}_{\bs{\omega}|\bs{x}}$ defined by \eqref{eq:appendix_returns_law} and \hbox{$\bs{\epsilon} \sim \text{uniform} \{\bs{\epsilon}^{(i)}\}_{i=1}^{n^*}$}. 

A dataset of prices consisting of $21$ trading periods is gathered for each ticker via the IEX endpoint available via the Tiingo API, resulting in $1638$ returns $\bs{r}^{(j)} \in \mathbb{R}^{n_{\text{assets}}}, \ j = 1,\hdots, 1638$. The sample $S^*$ used to train $\Phi_{\bar{\bs{r}}}$ and $\Phi_{\bs{\sigma}}$ is formed as follows:
\[
\bs{x}^{(i)} = \begin{pmatrix}
 \bs{r}^{(i)} & \dots & \bs{r}^{(i+W)} \\
    & (\textsf{PE}^{(1)}_{i \hdots i+W})^{\intercal} & \\
    & \vdots & \\
    & (\textsf{PE}^{(n_{\text{encoding}})}_{i \hdots i+W})^{\intercal}&
 \end{pmatrix}\  \text{and} \ \bs{\omega}^{(i)} = \bs{r}^{(i+W+1)}, \ i = 1,\hdots,  n^* \coloneq 1559,
\]
where $\textsf{PE}^{(j)}_{i \hdots i+W} \in \mathbb{R}^{W}$ is the $j$th sinusoidal positional encoding of the period associated with indices $i,\hdots, W+i$ (the period of the $i$th datum is given by $i$ modulo $78$) \citepSM{vaswani2017attentionSM}. We design $\Phi_{\bar{\bs{r}}}$ and $\Phi_{\boldsymbol{\sigma}}$ such that the resulting distribution can generate the dataset $S^*$ as a sample. Specifically,  $\Phi_{\bar{\bs{r}}}$ and $\Phi_{\boldsymbol{\sigma}}$ are parameterized by long short-term memory (LSTM) neural networks. The parameters of these networks are selected via random search, using a $20\%$ holdout for validation. 

The networks $\Phi_{\bar{\bs{r}}}$ and $\Phi_{\bs{\sigma}}$ are fit in a two-step training procedure. First we train $\Phi_{\bar{\bs{r}}}$ by optimizing \hbox{$\frac{1}{n}\sum_{i=1}^n\| \hat{\boldsymbol{r}}_i - \Phi_{\bar{\bs{r}}}(\boldsymbol{x}_i)\|_1$}. Then the vectors of residuals $\hat{\boldsymbol{\eta}}_i = |\hat{\boldsymbol{r}}_i - \hat{\Phi}_{\bar{\bs{r}}}(\boldsymbol{x}_i)|$ are formed and $\Phi_{\bs{\sigma}}$ is obtained by minimizing \hbox{$\frac{1}{n}\sum_{i=1}^n\| \hat{\boldsymbol{\eta}}_i - \Phi_{\bs{\sigma}}(\boldsymbol{x}_i)\|_1$}. Lastly, equation \eqref{eq:residuals} yields $\{\bs{\epsilon}^{(i)}\}_{i=1}^{n^*}$.

% \subsubsection{CVaR Results} \label{appendix:CVaR_results} 

\subsubsection{CVaR Optimality Gap CDFs}\label{appendix:cvar_performance}

    \begin{figure}[H]
        \FIGURE
        {\rotatebox{90}{\includegraphics[width=0.56\paperheight]{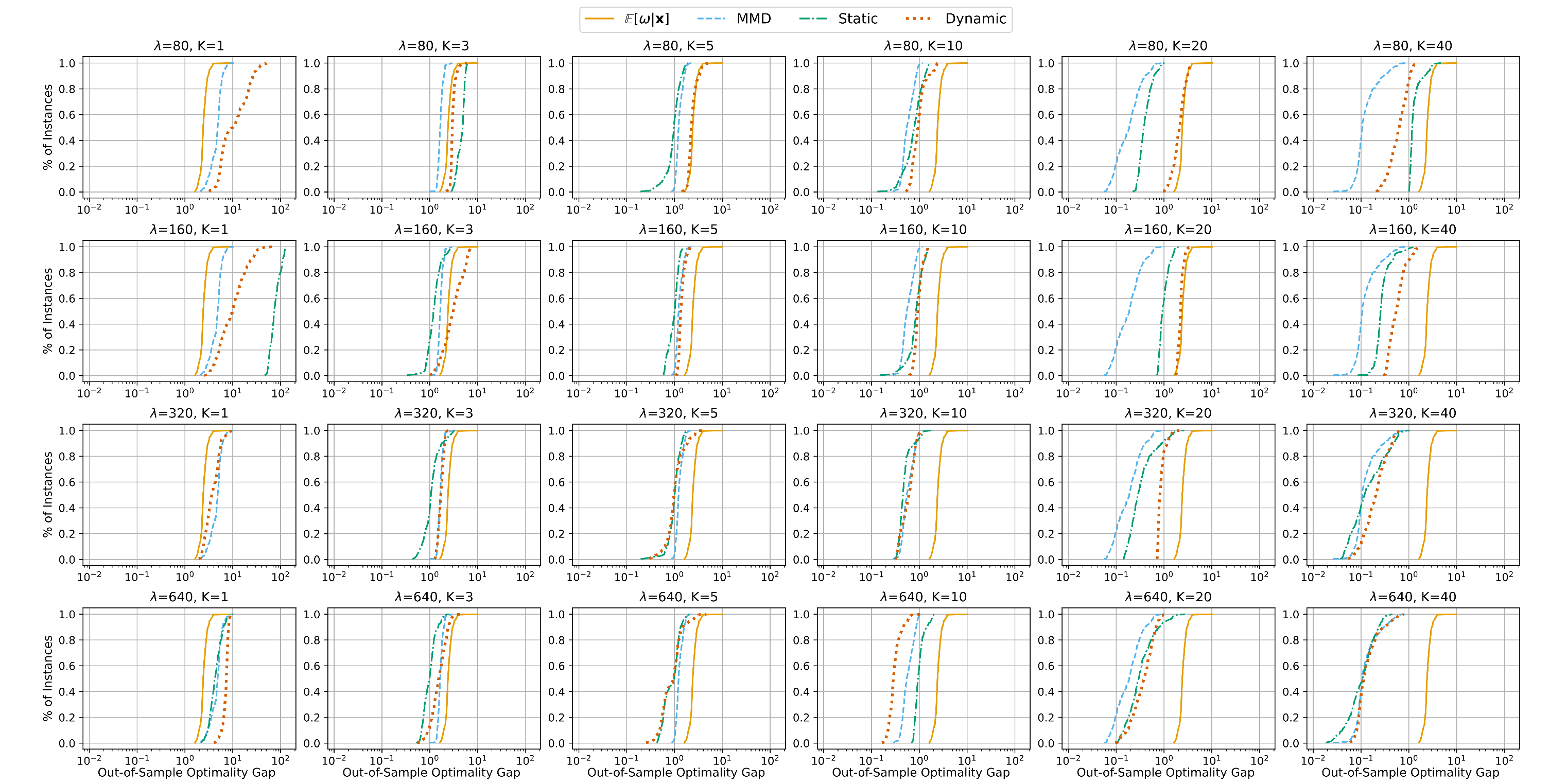}}}
        {Out of sample optimality gaps for the contextual CVaR problem \label{fig:cvar_gaps1}} {Figure continued on next page}
    \end{figure}

    \begin{figure}[H]
            \FIGURE
            {\rotatebox{90}{\includegraphics[width=0.56\paperheight]{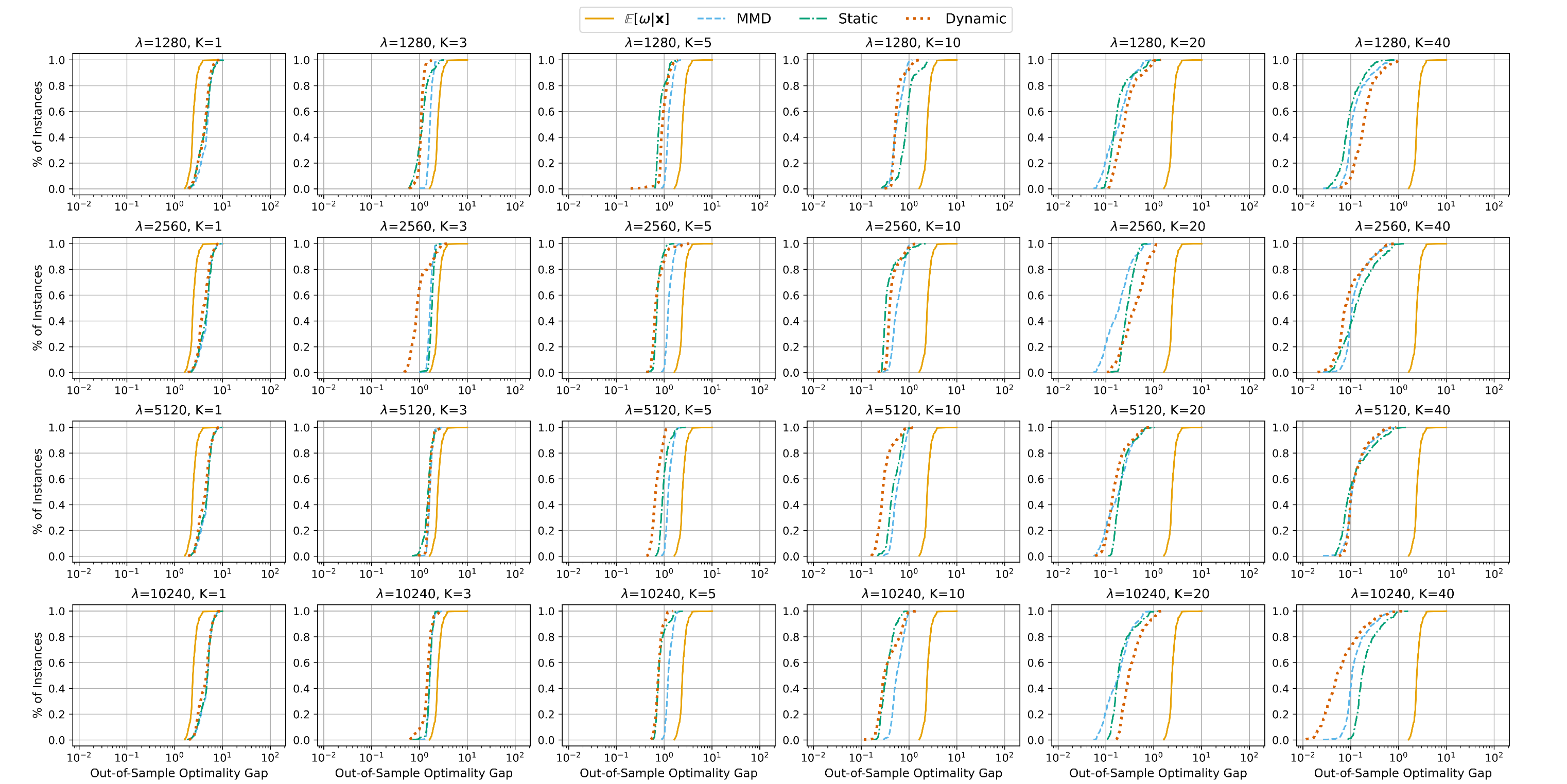}}}
            {Out of sample optimality gaps for the contextual CVaR problem \label{fig:cvar_gaps2}} {}
        \end{figure}

\subsection{MNV Experimental Supplement}
\subsubsection{MNV Formulation}\label{appendix:MNV_formulation}
This section presents the MNV problem as considered by \citetSM{narum2024problemSM}. There are $N$ products available for production. The decision-maker wishes to select a production plan to produce at most $P$ items while respecting the overall production capacity $C$. They wish to solve:
        \begin{align}\label{eq:mnv_firststage} 
        \max_{\bs{y}, \bs{t}} \quad \mathbb{E}_{\bs{\omega} \sim \mathbb{P}_{\bs{\omega}|\bs{x}}}\left[ Q(\bs{y}, \bs{\omega}) \right] - \sum_{i=1}^m c_i y_i  \tag{MNV} \\
        \textit{s.t.} \quad \sum_{i=1}^m y_i &\leq C  \notag \\
        \quad \sum_{i=1}^m t_i &\leq P  \notag \\
        \quad y_i &\leq M t_i & \forall i \in [m]  \notag \\
        \quad \bs{y} \in \mathbb{R}^m_+, & \ \bs{t} \in \{0, 1\}^m, \notag
        \end{align}
    where \(y_i\) represents the amount of product $i$ produced, \(t_i\) represents the binary decision indicating that product $i$ is eligible for production, and \(c_i\) is the product $i$ production cost. 

    The recourse profit is the total of the salvage values and the sales values, which come from both direct sales and customer-driven substitutions. The recourse profit is given by
        \begin{align}\label{eq:mnv_secondstage}
        Q(y, \bs{\omega}) = \max_{\bs{s}, \bs{z}, \bs{\bar{z}}, \bs{w}, \bs{\hat{z}}} \quad \sum_{i=1}^m \left( v_i s_i + v_i \bar{z}_i + g_i w_i \right)  \tag{MNV-Stage II} \\
        \textit{s.t.} \quad s_i + \sum_{j \in [m] :j \neq i} z_{ji} &\leq \omega_i & \forall i \in [m]  \notag \\
        \quad z_{ij} &\leq \alpha_{ij} \left( \omega_j - s_j \right) & \forall i \in [m], j \in [m], i \neq j  \notag \\
        \quad \bar{z}_i &= \sum_{j \in [m]:j \neq i} z_{ij} & \forall i \in [m]  \notag \\
        \quad M \left( \hat{z}_j - 1 \right) &\leq s_j - y_j & \forall j \in [m]  \notag \\
        \quad z_{ij} &\leq M \hat{z}_j & \forall i \in [m], j \in [m]  \notag \\
        \quad w_i &= y_i - \left( s_i + \bar{z}_i \right) & \forall i \in [m]  \notag \\
        \quad \hat{z}_{i} &\in \{0, 1\} &\ \forall i \in [m] \notag \\
        \quad \bs{s}, \bs{z}, \bs{\bar{z}}, \bs{w} &\in \mathbb{R}^m_+, \notag
        \end{align} 
where \(s_i\) represents the sales of product $i$, \(z_{ij}\) is the substitution sale amount of item $i$ to satisfy demand for item $j$, \(\bar{z}_i\) is the total amount of item $i$ substituted, \(w_i\) represents unsold inventory that is salvaged, and $\hat{z}_j$ is a binary variable indicating whether to start substitution sales satisfying demand for item $j$. 

 In the objective, parameters $v_i$ and $g_i$ denote the sales price and salvage value of item $i$. The substitution rate $\alpha_{ij}$ is the average probability that item $j$ can be replaced by item $i$. $\alpha_{ij}$ is not necessarily symmetric, e.g., pink t-shirts can be substituted with white t-shirts more often than white shirts can with pink. The big-$M$s are set to be the capacity $C$. Lastly, $\bs{\omega} \in \mathbb{R}^m_+$ denotes the uncertain demand for the $N$ products. We consider a setting with $m = 6$ potential products. The static data below is used to define the MNV instance considered in the manuscript. 
\[
\boldsymbol{v} = \left( 
46.111, \ 
44.691, \ 
46.448, \ 
48.406, \ 
44.476, \ 
44.476 
\right), 
\]
\[
\boldsymbol{c} = \left( 
18.980, \ 
25.618, \ 
24.163, \ 
26.217, \ 
17.974, \ 
26.418 
\right),
\]
\[
\boldsymbol{g} = \left( 
9.830, \ 
5.094, \ 
5.067, \ 
5.293, \ 
5.723, \ 
7.292 
\right)
\]
\[
[\alpha_{ij}] = \begin{bmatrix}
- & 0.087 & 0.184 & 0.042 & 0.088 & 0.110 \\
0.137 & - & 0.060 & 0.154 & 0.178 & 0.014 \\
0.182 & 0.051 & - & 0.285 & 0.290 & 0.243 \\
0.091 & 0.029 & 0.205 & - & 0.037 & 0.149 \\
0.010 & 0.273 & 0.078 & 0.199 & - & 0.156 \\
0.164 & 0.055 & 0.291 & 0.233 & 0.282 & -
\end{bmatrix},
\]
along with the constants $C = 70$ and  $P = 3$.

\subsubsection{MNV Stochastic Modelling}\label{appendix:MNV_data}

In the case of fashion retail, \citetSM{vaagen2011modellingSM} point out that demand is multi-modal with strong dependence. In light of this, \citetSM{narum2024problemSM} model demand with a multivariate mixture distribution with binary stochastic variables to determine the regime for the product, and two Normal distributions for the specific demand with a large or small mean depending on the regime. We reflect these considerations in the contextual setup for the problem, which is as follows. 

The context $\bs{x} \in \mathbb{R}^p$ is assumed to be Normally distributed $N(\bs{0}, \bs{\Sigma})$ where $p=10$ and $\Sigma_{ij} = \rho^{|i - j|}$ with $\rho = 0.8$. We randomly initialize a matrix $\bs{W} \in \mathbb{R}^{m \times p}$ where the elements of $\bs{W}$ come from a standard Normal distribution. The random demand state $Z \in \{-1,1\}^{m}$ is given by $\text{sign}\left[\bs{W} \bs{x} + \bs{\epsilon}\right]$, where $\bs{\epsilon} \in \mathbb{R}^m$ represents the idiosyncratic demand uncertainty for each product. $\bs{\epsilon}$ is assumed to be multivariate Normal $N(\bs{0}, \sigma \bs{I})$, where $\sigma$ is set such that the signal-to-noise ratio is $0.5$. Based on product $i$'s demand state $Z_i$, the demand $\omega_i$ is distributed according to $N(\mu_{Z_i}, \sigma^2_{Z_i})$. We set $\mu_{-1} = 5$, $\mu_{+1} = 15$, $\sigma_{-1} = 5/3$, and $\sigma_{+1} = 5/2$. The sampled demand is replaced with $0$ if it is negative. To define $\mathbb{P}_{\bs{\omega} | \bs{x}}$ we sample $80$ observations via the aforementioned process.

\newpage
\subsubsection{MNV Optimality Gap CDFs}\label{appendix:MNV_performance}
\ \ 
\begin{figure}[H]
\FIGURE
    {\includegraphics[width=0.99\textwidth]{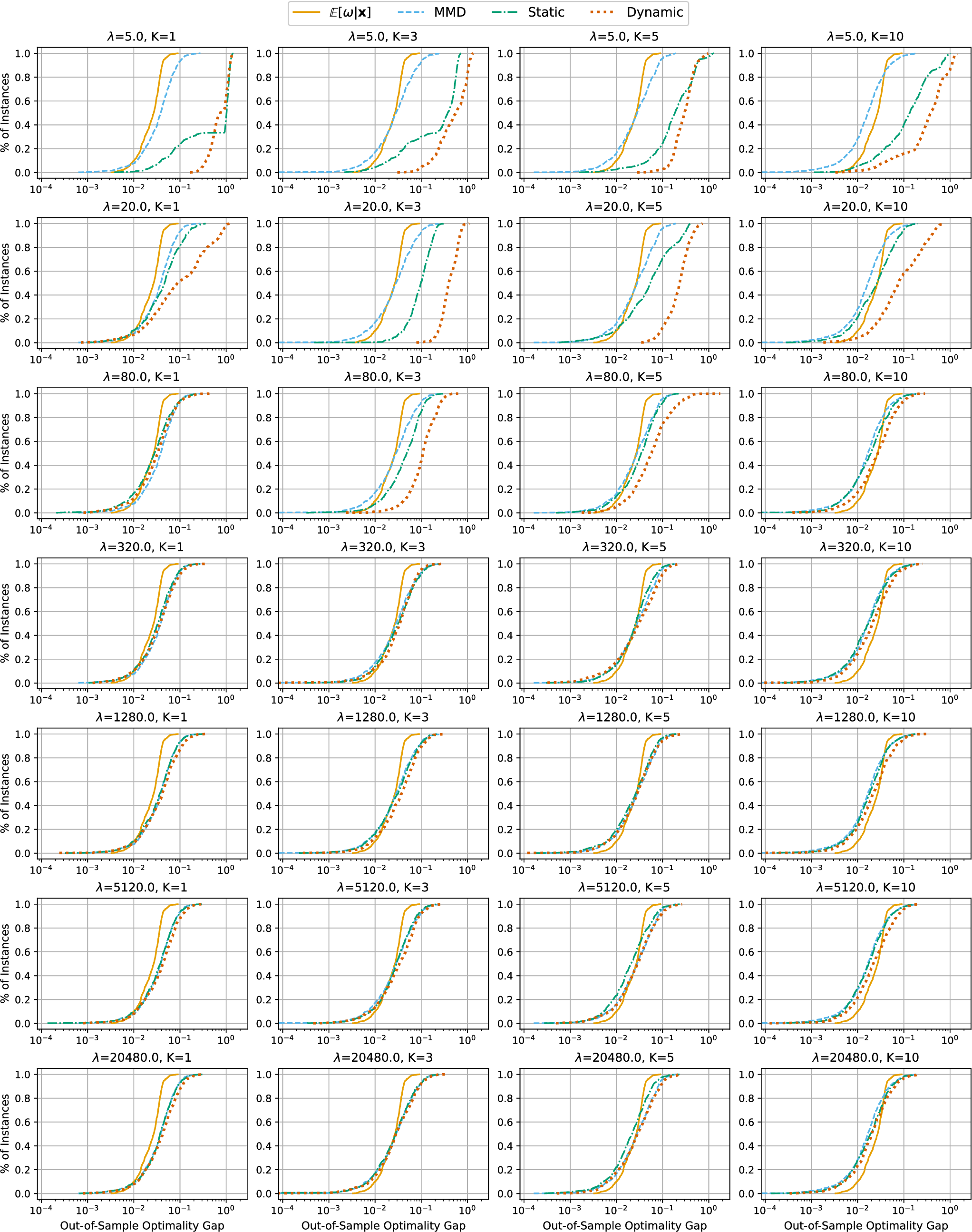}}
    {Out of sample optimality gaps for the contextual MNV problem \label{fig:MNV_gaps}}
    {}
\end{figure}

\subsection{Timing analysis }\label{appendix:timing}

Tables \ref{tab:CVaR_times} and \ref{tab:MNV_times} show the aforementioned metrics by $K$. Results for the newsvendor and CEP1 instances are not shown since directly solving the 2SPs is not as computationally challenging. We focus on the CVaR and MNV problems since the CVaR problem considers 1559 support points to define $\mathbb{P}_{\omega | \bs{x}}$ and MNV is complicated by mixed-binary variables in both stages. We do not show the results by $\lambda$ as it did not significantly impact timing. Furthermore, we do not show the time to compute solutions via the Task-Nets obtained from static and dynamic approaches since they did not vary relative to the MMD Task-Net (as their architecture is the same). 

\begin{table}[h]
\resizebox{\textwidth}{!}{%
    \TABLE
    {Computation times (in seconds) for contextual CVaR  \label{tab:CVaR_times} }
    {
    \begin{tabular}{cccccccc}
    \hline
    $K$ & \begin{tabular}[c]{@{}c@{}}MMD \\ Training\end{tabular} & \begin{tabular}[c]{@{}c@{}}Surrogate Solution \\ Calculation (MMD)\end{tabular} & \begin{tabular}[c]{@{}c@{}} 2SP Solve \\ $n_{\text{eval}} = 312$ \end{tabular} & \begin{tabular}[c]{@{}c@{}}MMD Loss\\ Evaluation\end{tabular} & \begin{tabular}[c]{@{}c@{}}Loss-Net \\ Training\end{tabular} & \begin{tabular}[c]{@{}c@{}}Static  \\  Training\end{tabular} & \begin{tabular}[c]{@{}c@{}}Dynamic \\ Training\end{tabular} \\ \hline
    1  & 167.6 & 1.2  & 1422.1 & 13.2  & 31.5 & 222.1 & 766.1  \\
3  & 172.4 & 3.0  & 1465.8 & 21.9  & 31.3 & 219.5 & 802.1  \\
5  & 166.7 & 3.0  & 1339.4 & 28.1  & 30.6 & 215.9 & 825.0  \\
10 & 168.8 & 5.6  & 1340.1 & 46.3  & 31.7 & 218.8 & 935.4  \\
20 & 167.1 & 11.1 & 1331.7 & 87.4  & 30.5 & 214.7 & 1123.6 \\
40 & 168.2 & 33.1 & 1392.0 & 273.8 & 31.5 & 217.7 & 1950.2 \\ \hline
    \end{tabular}%
    }{The median over $\lambda$ is displayed if the relevant method is dependent on $\lambda$ }
}
\end{table}

\begin{table}[h]
\resizebox{\textwidth}{!}{%
    \TABLE
    {Median computation times (in seconds) for contextual MNV
    \label{tab:MNV_times}}
    {
    \begin{tabular}{cccccccc}
    \hline
    $K$ & \begin{tabular}[c]{@{}c@{}}MMD \\ Training\end{tabular} & \begin{tabular}[c]{@{}c@{}}Surrogate Solution\\  Calculation (MMD)\end{tabular} & \begin{tabular}[c]{@{}c@{}} 2SP Solve \\ $n_{\text{eval}} = 300$ \end{tabular} & \begin{tabular}[c]{@{}c@{}}MMD Loss \\ Evaluation\end{tabular} & \begin{tabular}[c]{@{}c@{}}Loss-Net \\ Training\end{tabular} & \begin{tabular}[c]{@{}c@{}}Static  \\  Training\end{tabular} & \begin{tabular}[c]{@{}c@{}}Dynamic \\ Training\end{tabular} \\ \hline
    1   & 22.6                                                    & 0.7                                                                             & 233.0     & 11.3                                                           & 20.5                                                         & 29.5                                                                & 203.6                                                       \\
    3   & 23.6                                                    & 2.3                                                                             & 245.3     & 25.6                                                           & 21.3                                                         & 30.4                                                                & 256.6                                                       \\
    5   & 21.8                                                    & 4.3                                                                             & 236.1     & 40.1                                                           & 18.8                                                         & 27.5                                                                & 308.4                                                       \\
    10  & 20.9                                                    & 8.4                                                                             & 234.1     & 84.9                                                           & 19.1                                                         & 27.9                                                                & 491.6                                                       \\ \hline
    \end{tabular}%
    }{}
}
\end{table}

\section{Generalization Bounds}
\subsection{Statistical Learning Preliminaries}\label{section:appendix_learning_bounds}
We restate the classic result stemming from \citetSM{bartlett2002rademacherSM}. We aim to follow \citetSM{mohri2018foundationsSM}'s treatment of Rademacher complexities. For the reader's convenience, we recall the definition of Rademacher complexities below.

\begin{samepage}
\begin{definition}[Rademacher Averages]
  Let $\mathcal{G}$ be a family of functions mapping $\Xi \subseteq \mathbb{R}^{\hat{p}}$ to $\mathbb{R}$, for some $\hat{p} \in \mathbb{N}$ and $S = \{\bs{\xi}^{(i)}\}_{i=1}^n$ denote a fixed sample of size $n$, with $\bs{\xi}^{(i)} \in \Xi, \ i \in [n]$. Then, the empirical Rademacher complexity of $\mathcal{G}$ with respect to $S$ is defined as
  $$
  \hat{R}_{n}(\mathcal{G}, S) 
  \coloneq \mathbb{E}_{\bs{\sigma}} \left[ \sup_{g \in \mathcal{G}} \frac{1}{n}\sum_{i=1}^n \sigma_i g(\bs{\xi}^{(i)}) \right],
  $$
  where $\sigma_1, \hdots, \sigma_n$ denote independent uniform random variables taking values in $\{-1, +1\}$. 

  Furthermore, one can define the empirical vector Rademacher complexity by considering classes of vector-valued functions $\mathcal{F} \subseteq \{\bs{f}:\Xi \rightarrow \mathbb{R}^K\}$. The empirical Rademacher complexity of the class of vector-valued functions $\mathcal{F}$ with respect to $S$ is then defined as
  $$
  \hat{R}_{n}(\mathcal{F}, S) 
  \coloneq \mathbb{E}_{\bs{\sigma}} \left[ \sup_{\bs{f} \in \mathcal{F}} \frac{1}{n} \sum_{i=1}^n \sum_{j=1}^K \sigma_{ij} f_j(\bs{\xi}^{(i)}) 
  \right],$$
  where $\sigma_{ij}$ are also independent uniform variables taking values in $\{-1, +1\}$ \hbox{$\forall \ i \in [n]$}, $j\in [K]$, and $f_j(\bs{\xi}^{(i)})$ is the $j$th component of $\bs{f}(\bs{\xi}^{(i)})$. We refer to $\bs{\sigma}$ as the Rademacher variables.

  For both the scalar and vector-valued cases, the Rademacher complexity is given by
  $$
  R_{n}(\mathcal{G}) 
  \coloneq \mathbb{E}_{S} \left[\hat{R}_{n}(\mathcal{G}, S)\right].
  $$
\end{definition}
\end{samepage}

Next, we recall the following generalization bounds based on Rademacher averages, provided by \citetSM{mohri2018foundationsSM}.

\begin{theorem}[Rademacher Complexity Generalization Bounds (Theorem 3.3 \protect\citepSM{mohri2018foundationsSM})]
  \label{theorem:rademacher}
  Let $\mathcal{G}$ be a family of functions mapping from $\Xi$ to $\mathbb{R}$ that are bounded $|g(\bs{\xi})|\leq M \ \forall g \in \mathcal{G}, \ \bs{\xi} \in \Xi$. Then, for any $\delta > 0$, with probability at least $1-\delta$ over the draw of an iid sample $S = \{\bs{\xi}^{(i)}\}_{i=1}^n$ of size $n$, each of the following holds for all $g \in \mathcal{G}$:
  \begin{equation*}
  \mathbb{E}_{\bs{\xi}}[g(\bs{\xi})] \leq \frac{1}{n} \sum_{i=1}^{n} g(\bs{\xi}^{(i)}) + 2 R_n(\mathcal{G}) + M\sqrt{\frac{\log(1/\delta)}{2n}}
  \end{equation*}
  and
  \begin{equation*}
  \mathbb{E}_{\bs{\xi}}[g(\bs{\xi})] \leq \frac{1}{n} \sum_{i=1}^{n} g(\bs{\xi}^{(i)}) + 2 \hat{R}_n(\mathcal{G}, S) + 3 M \sqrt{\frac{\log(2/\delta)}{2n}}.
  \end{equation*}
  \end{theorem}

In our setting, we aim to upper-bound the expected losses associated with a particular choice of task mapping. In general, the loss of a given $(\hat{\bs{x}}, \hat{\bs{\omega}})$ is given by applying a particular loss function to the result of the task-mapping applied to $\hat{\bs{x}}$.  

% \begin{samepage}
% \begin{lemma}[Talagrand's Lemma (Lemma 5.7 \citepSM{mohri2018foundationsSM})]
%   \label{lemma:talagrand}
%   Let $S = \{\bs{\xi}^{(i)}\}_{i=1}^N$ denote a fixed sample of size $n$, with $\bs{\xi}^{(i)} \in \Xi, \ i \in [n]$, $\ell_1, \ldots, \ell_n$ be $L$-Lipschitz functions from $\mathbb{R}$ to $\mathbb{R}$ and $\sigma_1, \ldots, \sigma_n$ be Rademacher random variables. Then, for any hypothesis set $\mathcal{G}$ of real-valued functions, the following inequality holds:
%   \[
%    \mathbb{E}_{\bs{\sigma}} \left[ \sup_{g \in \mathcal{G}} \frac{1}{n} \sum_{i=1}^{n} \sigma_i (\ell_i \circ g)(\bs{\xi}^{(i)}) \right] \leq L \mathbb{E}_{\bs{\sigma}}\left[ \sup_{g \in \mathcal{G}}  \frac{1}{n} \sum_{i=1}^{n} \sigma_i g(\bs{\xi}^{(i)}) \right] = L \hat{R}_n(\mathcal{G}, S).
%   \]
%   In particular, if $\ell_i = \ell$ for all $i \in [n]$, then the following holds:
%   \[
%   \hat{R}_n(\ell \circ \mathcal{G}, S) \leq L \hat{R}_n(\mathcal{G}, S).
%   \]
%   \end{lemma}
% \end{samepage}
  \citetSM{mohri2018foundationsSM} present a version of Talagrand's Lemma which bounds the empirical Rademacher complexities of the function class obtained by composition of Lipschitz functions with a given function class. The following Lemma extends Talagrand's to the case where the $\ell_1,\hdots \ell_n$ are vector-valued. We consider the finite-dimensional version of \citetSM{maurer2016vectorSM}'s vector contraction inequality for Rademacher complexities. 
\begin{lemma}[Vector Contraction Inequality] 
\label{theorem:contraction}
Let $S = \{\bs{\xi}^{(i)}\}_{i=1}^n$ denote a fixed sample of size $n$, with $\bs{\xi}^{(i)} \in \Xi, \ i \in [n]$, $\ell_1, \ldots, \ell_n$ be $L$-Lipschitz functions from $\mathbb{R}^K$ to $\mathbb{R}$. Then, for any hypothesis set $\mathcal{F}$ being a class of functions $\bs{f}: \Xi \rightarrow \mathbb{R}^K$:
$$
\mathbb{E}_{\bs{\sigma}} \left[ \sup_{\bs{f} \in \mathcal{F}} \frac{1}{n} \sum_{i=1}^n \sigma_i \ell_i (\bs{f}(\bs{\xi}^{(i)})) \right] \leq \sqrt{2}L \hat{R}_{n}(\mathcal{F}, S),
$$
where $\bs{\sigma}$ are the Rademacher variables.
%, and $f_j(\bs{\xi}^{(i)})$ is the $j$th component of $\bs{f}(\bs{\xi}^{(i)})$.
\end{lemma}

Finally, we apply Lemma \ref{theorem:contraction} and Theorem \ref{theorem:rademacher} to construct a useful generalization bound for our setting. Also recall that $\Omega \subset \mathbb{R}^p$.

\begin{samepage}
\multivariaterademacher*
\end{samepage}
\begin{proof}
Define the family of functions $\mathcal{G} = \{(\bs{x}, \bs{\omega}) \mapsto \ell(\bs{f}(\bs{x}), \bs{\omega}): \bs{f} \in \mathcal{F} \}$. Lemma \ref{theorem:contraction} implies 
\begin{equation}\label{thm_proof:vector_inequality}
\hat{R}_{n}(\mathcal{G}, S)  = \mathbb{E}_{\bs{\sigma}} \left[ \sup_{\bs{f} \in \mathcal{F}} \frac{1}{n}\sum_{i=1}^n \sigma_i \ell_{\bs{\omega}^{(i)}}(\bs{f}(\bs{x}^{(i)})) \right] \leq \sqrt{2} L \hat{R}_{n}(\mathcal{F}, S),
\end{equation}
where
$$
\hat{R}_{n}(\mathcal{F}, S) 
= \mathbb{E}_{\bs{\sigma}} \left[ \sup_{\bs{f} \in \mathcal{F}} \frac{1}{n} \sum_{i=1}^n \sum_{j=1}^K \sum_{t = 1}^p \sigma_{ijt} f_{jt}(\bs{x}^{(i)}) 
\right]$$
denotes the empirical Rademacher average obtained by considering $\Omega^K$ as vectors in $\mathbb{R}^{pK}$. $f_{jt}(\bs{x}^{(i)})$ is the $t$th component of the $j$th vector in $\bs{f}(\bs{x}^{(i)}) = \left(\bs{f}_{1}(\bs{x}^{(i)}),\hdots,\bs{f}_{K}(\bs{x}^{(i)})\right)$. As before, $\bs{\sigma}$ are Rademacher variables. Applying Theorem \ref{theorem:rademacher} to $\mathcal{G}$ and inequality (\ref{thm_proof:vector_inequality}) yields the result. 
\end{proof}
Thus, if (i) $\hat{R}_{n}(\mathcal{F}, S)$ asymptotically vanishes as $n \rightarrow \infty$, (ii) $\ell$ is bounded, and (iii) $\ell_{\bs{\omega}}$ is Lipschitz $\forall \ \bs{\omega} \in \Omega$, then the task-mapping achieves generalization on unseen data.

\subsection{Rademacher Complexity of Bounded Linear Task Mapping}\label{appendix:taskmapping}

To help illustrate Theorem \ref{theorem:multivariate}, this subsection explores the hypothesis class for the task mapping consisting of $K$ linear transformations of the context. Based on this choice of hypothesis class, upper bounds on the Rademacher complexities are presented. The analysis proceeds similarly to Proposition 5 by \citetSM{oneto2020exploitingSM}. 

\begin{definition}[Bounded Linear Task Mapping]\label{defn:bounded}
Consider a bounded context space $\mathcal{X} = \{\bs{x} \in \mathbb{R}^d : \norm{\bs{x}}_2 \leq B \}$ along with the class of task mappings given by 
$$\mathcal{F} = \left\{\bs{f} : \mathcal{X} \rightarrow \Omega^K, \ \bs{f}_{j}(\bs{x}) = \rho(\bs{W}^{(j)} \bs{x}) : (\bs{W}^{(1)},\hdots \bs{W}^{(K)}) \in \mathcal{W} \right\},$$
with 
$$
\mathcal{W} = \left\{(\bs{W}^{(1)},\hdots \bs{W}^{(K)}) : \norm{\bs{W}^{(j)}}_{F} \leq 1, \ j = 1,\hdots ,K \right\},
$$
where the activation function $\rho : \mathbb{R} \rightarrow \mathbb{R}$ is $L_{\rho}$ Lipschitz and is applied element-wise. $\mathcal{F}$ is referred to as a bounded linear task mapping. 
\end{definition}
In the case of the bounded linear task mapping, we provide generalization bounds shown in the following theorem.  Although in practice, we do not parametrize $\mathcal{F}$ by $K$ bounded linear transformations, and instead elect for neural networks, the ability to generalize provides a solid footing for our proposed approach. 

\begin{restatable}{proposition}{rademacherlinear}\label{proposition:linear}
  Consider a bounded linear task mapping as defined in Definition \ref{defn:bounded}. Given a sample $S = \{\bs{x}^{(i)}\}_{i=1}^n$ of size $n$, it follows that
$$
\hat{R}_{n}(\mathcal{F}, S) \leq \frac{1}{n} L_{\rho} K  \sqrt{p \sum_{i=1}^n \norm[\big]{ \bs{x}^{(i)}}_2^2}
$$
and 
$$
R_{n}(\mathcal{F}) \leq L_{\rho} K B \sqrt{\frac{p}{n}}.
$$
\end{restatable}

The following proposition holds.

\begin{samepage}
\rademacherlinear*
\end{samepage}
\begin{proof}
Given a sample $S = \{\bs{x}^{(i)}\}_{i=1}^n$ of size $n$, it follows that 
\begin{equation}
    \hat{R}_{n}(\mathcal{F}, S) = \mathbb{E}_{\bs{\sigma}} \left[ \sup_{(\bs{W}^{(1)},\hdots \bs{W}^{(K)}) \in \mathcal{W}} \frac{1}{n} \sum_{i=1}^n \sum_{j=1}^K \sum_{t = 1}^p \sigma_{ijt} \rho\left( \langle \bs{w}_t^{(j)}, \bs{x}^{(i)} \rangle \right) \right],
  \end{equation}
  where $\bs{W} = (\bs{W}^{(1)},\hdots \bs{W}^{(K)})$ and $\bs{w}^{(j)}_t \in \mathbb{R}^d$ is the vector corresponding to the $t$th row of $\bs{W}^{(j)}$. The following inequalities hold:
    \begin{align}
      \hat{R}_{n}(\mathcal{F}, S) \leq \ &  L_{\rho}\mathbb{E}_{\bs{\sigma}} \left[ \sup_{(\bs{W}^{(j)})_{j=1}^K \in \mathcal{W}} \frac{1}{n} \sum_{i=1}^n \sum_{j=1}^K \sum_{t = 1}^p \sigma_{ijt} \langle \bs{w}_t^{(j)}, \bs{x}^{(i)} \rangle \right] & (\text{Lemma 5.7 \citepSM{mohri2018foundationsSM}})\nonumber \\
      = \ & L_{\rho} \sum_{j=1}^K \mathbb{E}_{\bs{\sigma}} \left[ \sup_{\norm{\bs{W}^{(j)}}_F \leq 1} \frac{1}{n} \sum_{i=1}^n \sum_{t = 1}^p \sigma_{ijt} \langle \bs{w}_t^{(j)}, \bs{x}^{(i)} \rangle  \right] & (\text{definition of} \ \mathcal{W}) \nonumber\\
      = \ & L_{\rho} K \mathbb{E}_{\bs{\sigma}} \left[ \sup_{\norm{\bs{W}}_F \leq 1} \frac{1}{n} \sum_{i=1}^n \sum_{t = 1}^p \sigma_{it} \langle \bs{w}_t, \bs{x}^{(i)} \rangle  \right] \nonumber \\
      = \ & \frac{1}{n} L_{\rho} K \mathbb{E}_{\bs{\sigma}} \left[ \sup_{\norm{\bs{W}}_F \leq 1} \sum_{t = 1}^p  \langle \bs{w}_t, \sum_{i=1}^n \sigma_{it} \bs{x}^{(i)} \rangle  \right]  & (\text{bilinearity of }\langle\cdot, \cdot \rangle).  \label{eq:linear_hypothesis_placeholder}
    \end{align}
Let $\bs{X}_{\bs{\sigma}} = \left(\sum_{i=1}^n \sigma_{i0} \bs{x}^{(i)}, \hdots, \sum_{i=1}^n \sigma_{ip} \bs{x}^{(i)} \right)^{\intercal} \in \mathbb{R}^{p \times d}$.
  \begin{align}
    \eqref{eq:linear_hypothesis_placeholder} = \ & \frac{1}{n} L_{\rho} K \mathbb{E}_{\bs{\sigma}} \left[ \sup_{\norm{\bs{W}}_F \leq 1}  \langle \bs{W}, \bs{X}_{\bs{\sigma}} \rangle_{F}  \right] \nonumber \\
    \leq \ & \frac{1}{n} L_{\rho} K \mathbb{E}_{\bs{\sigma}} \left[ \sup_{\norm{\bs{W}}_F \leq 1}  \norm{\bs{W}}_F  \norm{\bs{X}_{\bs{\sigma}}}_F   \right] & (\text{Cauchy inequality})\nonumber \\
    = \ & \frac{1}{n} L_{\rho} K \mathbb{E}_{\bs{\sigma}} \left[ \norm{\bs{X}_{\bs{\sigma}}}_F  \right]  & (\text{Evaluating the}\ \sup) \nonumber \\
    = \ & \frac{1}{n} L_{\rho} K \mathbb{E}_{\bs{\sigma}} \left[ \sqrt{\sum_{t = 1}^p \norm[\bigg]{\sum_{i=1}^n \sigma_{it} \bs{x}^{(i)}}_2^2} \right] \nonumber \\
    \leq \ & \frac{1}{n} L_{\rho} K  \sqrt{\sum_{t = 1}^p  \mathbb{E}_{\bs{\sigma}} \left[ \norm[\bigg]{\sum_{i=1}^n \sigma_{it} \bs{x}^{(i)}}_2^2 \right]}  & (\text{Jensen's inequality}) \nonumber \\
    = \ & \frac{1}{n} L_{\rho} K  \sqrt{\sum_{t = 1}^p \sum_{i=1}^n \norm[\big]{ \bs{x}^{(i)}}_2^2} & ((i',t') \neq (i,t) \implies  \mathbb{E}_{\bs{\sigma}}[\sigma_{it}\sigma_{i't'}] = 0)  \nonumber \\
    = \ & \frac{1}{n} L_{\rho} K  \sqrt{p \sum_{i=1}^n \norm[\big]{ \bs{x}^{(i)}}_2^2}   \nonumber \\
    \leq \ &  L_{\rho} K B \sqrt{\frac{p}{n}}.  & \left(\norm[\big]{ \bs{x}^{(i)}}_2 \leq B\right) \nonumber
  \end{align}
  Since $\hat{R}_{n}(\mathcal{F}, S) \leq L_{\rho} K B \sqrt{p/n}$, it follows that $\mathbb{E}_S \left[\hat{R}_{n}(\mathcal{F}, S)\right] \leq \mathbb{E}_S \left[L_{\rho} K B \sqrt{p/n}\right] = L_{\rho} K B \sqrt{p/n}$. Thus, the result holds. 
\end{proof}
Proposition \ref{proposition:linear} suggests that the number of surrogate scenarios $K$ and scenario dimension $p$ have an adverse impact on generalization ability. Applying Proposition \ref{proposition:linear} in combination with Theorem \ref{theorem:multivariate} implies the classic $O(1/\sqrt{n})$ generalization bounds, commonly observed by many machine learning algorithms. Although in practice, we do not select linear hypothesis classes for the task-mapping, it is reassuring that the CSG approach admits function classes that are capable of generalization. 

\subsection{MMD Loss Generalization Guarantees}\label{subsection:mmdloss}

This section considers Lipschitz continuity and boundedness of $\ell_{\text{MMD}}$. The following Lemma provides an upper bound on the Lipschitz constant of $\ell_{\text{MMD}, \bs{\omega}} : \bs{\zeta}_{1 \hdots K} \mapsto \ell_{\text{MMD}}(\bs{\zeta}_{1 \hdots K}, \bs{\omega})$. Let $L_k$ denote a bound on the Lipschitz constant of the kernel function. That is, the following holds for any fixed $\bs{\omega}'$:
$$
|k(\bs{\omega}, \bs{\omega}') - k(\eta, \bs{\omega}')| \leq L_k \norm{\bs{\omega} - \eta}_2 \ \forall\  \bs{\omega} \in \Omega, \ \eta \in \Omega.
$$
There are many examples of kernels that are Lipschitz. For example, the Gaussian kernel: $$k(\bs{\omega}, \bs{\omega}') = \sigma^2 \exp(-\norm{\bs{\omega} - \bs{\omega}'}_2^2/(2\alpha^2)),$$ with scalar parameters $\alpha$ and $\sigma$ has a Lipschitz constant $\sigma^2 e^{-0.5}/\alpha$. Since for any fixed $\bs{\omega}'$, $k(\bs{\omega}, \bs{\omega}')$ is the composition of  $\gamma: t \mapsto \sigma^2 \exp(-t^2/(2\alpha^2))$ with $d_{\bs{\omega}'}: \bs{\omega} \mapsto \norm{\bs{\omega} - \bs{\omega}'}_2$. By the reverse triangle inequality $d_{\bs{\omega}'}$ has Lipschitz constant $1$, and $\gamma$ has Lipschitz constant $\sigma^2 e^{-0.5}/\alpha$. Hence, by composition of Lipschitz functions,  $k_{\bs{\omega}'} : \bs{\omega} \mapsto k(\bs{\omega}, \bs{\omega}')$ has Lipschitz constant $\sigma^2 e^{-0.5}/\alpha$.

\kernellipschitz*
\begin{proof}
  \noindent \textbf{Lipschitz continuity of $\ell_{\text{MMD}, \bs{\omega}}(\cdot)$:} The proof proceeds similarly to Lemma 4 by \citetSM{oneto2020exploitingSM}. For any $\bs{\omega} \in \Omega$, 
  $$\ell_{\text{MMD}, \bs{\omega}}(\bs{\zeta}_{1 \hdots K}) = \frac{1}{K^2} \sum_{i, j}^K  k{\big (} \bs{\zeta}_i, \bs{\zeta}_j {\big )} - \frac{2}{K} \sum_{i=1}^K  k {\big (} \bs{\omega}, \bs{\zeta}_i {\big )}.$$
We consider one term at a time in computing a bound of the Lipschitz constant. Consider $\bs{\zeta}_{1 \hdots K}$ and $\tilde{\bs{\zeta}}_{1 \hdots K}$, it follows that 
  \begin{equation}\label{eqn:kernel_lipschitz}
  \frac{1}{K^2} \sum_{i, j} k{\big (} \bs{\zeta}_i, \bs{\zeta}_j {\big )} - \frac{1}{K^2} \sum_{i, j} k{\big (} \tilde{\bs{\zeta}}_i, \tilde{\bs{\zeta}}_j {\big )} = \frac{1}{K^2} \sum_{i, j} \left(k{\big (} \bs{\zeta}_i, \bs{\zeta}_j {\big )} -  k{\big (} \tilde{\bs{\zeta}}_i, \tilde{\bs{\zeta}}_j {\big )} \right).
  \end{equation}
Adding and subtracting $k{\big (}\bs{\zeta}_i, \tilde{\bs{\zeta}}_j{\big )}$ inside the sum implies 
\begin{align}
  k{\big (} \bs{\zeta}_i, \bs{\zeta}_j {\big )} -  k{\big (} \tilde{\bs{\zeta}}_i, \tilde{\bs{\zeta}}_j {\big )} = \ & k{\big (} \bs{\zeta}_i, \bs{\zeta}_j {\big )} - k{\big (}\bs{\zeta}_i, \tilde{\bs{\zeta}}_j{\big )} + k{\big (}\bs{\zeta}_i, \tilde{\bs{\zeta}}_j{\big )} -   k{\big (} \tilde{\bs{\zeta}}_i, \tilde{\bs{\zeta}}_j {\big )} \nonumber \\
  \leq \ & L_k \norm{\bs{\zeta}_j - \tilde{\bs{\zeta}}_j}_2 + L_k \norm{\bs{\zeta}_i - \tilde{\bs{\zeta}}_i}_2 \nonumber
\end{align}
Therefore, 
\begin{align}
\eqref{eqn:kernel_lipschitz} \leq  \ & \frac{1}{K^2} \sum_{i,j} L_k \norm{\bs{\zeta}_j - \tilde{\bs{\zeta}}_j}_2 + L_k \norm{\bs{\zeta}_i - \tilde{\bs{\zeta}}_i}_2 \nonumber \\
= \ & \frac{2L_k}{K} \sum_{i}  \norm{\bs{\zeta}_i - \tilde{\bs{\zeta}}_i}_2 \nonumber\\
\leq \ & \frac{2 L_k}{\sqrt{K}} \norm[\big]{\text{vec}(\bs{\zeta}_{1 \hdots K})- \text{vec}(\tilde{\bs{\zeta}}_{1 \hdots K})}_2,\nonumber
\end{align}
where the last line holds because $\bs{\alpha} \in \mathbb{R}_+^K \implies \sum_{i} \alpha_i = \norm{\bs{\alpha}}_1 \leq \sqrt{K} \norm{\bs{\alpha}}_2$. Setting $\alpha_i = \norm{\bs{\zeta}_i - \tilde{\bs{\zeta}}_i}_2$, implies $\norm{\bs{\alpha}}_2 = \norm[\big]{\text{vec}(\bs{\zeta}_{1 \hdots K})- \text{vec}(\tilde{\bs{\zeta}}_{1 \hdots K})}_2$. 
Similarly, 
\begin{align}
  \frac{2}{K} \sum_{i=1}^K  k {\big (} \bs{\omega}, \bs{\zeta}_i {\big )} - \frac{2}{K} \sum_{i=1}^K  k {\big (} \bs{\omega}, \tilde{\bs{\zeta}}_i {\big )} = \ & \frac{2}{K} \sum_{i=1}^K  k {\big (} \bs{\omega}, \bs{\zeta}_i {\big )} -  k {\big (} \bs{\omega}, \tilde{\bs{\zeta}}_i {\big )} \nonumber \\
  \leq  \ & \frac{2 L_k}{K} \sum_{i=1}^K  \norm{\bs{\zeta}_i - \tilde{\bs{\zeta}}_i}_2 \nonumber \\
  \leq \ & \frac{2 L_k}{\sqrt{K}} \norm[\big]{\text{vec}(\bs{\zeta}_{1 \hdots K})- \text{vec}(\tilde{\bs{\zeta}}_{1 \hdots K})}_2 \nonumber .
\end{align}
Therefore, $4L_k/\sqrt{K}$ is an upper bound for $\ell_{\text{MMD}, \bs{\omega}}$'s Lipschitz constant.

\noindent\textbf{Uniform Bound:}
Lastly, the boundedness of $k(\cdot, \cdot)$ implies $|\ell_{\text{MMD}}(\bs{f}(\bs{x}), \bs{\omega})| \leq 2U$, since there are two cases:
\begin{enumerate}[leftmargin=0pt, itemindent=*, align=left]
  \item $\ell_{\text{MMD}}\left(\bs{f}(\bs{x}), \bs{\omega}\right) \leq 0$, implying 
  $\frac{2}{K} \sum_{i = 1}^K k\left(\bs{\omega}, \bs{f}_i(\bs{x})\right) \geq \frac{1}{K^2} \sum_{i,j} k\left(\bs{f}_i(\bs{x}), \bs{f}_j(\bs{x}) \right).$
  \noindent Therefore, 
  \begin{align*}
    |\ell_{\text{MMD}}\left(\bs{f}(\bs{x}), \bs{\omega}\right)| &=  \frac{2}{K} \sum_{i = 1}^K k\left(\bs{\omega}, \bs{f}_i(\bs{x})\right) - \frac{1}{K^2} \sum_{i,j} k\left(\bs{f}_i(\bs{x}), \bs{f}_j(\bs{x}) \right) \\
    & \leq \frac{2}{K} \sum_{i = 1}^K k\left(\bs{\omega}, \bs{f}_i(\bs{x})\right) \\
    & \leq 2U.
  \end{align*}
  \item $\ell_{\text{MMD}}\left(\bs{f}(\bs{x}), \bs{\omega}\right) > 0$, implying 
  $\frac{2}{K} \sum_{i = 1}^K k\left(\bs{\omega}, \bs{f}_i(\bs{x})\right) < \frac{1}{K^2} \sum_{i,j} k\left(\bs{f}_i(\bs{x}), \bs{f}_j(\bs{x}) \right).$
  \noindent Therefore, 
  \begin{align*}
    |\ell_{\text{MMD}}\left(\bs{f}(\bs{x}), \bs{\omega}\right)| &=  \frac{1}{K^2} \sum_{i,j} k\left(\bs{f}_i(\bs{x}), \bs{f}_j(\bs{x}) \right) -\frac{2}{K} \sum_{i = 1}^K k\left(\bs{\omega}, \bs{f}_i(\bs{x})\right) \\
    & \leq \frac{1}{K^2} \sum_{i,j} k\left(\bs{f}_i(\bs{x}), \bs{f}_j(\bs{x}) \right)  \\
    & \leq U.
  \end{align*}
\end{enumerate}

\end{proof}

\subsection{Task Loss Generalization Guarantees}\label{subsection:optloss}
This section considers finite-sample guarantees for $\ell_{\text{opt}}$. Furthermore, we provide finite-sample guarantees that apply when $\ell_{\text{opt}}$ is approximated using a class of functions $\mathcal{G}$. Our proof structure mirrors that of \citetSM{grigas2021integratedSM} as they also use approximations (of a unique optimal solution map) to enable their problem-driven approach to predict probabilities over a fixed-support. We begin with the following corollary by \citetSM{still2018lecturesSM} regarding the Lipschitz continuity of the LP value function under right-hand side perturbations. 

\begin{corollary}[\protect\citepSM{still2018lecturesSM}]\label{cor:lp_rhs_lip}
Let \(\bs{E}\in\mathbb{R}^{\bar{m}\times \bar{n}}\) and \(\bs{c}\in\mathbb{R}^{\bar{n}}\) be fixed. For each
parameter \(\bs{\theta}\in\mathbb{R}^{\bar{m}}\), consider the LP
\begin{equation}\label{eq:P_t}
v(\bs{\theta})\;:=\;\max_{\bs{x}\in\mathbb{R}^{\bar{n}}}\ \bs{c}^{\intercal} \bs{x}
\quad \textit{s.t.}\quad
\bs{E}\bs{x}\le \bs{\theta}, \tag{$P(\bs{\theta})$}
\end{equation}
and define its feasible set and effective domain by
\[
F(\bs{\theta})\;:=\;\{\bs{x}\in\mathbb{R}^{\bar{n}}:\ \bs{E}\bs{x}\le \bs{\theta}\},\qquad
\text{dom}(F)\;:=\;\{\bs{\theta}\in\mathbb{R}^{\bar{m}}:\ F(\bs{\theta})\neq\emptyset\}.
\]
The dual of \(P(\bs{\theta})\) is
\begin{equation}\label{eq:D_t}
\min_{\bs{\pi}\in\mathbb{R}^{\bar{m}}}\ \bs{\theta}^{\intercal} \bs{\pi}
\quad \textit{s.t.}\quad
\bs{E}^{\intercal} \bs{\pi}=\bs{c},\ \ \bs{\pi}\ge \bs{0} . \tag{$D(\bs{\theta})$}
\end{equation}
Assume that \(D(\bs{\theta})\) is feasible. Then there exists a constant
\(L>0\) such that for all \(\bs{\theta}_1,\bs{\theta}_2 \in \text{dom}(F)\),
\[
\big|v(\bs{\theta}_1)-v(\bs{\theta}_2)\big|\ \le\ L\,\|\bs{\theta}_1-\bs{\theta}_2\|.
\]
\end{corollary}
Clearly, Corollary \ref{cor:lp_rhs_lip} applies to both minimization and maximization LPs by negating the objective function. Furthermore, the result extends to LPs with non-negativity constraints on a subset of the decision variables by applying the Corollary with 
\[
\bs{E}^{\prime} = \begin{bmatrix} \bs{E} \\ -\bs{I}_{\bar{n}} \end{bmatrix}, 
\] 
where $\bs{I}_{\bar{n}}$ is the $\bar{n} \times \bar{n}$ identity matrix and fixing the corresponding components of $\bs{\theta}$ to zero in the Lipschitz bound.

\textbf{Notation:} We define the following notation to simplify the presentation. Consider the second-stage recourse value function:
\[
\tilde{Q}(\bs{\eta})
:=\min_{\bs{z}\ge \bs{0}}\ \{\bs{q}^{\intercal} \bs{z}:\ \bs{D}\bs{z}=\bs{\eta}\},
\qquad \bs{\eta}\in\mathbb{R}^p.
\]
The dual of $\tilde{Q}(\bs{\eta})$ is
\begin{equation}\label{eq:recourse_dual}
\underline{Q}(\bs{\eta}) = \max_{\bs{\pi}\in\Pi}\ \bs{\pi}^{\intercal} \bs{\eta},
\qquad
\Pi := \{\bs{\pi}\in\mathbb{R}^p:\ \bs{D}^{\intercal} \bs{\pi}\le \bs{q}\}.
\end{equation} 
Note that by Assumption \ref{ass:complete_recourse}, $\tilde{Q}(\bs{\eta}) < \infty$ for all $\bs{\eta}\in\mathbb{R}^p$. Similarly, Assumption \ref{ass:dual_feasibility} implies $\Pi \neq \emptyset$, hence weak duality implies $\tilde{Q}(\bs{\eta}) > -\infty$. Therefore, $\tilde{Q}(\bs{\eta})$ is finite for all $\bs{\eta}\in\mathbb{R}^p$. Under the assumptions, both the primal and dual LPs are feasible implying strong duality: $\tilde{Q}(\bs{\eta}) = \underline{Q}(\bs{\eta})$.

Furthermore, the dual of the \eqref{eq:zeta_saa_rhs_only} is given by 
\begin{align}
\max_{\bs{\alpha} , \bs{\beta}, \bs{\pi}_1,\dots,\bs{\pi}_K}\quad
& \bs{b}^{\intercal} \bs{\alpha} - \bs{e}^{\intercal} \bs{\beta} + \sum_{i=1}^K \bs{\zeta}_i^{\intercal} \bs{\pi}_i \tag{\text{D-$\zeta$-SAA-RHS}} \label{eq:dual_zeta_saa_rhs_only}\\
\textit{s.t.}\quad
& \bs{A}^{\intercal} \bs{\alpha} - \bs{\beta} +  \bs{C}^{\intercal} \sum_{i=1}^K \bs{\pi}_i \le \bs{h}, \notag\\
& \bs{D}^{\intercal} \bs{\pi}_i \le \frac{1}{K}\bs{q},\quad \forall i\in\{1,\dots,K\}, \notag\\
& \bs{\alpha} \in \mathbb{R}^{p^{\prime}}_+, \quad \bs{\beta}\in\mathbb{R}^{s_1}_+, \quad	\bs{\pi}_i \in \mathbb{R}^p, \quad 	\forall i\in\{1,\dots,K\}, 	\notag
\end{align}
where $\bs{\alpha}$, $\bs{\beta}$ and $\bs{\pi}_i$ are the dual variables associated with the constraints in \eqref{eq:zeta_saa_rhs_only}. Likewise, the dual of \eqref{eq:opt_search_rhs_only} is given by
\begin{align}
\max_{\substack{\bs{\alpha},\,\bs{\beta},\,\varpi\\ \bs{\pi},\,\bs{\pi}_1,\dots,\bs{\pi}_K}}
 \quad
& \bs{b}^{\intercal} \bs{\alpha}
  - \bs{e}^{\intercal} \bs{\beta} + \bs{\omega}^{\intercal} \bs{\pi}
  + \sum_{i=1}^K \bs{\zeta}_i^{\intercal} \bs{\pi}_i
  - \varpi\, v^*(\bs{\zeta}_{1 \hdots K}) \tag{D-Opt-Search-RHS} \label{eq:dual_opt_search_rhs_only}\\
\textit{s.t.}\quad
& \bs{A}^{\intercal} \bs{\alpha}
  - \bs{\beta} + \bs{C}^{\intercal}\!\Big(\bs{\pi} + \sum_{i=1}^K \bs{\pi}_i\Big)
  \le (1+\varpi)\,\bs{h}, \notag\\
& \bs{D}^{\intercal} \bs{\pi} \le \bs{q}, \notag\\
& \bs{D}^{\intercal} \bs{\pi}_i \le \frac{\varpi}{K}\bs{q},\quad \forall i\in\{1,\dots,K\}, \notag\\
& \bs{\alpha} \in \mathbb{R}^{p^{\prime}}_+, \ \bs{\beta}\in\mathbb{R}^{s_1}_+, \ \varpi \geq 0, \ \bs{\pi} \in \mathbb{R}^p, \ \bs{\pi}_i \in \mathbb{R}^p,\  \forall i\in\{1,\dots,K\}, \notag
\end{align}
where $\bs{\alpha}$, $\bs{\beta}$, $\varpi$, $\bs{\pi}$ and $\bs{\pi}_1, \hdots, \bs{\pi}_K$ are the dual variables associated with the constraints in \eqref{eq:opt_search_rhs_only}.

\optlipschitz*
\begin{proof}
% Fix $\bs{\omega}\in\Omega$. We show that the map
% $\ell_{\text{opt},\bs{\omega}}:\bs{\zeta}_{1 \hdots K}\mapsto
% \ell_{\text{opt}, \bs{\omega}}(\bs{\zeta}_{1 \hdots K})$
% is Lipschitz. 

We begin by showing three useful facts that will be used to establish the Lipschitz continuity of $\ell_{\text{opt}, \bs{\omega}}(\cdot)$ and boundedness of $\ell_{\text{opt}}(\cdot)$.

\noindent \textit{Fact 1: Boundedness of $\Pi$ and Lipschitz continuity of $\tilde{Q}$. } 
$\Pi$ is bounded under complete recourse. Suppose for contradiction that $\Pi$ is unbounded. Since $\Pi$ is a polyhedron, unboundedness implies the existence of a nonzero recession direction
$\bs{d}\neq \bs{0}$ with $\bs{D}^{\intercal} \bs{d}\le \bs{0}$.
By complete recourse, there exists $\bs{z}_{\bs{d}}\ge \bs{0}$ such that $\bs{D}\bs{z}_{\bs{d}}=\bs{d}$.
Then
\[
\|\bs{d}\|_2^2
= \bs{d}^{\intercal} \bs{d}
= \bs{d}^{\intercal} \bs{D}\bs{z}_{\bs{d}}
= \bs{z}_{\bs{d}}^{\intercal} \bs{D}^{\intercal} \bs{d}
\le 0,
\]
which forces $\bs{d}=\bs{0}$, a contradiction. Hence $\Pi$ is bounded and
\[
U_{\Pi} \coloneq \max_{\bs{\pi}\in\Pi}\|\bs{\pi}\|_2 < \infty.
\]

We also note that $\tilde{Q}$ is Lipschitz continuous with constant $\max_{\bs{\pi}\in\Pi}\|\bs{\pi}\|_2$. For any $\bs{\eta}_1,\bs{\eta}_2\in\mathbb{R}^p$,
\begin{align*}
\tilde{Q}(\bs{\eta}_1)-\tilde{Q}(\bs{\eta}_2)
&=\max_{\bs{\pi}\in\Pi}\bs{\pi}^{\intercal} \bs{\eta}_1-\max_{\bs{\pi}\in\Pi}\bs{\pi}^{\intercal} \bs{\eta}_2 \\
&\le \max_{\bs{\pi}\in\Pi}\bs{\pi}^{\intercal}(\bs{\eta}_1-\bs{\eta}_2)
\le \max_{\bs{\pi}\in\Pi}\|\bs{\pi}\|_2\,\|\bs{\eta}_1-\bs{\eta}_2\|_2 = U_{\Pi} \|\bs{\eta}_1-\bs{\eta}_2\|_2
\end{align*}
where the last inequality follows from Cauchy--Schwarz. Swapping $\bs{\eta}_1$ and $\bs{\eta}_2$ yields the reverse bound, hence
\begin{equation} \label{eq:Q_lipschitz}
\big|\tilde{Q}(\bs{\eta}_1)-\tilde{Q}(\bs{\eta}_2)\big|
\le U_{\Pi} \|\bs{\eta}_1-\bs{\eta}_2\|_2,
\end{equation}
so $\tilde{Q}$ is Lipschitz continuous with constant $U_{\Pi}$.

\noindent \textit{Fact 2: Dual feasibility of \eqref{eq:zeta_saa_rhs_only} and \eqref{eq:opt_search_rhs_only}.} Under Assumption \ref{ass:dual_feasibility}, there exists $\bar{\bs{\pi}}$ in $\Pi$. It follows that one can construct a feasible solution to \eqref{eq:dual_zeta_saa_rhs_only} and \eqref{eq:dual_opt_search_rhs_only}. Indeed, let $\bar{\bs{\pi}} \in \Pi$. Then it follows that $$(\bs{\alpha}, \bs{\beta}, \bs{\pi}_1, \hdots, \bs{\pi}_K) = (\bs{0}, \bar{\bs{\beta}}, \bar{\bs{\pi}}/K, \hdots, \bar{\bs{\pi}}/K)$$ is a feasible solution to \eqref{eq:dual_zeta_saa_rhs_only}, where $\bar{\bs{\beta}} \coloneq (\bar{\beta}_j)_{j=1}^{s_1}$ with $\bar{\beta}_j = \max\{\bs{C}_j^{\intercal} \bar{\bs{\pi}} - h_j, 0 \}$ and $\bs{C}_j \in \mathbb{R}^{p}$ corresponding to column $j \in \{1,\dots, s_1\}$ of $\bs{C}$. Similarly,  $$(\bs{\alpha}, \bs{\beta}, \varpi, \bs{\pi}, \bs{\pi}_1, \hdots, \bs{\pi}_K) = (\bs{0}, 2\bar{\bs{\beta}}, 1, \bar{\bs{\pi}}, \bar{\bs{\pi}}/K, \hdots, \bar{\bs{\pi}}/K)$$ is a feasible solution to \eqref{eq:dual_opt_search_rhs_only}. 

\noindent \textit{Fact 3: Primal Solutions and Optimality of \eqref{eq:zeta_saa_rhs_only} and \eqref{eq:opt_search_rhs_only}}
Note by Fact 2, \eqref{eq:zeta_saa_rhs_only} is dual feasible. Furthermore, by Assumption \ref{ass:complete_recourse}, and $\mathcal{Y} \neq \emptyset$, \eqref{eq:zeta_saa_rhs_only} is primal feasible. Therefore, by strong duality for LPs, there exists an optimal solution $(\bs{y}, \bs{z}_1, \hdots, \bs{z}_K)$ to \eqref{eq:zeta_saa_rhs_only} such that $\bs{h}^{\intercal} \bs{y} + \frac{1}{K} \sum_{i=1}^K \bs{q}^{\intercal} \bs{z}_i \leq v^*(\bs{\zeta}_{1 \hdots K})$ holds with equality. Furthermore, for such $\bs{y}$, again by Assumption \ref{ass:complete_recourse}, there exists $\bs{z} \geq \bs{0}$ such that $\bs{C}\bs{y} + \bs{D} \bs{z} = \bs{\omega}$ for any $\bs{\omega} \in \Omega$. Thus, \eqref{eq:opt_search_rhs_only} has a primal feasible solution. Lastly by Fact 2. \eqref{eq:opt_search_rhs_only} has a dual feasible solution. Therefore, by strong duality for LPs there exists an optimal solution  $(\bs{y}, \bs{z}, \bs{z}_1, \hdots, \bs{z}_K)$ to \eqref{eq:opt_search_rhs_only}.

% Indeed let $(\bar{\bs{\alpha}}, \bar{\bs{\pi}}_1, \hdots, \bar{\bs{\pi}}_K)$ be a feasible solution to \eqref{eq:dual_zeta_saa_rhs_only}. Then setting $\varpi = 1$, $\bs{\alpha} = 2\bar{\bs{\alpha}}$, $\bs{\pi} = \sum_{i=1}^K \bar{\bs{\pi}}_i$, and $\bs{\pi}_i = \bar{\bs{\pi}}_i$ for all $i\in\{1,\dots,K\}$ gives a feasible solution to \eqref{eq:dual_opt_search_rhs_only}.

We now proceed to the main proof.

\noindent \textbf{Lipschitz continuity of $v^*:$}
Consider \eqref{eq:zeta_saa_rhs_only} and view the right-hand side coefficients as a single parameter vector, e.g.
\[
\bs{\theta} \;:=\;
\begin{bmatrix}
\bs{b}\\
\bs{e} \\ 
\text{vec}(\bs{\zeta}_{1 \hdots K})
\end{bmatrix}.
\]
In this case, $\text{dom}(F)$ corresponds to the values of $\bs{b}$, $\bs{e}$ and $\bs{\zeta}_{1 \hdots K}$ such that \eqref{eq:zeta_saa_rhs_only} is feasible. Then \eqref{eq:zeta_saa_rhs_only} is an LP whose parameter $\bs{\theta}$
appears only in the right-hand side of the feasibility system (in the sense of
Corollary \ref{cor:lp_rhs_lip}). Since \eqref{eq:zeta_saa_rhs_only} is dual feasible (Fact 2.), Corollary \ref{cor:lp_rhs_lip} yields the
existence of $L_v>0$ such that for any $\bs{\theta}_1,\bs{\theta}_2 \in \text{dom}(F)$,
\[
\big|v(\bs{\theta}_1)-v(\bs{\theta}_2)\big|\le L_v\|\bs{\theta}_1-\bs{\theta}_2\|.
\]
Now fix $\bs{b}$ and $\bs{e}$ in $\bs{\theta}_1$ and $\bs{\theta}_2$ such that $\bs{A}\bs{y} \geq \bs{b}$ and $\bs{0} \leq \bs{y} \leq \bs{e}$, i.e., $\mathcal{Y} \neq \emptyset$, and let $\bs{\theta}_1,\bs{\theta}_2$ differ only through the components corresponding to
$\text{vec}(\bs{\zeta}_{1 \hdots K})$. This gives a Lipschitz bound with respect
to $\bs{\zeta}_{1 \hdots K}$ i.e., there exists $L_v>0$ such that for any
$\bs{\zeta}_{1 \hdots K},\tilde{\bs{\zeta}}_{1\hdots K} \in \Omega^K$,
\begin{equation}\label{eq:v_lip_rephrased}
\big|v^*(\bs{\zeta}_{1 \hdots K})-v^*(\tilde{\bs{\zeta}}_{1\hdots K})\big|
\le L_v\,
\big\|\text{vec}(\bs{\zeta}_{1 \hdots K})-\text{vec}(\tilde{\bs{\zeta}}_{1\hdots K})\big\|_2.
\end{equation}
\eqref{eq:v_lip_rephrased} holds for any $\bs{\zeta}_{1 \hdots K}, \tilde{\bs{\zeta}}_{1\hdots K} \in \Omega^K$ since complete recourse (Assumption \ref{ass:complete_recourse}) guarantees a feasible solution to \eqref{eq:zeta_saa_rhs_only} for all $\bs{\zeta}_{1 \hdots K} \in \Omega^K$ and any $\bs{b}$ and $\bs{e}$ as selected above.

\noindent \textbf{Lipschitz continuity of }$\ell_{\text{opt}, \bs{\omega}}(\cdot)$\textbf{:}
For \eqref{eq:opt_search_rhs_only}, again treat
all right-hand side quantities as a parameter:
\[
\hat{\bs{\theta}}
\;:=\;
\begin{bmatrix}
\bs{b}\\
\bs{e}\\
\bs{\omega}\\
\text{vec}(\bs{\zeta}_{1 \hdots K})\\
v^*(\bs{\zeta}_{1 \hdots K})
\end{bmatrix}.
\]
In this case, $\text{dom}(F)$ corresponds to the values of $\bs{b}$, $\bs{e}$, $\bs{\omega}$, $\bs{\zeta}_{1 \hdots K}$, and $v^*(\bs{\zeta}_{1 \hdots K})$ such that \eqref{eq:opt_search_rhs_only} is feasible. Then \eqref{eq:opt_search_rhs_only} is an LP whose parameter $\hat{\bs{\theta}}$ appears only in the right-hand side of the feasibility system (in the sense of Corollary \ref{cor:lp_rhs_lip}). Indeed, by Fact 2, \eqref{eq:opt_search_rhs_only} is dual feasible.
Thus, Corollary \ref{cor:lp_rhs_lip} yields a constant $L'>0$ such that for any $\hat{\bs{\theta}}_1,\hat{\bs{\theta}}_2 \in \text{dom}(F)$,
\[
\big|v(\hat{\bs{\theta}}_1)-v(\hat{\bs{\theta}}_2)\big|
\le L'\|\hat{\bs{\theta}}_1-\hat{\bs{\theta}}_2\|.
\]
Now fix the components of $\hat{\bs{\theta}}_1$ and $\hat{\bs{\theta}}_2$ corresponding to $\bs{b}$ and $\bs{e}$ such that $\bs{A}\bs{y} \geq \bs{b}$, $ \bs{0} \leq \bs{y} \leq \bs{e}$, i.e. $\mathcal{Y} \neq \emptyset$. Let $\hat{\bs{\theta}}_1,\hat{\bs{\theta}}_2$ differ only through $\text{vec}(\bs{\zeta}_{1 \hdots K})$, the scalar $v^*(\bs{\zeta}_{1 \hdots K})$ and $\bs{\omega}$.
Thus, for any $\bs{\zeta}_{1 \hdots K},\tilde{\bs{\zeta}}_{1\hdots K} \in \Omega^K$, $\bs{\omega} \in \Omega$
\begin{align}
&\big|\ell_{\text{opt}}(\bs{\zeta}_{1 \hdots K},  \bs{\omega})
-\ell_{\text{opt}}(\tilde{\bs{\zeta}}_{1\hdots K}, \tilde{\bs{\omega}})\big|
\qquad \notag \\
& \leq L'\,
\sqrt{
\big\|\text{vec}(\bs{\zeta}_{1 \hdots K})-\text{vec}(\tilde{\bs{\zeta}}_{1\hdots K})\big\|_2^2
+\big|v^*(\bs{\zeta}_{1 \hdots K})-v^*(\tilde{\bs{\zeta}}_{1\hdots K})\big|^2 + \big\| \bs{\omega} - \tilde{\bs{\omega}} \big\|_2^2
}.
\label{eq:opt_lip_intermediate}
\end{align}
\eqref{eq:opt_lip_intermediate} holds for any $\bs{\zeta}_{1 \hdots K}, \, \tilde{\bs{\zeta}}_{1\hdots K} \in \Omega^K$ and $\bs{\omega}, \, \tilde{\bs{\omega}}$ since Fact 3. guarantees that \eqref{eq:opt_search_rhs_only} is feasible for all $\bs{\zeta}_{1 \hdots K} \in \Omega^K$ and $\bs{\omega} \in \Omega$ for any $\bs{b}$, and $\bs{e}$ as selected above. Finally, applying \eqref{eq:opt_lip_intermediate} yields
\[
\bigl|\ell_{\mathrm{opt}}(\bs{\zeta}_{1 \hdots K},\bs{\omega})
-\ell_{\mathrm{opt}}(\tilde{\bs{\zeta}}_{1\hdots K},\tilde{\bs{\omega}})\bigr|
\le
L'\sqrt{1+L_v^2}\,
\Bigl\|
\mathrm{vec}\!\bigl([\bs{\zeta}_{1 \hdots K},\bs{\omega}]\bigr)
-
\mathrm{vec}\!\bigl([\tilde{\bs{\zeta}}_{1\hdots K},\tilde{\bs{\omega}}]\bigr)
\Bigr\|_2 .
\]
Therefore an upper bound on the Lipschitz constant is
$L_{\text{opt}}=L'\sqrt{1+L_v^2}$. Furthermore fixing $\bs{\omega} = \tilde{\bs{\omega}}$ implies
\[
\big|\ell_{\text{opt}, \bs{\omega}}(\bs{\zeta}_{1 \hdots K})
-\ell_{\text{opt}, \bs{\omega}}(\tilde{\bs{\zeta}}_{1\hdots K})\big|
\le
L'\sqrt{1+L_v^2}\,
\big\|\text{vec}(\bs{\zeta}_{1 \hdots K})-\text{vec}(\tilde{\bs{\zeta}}_{1\hdots K})\big\|_2.
\]

% As discussed, complete recourse and dual feasibility of \eqref{eq:zeta_saa_rhs_only} guarantees \eqref{eq:opt_search_rhs_only} is feasible for all $\bs{\zeta}_{1 \hdots K}$ and $\bs{\omega}$ (and appropriate $\bs{b}$). Furthermore, since a dual solution to \eqref{eq:opt_search_rhs_only} can be constructed from a dual solution to \eqref{eq:zeta_saa_rhs_only}, we conclude that \eqref{eq:opt_search_rhs_only} is lower bounded. Hence
% $\ell_{\text{opt}, \bs{\omega}}(\bs{\zeta}_{1 \hdots K})$ is finite.

\noindent\textbf{Uniform bound for $v^*$:} 
By Fact 3., for any $\bs{\zeta}_{1 \hdots K} \in \Omega^K$ there exists $(\bar{\bs{y}},\bar{\bs{z}}_1,\dots,\bar{\bs{z}}_K)$
that is feasible for \eqref{eq:zeta_saa_rhs_only}. By optimality, and the definition of $\tilde{Q}(\cdot)$:
\[
v^*(\bs{\zeta}_{1 \hdots K})
\le \bs{h}^{\intercal} \bar{\bs{y}} + \frac{1}{K}\sum_{i=1}^K \tilde{Q}(\bs{\zeta}_i-\bs{C}\bar{\bs{y}})
= \bs{h}^{\intercal} \bar{\bs{y}} + \frac{1}{K}\sum_{i=1}^K \max_{\bs{\pi} \in \Pi} \bs{\pi}^{\intercal}(\bs{\zeta}_i-\bs{C}\bar{\bs{y}}).
\]
Since $\bs{\pi}^{\intercal} \bs{v} \leq \|\bs{\pi}\|_2 \|\bs{v}\|_2$ for any $\bs{\pi}, \bs{v} \in \mathbb{R}^p$, we have
\[
\tilde{Q}(\bs{\zeta}_i-\bs{C}\bar{\bs{y}}) \leq \ U_{\Pi} \|\bs{\zeta}_i-\bs{C}\bar{\bs{y}}\|_2 \leq U_{\Pi} \bigl(\|\bs{\zeta}_i\|_2 + \|\bs{C}\bar{\bs{y}}\|_2\bigr) \leq U_{\Pi} \bigl(U_{\Omega}+ \|\bs{C}\bar{\bs{y}}\|_2\bigr)  ,
\]
where we have used Fact 1. ($\Pi$ is bounded) and $\Omega$ is bounded. This implies 
\[v^*(\bs{\zeta}_{1 \hdots K})
\le \bs{h}^{\intercal} \bar{\bs{y}} + U_{\Pi} \bigl(U_{\Omega}+ \|\bs{C}\bar{\bs{y}}\|_2\bigr)
=: \overline{U}_v,
\]
Let $(\bs{\alpha}, \bs{\beta}, \bs{\pi}_1, \hdots, \bs{\pi}_K) = (\bs{0}, \bar{\bs{\beta}}, \bar{\bs{\pi}}/K, \hdots, \bar{\bs{\pi}}/K)$ denote the feasible solution to \eqref{eq:dual_zeta_saa_rhs_only} used in Fact 2. Then by weak duality,
\[v^*(\bs{\zeta}_{1 \hdots K})
\geq -\sum_{j=1}^{s_1} e_j \max\{\bs{C}_j^{\intercal} \bar{\bs{\pi}} - h_j, 0 \} + \frac{1}{K}\sum_{i=1}^K \bs{\zeta}_i^{\intercal} \bar{\bs{\pi}}.
\]
Furthermore, since $\bs{C}_j^{\intercal} \bar{\bs{\pi}} - h_j \leq |\bs{C}_j^{\intercal} \bar{\bs{\pi}} |  +  |h_j| \leq \| \bs{C}_j\|_2 \|\bar{\bs{\pi}}\|_2 + |h_j|$, and $\| \bs{C}_j\|_2 \|\bar{\bs{\pi}}\|_2 + |h_j| \geq 0$, it follows that 
\begin{equation*}
\max\{\bs{C}_j^{\intercal} \bar{\bs{\pi}} - h_j, 0 \} \leq \| \bs{C}_j\|_2 \|\bar{\bs{\pi}}\|_2 + |h_j| \leq  \| \bs{C}_j\|_2 U_{\Pi} + |h_j|, \quad \forall j \in \{1,\dots, s_1\},
\end{equation*}
where we have used Fact 1. Therefore, since $\bs{e} > \bs{0}$, we have
\begin{equation}\label{eq:dual_bound_intermediate}
-\sum_{j=1}^{s_1} e_j \max\{\bs{C}_j^{\intercal} \bar{\bs{\pi}} - h_j, 0 \} \geq - \sum_{j=1}^{s_1} e_j \bigl( \| \bs{C}_j\|_2 U_{\Pi} + |h_j| \bigr),
\end{equation}
In addition, since $\bs{\pi}^{\intercal} \bs{v} \geq - \|\bs{\pi}\|_2 \|\bs{v}\|_2$ for any $\bs{\pi}, \bs{v} \in \mathbb{R}^p$, we have $
\bs{\zeta}_i^{\intercal} \bar{\bs{\pi}} \geq - \|\bs{\zeta}_i\|_2 \|\bar{\bs{\pi}}\|_2 $.
Therefore,
\begin{equation}\label{eq:dual_bound_intermediate2}
 \frac{1}{K}\sum_{i=1}^K \bs{\zeta}_i^{\intercal} \bar{\bs{\pi}} \geq -U_{\Omega} U_{\Pi}
\end{equation}
Combining \eqref{eq:dual_bound_intermediate} and \eqref{eq:dual_bound_intermediate2} gives the following lower bound on $v^*(\bs{\zeta}_{1 \hdots K})$ 
$$
v^*(\bs{\zeta}_{1 \hdots K}) \geq - \sum_{j=1}^{s_1} e_j \bigl( \| \bs{C}_j\|_2 U_{\Pi} + |h_j| \bigr) - U_{\Omega} U_{\Pi} =: \underline{U}_v.
$$
Therefore, for all $\bs{\zeta}_{1 \hdots K} \in \Omega^K$, $v^*(\bs{\zeta}_{1 \hdots K}) \in [\underline{U}_v, \overline{U}_v]$. Next, we use this bound on $v^*(\cdot)$ to bound $\ell_{\text{opt}}(\cdot)$.

\noindent \textbf{Uniform bound for} $\ell_{\text{opt}}$\textbf{:}
First note that the constraint $\bs{h}^{\intercal} \bs{y} + \frac{1}{K} \sum_{i=1}^K \bs{q}^{\intercal} \bs{z}_i \leq v^*(\bs{\zeta}_{1 \hdots K})$ in \eqref{eq:opt_search_rhs_only} forces
$(\bs{y},\bs{z}_1,\dots,\bs{z}_K)$ to be optimal for \eqref{eq:zeta_saa_rhs_only}.
Consequently, $\ell_{\mathrm{opt}}(\bs{\zeta}_{1 \hdots K},\bs{\omega})$ can be written as
\[
\ell_{\text{opt}}(\bs{\zeta}_{1 \hdots K},\bs{\omega})
= \min_{\bs{y}\in\mathcal{Y}^*(\bs{\zeta}_{1 \hdots K})}
\bs{h}^{\intercal} \bs{y} + \tilde{Q}(\bs{\omega}-\bs{C}\bs{y}),
\]
where recall $\mathcal{Y}^*(\bs{\zeta}_{1 \hdots K})$ denotes the set of $\bs{y}$ forming a part of an optimal solution to \eqref{eq:zeta_saa_rhs_only}.
Take any $\bar{\bs{y}} \in\mathcal{Y}^*(\bs{\zeta}_{1 \hdots K})$. Then
\[
v^*(\bs{\zeta}_{1 \hdots K})
= \bs{h}^{\intercal} \bar{\bs{y}} + \frac{1}{K}\sum_{i=1}^K \tilde{Q}(\bs{\zeta}_i-\bs{C}\bar{\bs{y}}).
\]
Using \eqref{eq:Q_lipschitz} with $\bs{\eta}_1=\bs{\omega}-\bs{C}\bar{\bs{y}}$ and $\bs{\eta}_2=\bs{\zeta}_i-\bs{C}\bar{\bs{y}}$,
\[
\tilde{Q}(\bs{\omega}-\bs{C}\bar{\bs{y}})
\le \tilde{Q}(\bs{\zeta}_i-\bs{C}\bar{\bs{y}}) + U_{\Pi}\|\bs{\omega}-\bs{\zeta}_i\|_2.
\]
Since $\|\bs{\omega}\|_2\le U_\Omega$ and $\|\bs{\zeta}_i\|_2\le U_\Omega$,
\[
\|\bs{\omega}-\bs{\zeta}_i\|_2 \le \|\bs{\omega}\|_2 + \|\bs{\zeta}_i\|_2 \le 2U_\Omega.
\]
Averaging over $i=1,\dots,K$ gives
\[
\tilde{Q}(\bs{\omega}-\bs{C}\bar{\bs{y}})
\le \frac{1}{K}\sum_{i=1}^K \tilde{Q}(\bs{\zeta}_i-\bs{C}\bar{\bs{y}}) + 2U_{\Pi}U_{\Omega}.
\]
Adding $\bs{h}^{\intercal} \bar{\bs{y}}$ yields
\[
\bs{h}^{\intercal} \bar{\bs{y}} + \tilde{Q}(\bs{\omega}-\bs{C}\bar{\bs{y}})
\le v^*(\bs{\zeta}_{1 \hdots K}) + 2U_{\Pi}U_{\Omega}.
\]
Since this holds for every $\bar{\bs{y}} \in\mathcal{Y}^*(\bs{\zeta}_{1 \hdots K})$, it holds in particular for the minimizer
defining $\ell_{\mathrm{opt}}$, hence
\[
\ell_{\mathrm{opt}}(\bs{\zeta}_{1 \hdots K},\bs{\omega})
\le v^*(\bs{\zeta}_{1 \hdots K}) + 2U_{\Pi}U_{\Omega}
\le \overline{U}_v + 2U_{\Pi}U_{\Omega}
=: \overline{U}_\ell,
\]
where we used the uniform upper bound for $v^*$.

Likewise, for the lower bound, fix any $\bar{\bs{y}}\in\mathcal{Y}^*(\bs{\zeta}_{1 \hdots K})$. Using \eqref{eq:Q_lipschitz} with
$\bs{\eta}_1=\bs{\zeta}_i-\bs{C}\bar{\bs{y}}$ and $\bs{\eta}_2=\bs{\omega}-\bs{C}\bar{\bs{y}}$ yields
\[
\tilde{Q}(\bs{\omega}-\bs{C}\bar{\bs{y}})
\ge \tilde{Q}(\bs{\zeta}_i-\bs{C}\bar{\bs{y}}) - U_{\Pi}\|\bs{\omega}-\bs{\zeta}_i\|_2.
\]
Since $\|\bs{\omega}-\bs{\zeta}_i\|_2\le \|\bs{\omega}\|_2+\|\bs{\zeta}_i\|_2\le 2U_{\Omega}$, we obtain
\[
\tilde{Q}(\bs{\omega}-\bs{C}\bar{\bs{y}})
\ge \tilde{Q}(\bs{\zeta}_i-\bs{C}\bar{\bs{y}}) - 2U_{\Pi}U_{\Omega}.
\]
Averaging over $i=1,\dots,K$ gives
\[
\tilde{Q}(\bs{\omega}-\bs{C}\bar{\bs{y}})
\ge \frac{1}{K}\sum_{i=1}^K \tilde{Q}(\bs{\zeta}_i-\bs{C}\bar{\bs{y}}) - 2U_{\Pi}U_{\Omega}.
\]
Adding $\bs{h}^{\intercal} \bar{\bs{y}}$ yields
\[
\bs{h}^{\intercal} \bar{\bs{y}}+\tilde{Q}(\bs{\omega}-\bs{C}\bar{\bs{y}})
\ge v^*(\bs{\zeta}_{1 \hdots K}) - 2U_{\Pi}U_{\Omega}.
\]
Since this holds for every $\bar{\bs{y}}\in\mathcal{Y}^*(\bs{\zeta}_{1 \hdots K})$, it holds in particular for the minimizer
defining $\ell_{\mathrm{opt}}$, hence
\[
\ell_{\mathrm{opt}}(\bs{\zeta}_{1 \hdots K},\bs{\omega})
\ge v^*(\bs{\zeta}_{1 \hdots K}) - 2U_{\Pi}U_{\Omega}
\ge \underline{U}_v - 2U_{\Pi}U_{\Omega}
=: \underline{U}_\ell,
\]
where we used the uniform lower bound for $v^*$.
Therefore, for all $\bs{\zeta}_{1 \hdots K}\in\Omega^K$ and $\bs{\omega}\in\Omega$,
$\ell_{\mathrm{opt}}(\bs{\zeta}_{1 \hdots K},\bs{\omega})\in[\underline{U}_\ell,\overline{U}_\ell]$.

\end{proof}

Next, we provide the proof of Theorem \ref{thm:appxlossguarantee}. 
\begin{samepage}
  \appxlossguarantee*
\end{samepage}

\begin{proof}
Fix $\delta > 0$. Under Assumptions \ref{ass:complete_recourse} and \ref{ass:dual_feasibility}, and boundedness of $\Omega$,
Lemma \ref{lemma:opt_lipschitz} implies that $\ell_{\text{opt}}$ is uniformly bounded by $U_{\text{opt}}$ and, for any fixed $\bs{\omega}\in\Omega$,
$\ell_{\text{opt},\bs{\omega}}$ is $L_{\text{opt}}$-Lipschitz.

Applying Theorem \ref{theorem:multivariate} to $\ell=\ell_{\text{opt}}$ with confidence parameter $\delta/2$
yields that, with probability at least $1-\delta/2$ over the draw of $S$,
\begin{align}
\mathbb{E}_{(\bs{x},\bs{\omega})\sim \mathbb{P}_{\bs{x},\bs{\omega}}}
\big[\ell_{\text{opt}}(\bs{f}(\bs{x}),\bs{\omega})\big]
&\le
\frac{1}{n}\sum_{i=1}^n \ell_{\text{opt}}\!\big(\bs{f}(\bs{x}^{(i)}),\bs{\omega}^{(i)}\big) \notag \\
& +2\sqrt{2}\,L_{\text{opt}}\,R_n(\mathcal{F})
+U_{\text{opt}}\sqrt{\frac{\log(2/\delta)}{2n}}. \label{eq:step1_multivar}
\end{align}
By the triangle inequality,
\begin{align}
\frac{1}{n}\sum_{i=1}^n \ell_{\text{opt}}\!\big(\bs{f}(\bs{x}^{(i)}),\bs{\omega}^{(i)}\big)
&\le
\frac{1}{n}\sum_{i=1}^n \hat{E}\!\big(\bs{f}(\bs{x}^{(i)}),\bs{\omega}^{(i)}\big)
\notag\\
&+
\frac{1}{n}\sum_{i=1}^n
\Big|\ell_{\text{opt}}\!\big(\bs{f}(\bs{x}^{(i)}),\bs{\omega}^{(i)}\big)
      -\hat{E}\!\big(\bs{f}(\bs{x}^{(i)}),\bs{\omega}^{(i)}\big)\Big|.
\label{eq:step2_triangle}
\end{align}
By Lemma \ref{lemma:losserror}, for any $g\in\mathcal{G}_{\varPsi}$ we also have
\[
\Big|\ell_{\text{opt}}(\bs{\zeta}_{1 \hdots K},\bs{\omega})-g(\bs{\zeta}_{1 \hdots K},\bs{\omega})\Big|
\le U_{\text{opt}}+U_{\varPsi},
\qquad \forall\,(\bs{\zeta}_{1 \hdots K},\bs{\omega})\in\Omega^{K+1}.
\]
Moreover, since $\hat{E}\in\mathcal{G}_{\varPsi}$ is selected using the independent sample $S''$,
we may condition on $S''$ and treat $\hat{E}$ as fixed when analyzing the draw of $S$.
Applying Lemma \ref{lemma:losserror}
(i.e., Theorem 11.3 in \citetSM{mohri2018foundationsSM}) to the distribution
$\mathbb{P}_{\bs{f}(\bs{x}),\bs{\omega}}$ with sample size $n$ and confidence parameter $\delta/4$
yields that, with probability at least $1-\delta/4$ over the draw of $S$,
\begin{align}
\frac{1}{n}\sum_{i=1}^n
\Big|\ell_{\text{opt}}\!\big(\bs{f}(\bs{x}^{(i)}),\bs{\omega}^{(i)}\big)
      & -\hat{E}\!\big(\bs{f}(\bs{x}^{(i)}),\bs{\omega}^{(i)}\big)\Big|
 \notag \\
& \leq \mathbb{E}_{(\bs{\zeta}_{1 \hdots K},\bs{\omega})\sim \mathbb{P}_{\bs{f}(\bs{x}),\bs{\omega}}}
\Big[
\Big|\ell_{\text{opt}}(\bs{\zeta}_{1 \hdots K},\bs{\omega})
      -\hat{E}(\bs{\zeta}_{1 \hdots K},\bs{\omega})\Big|
\Big] \notag \\ 
& + 2R_n(\mathcal{G}_{\varPsi}) \notag 
+(U_{\text{opt}}+U_{\varPsi})\sqrt{\frac{\log(4/\delta)}{2n}}.
\label{eq:step3_emp_to_exp}
\end{align}
By Lemma \ref{lemma:opt_lipschitz} and the assumption that every $g\in\mathcal{G}_{\varPsi}$ is $L_{\varPsi}$-Lipschitz,
the map
\[
(\bs{\zeta}_{1 \hdots K},\bs{\omega})\ \mapsto\
\Big|\ell_{\text{opt}}(\bs{\zeta}_{1 \hdots K},\bs{\omega})-\hat{E}(\bs{\zeta}_{1 \hdots K},\bs{\omega})\Big|
\]
is $(L_{\text{opt}}+L_{\varPsi})$-Lipschitz (with respect to the same norm used to define $d_{\text{Lip}}$).
Hence, by the definition of the IPM over $1$-Lipschitz test functions,
\begin{align}
\mathbb{E}_{(\bs{\zeta}_{1 \hdots K},\bs{\omega})\sim \mathbb{P}_{\bs{f}(\bs{x}),\bs{\omega}}}
\Big[
\Big|\ell_{\text{opt}}(\bs{\zeta}_{1 \hdots K},\bs{\omega})
      & -\hat{E}(\bs{\zeta}_{1 \hdots K},\bs{\omega})\Big|
\Big]
  \notag \\ 
& \leq \mathbb{E}_{(\bs{\zeta}_{1 \hdots K},\bs{\omega})\sim \mathbb{P}_{\bs{\zeta}_{1 \hdots K},\bs{\omega}}}
\Big[
\Big|\ell_{\text{opt}}(\bs{\zeta}_{1 \hdots K},\bs{\omega})
      -\hat{E}(\bs{\zeta}_{1 \hdots K},\bs{\omega})\Big|
\Big]
\notag\\
&+
(L_{\text{opt}}+L_{\varPsi})\,
d_{\text{Lip}}\!\Big(
\mathbb{P}_{\bs{\zeta}_{1 \hdots K}, \bs{\omega}},
\mathbb{P}_{\bs{f}(\bs{x}), \bs{\omega}}
\Big).
\label{eq:step4_ipm}
\end{align}
Applying Lemma \ref{lemma:losserror} to $\mathbb{P}_{\bs{\zeta}_{1 \hdots K},\bs{\omega}}$
with confidence parameter $\delta/4$ yields that, with probability at least $1-\delta/4$ over the draw of $S''$,
\begin{align}
\mathbb{E}_{(\bs{\zeta}_{1 \hdots K},\bs{\omega})\sim \mathbb{P}_{\bs{\zeta}_{1 \hdots K},\bs{\omega}}}
\Big[
\Big|\ell_{\text{opt}}(\bs{\zeta}_{1 \hdots K},\bs{\omega})
      &-\hat{E}(\bs{\zeta}_{1 \hdots K},\bs{\omega})\Big|
\Big] \notag \\
&\le
\frac{1}{m}\sum_{i=1}^m
\Big|\ell_{\text{opt}}\!\big(\bs{\zeta}^{(i)}_{1 \hdots K},\bs{\omega}^{(i)}\big)
      -\hat{E}\!\big(\bs{\zeta}^{(i)}_{1 \hdots K},\bs{\omega}^{(i)}\big)\Big|
\notag\\
&\quad+
2R_m(\mathcal{G}_{\varPsi})
+(U_{\text{opt}}+U_{\varPsi})\sqrt{\frac{\log(4/\delta)}{2m}}.
\label{eq:step5_exp_to_emp}
\end{align}

Combining \eqref{eq:step1_multivar}--\eqref{eq:step5_exp_to_emp} and using a union bound over the three events
(with failure probabilities $\delta/2$, $\delta/4$, and $\delta/4$) shows that, with probability at least $1-\delta$, \eqref{eq:appx:bound} holds.
\end{proof}

% \subsection*{Task-Based Loss can be Bounded and Lipschitz}
% In this subsection we illustrate an example of a class of 2SPs that admit task-based losses $\ell_{\text{P}}$ that bounded and Lipschitz for all $\bs{\omega} \in \Omega$. To begin, we rely on the following Lemma due to \citetSM{dontchev2009implicit}

% \begin{lemma}[Tilted Minimization of Strongly Convex Functions]
%   Let $g: \mathbb{R}^m \to \mathbb{R}$ be continuously differentiable on an open set $O$, and let $\mathcal{U} \subset O$ be a nonempty, closed, convex set on which $g$ is strongly convex with constant $\mu > 0$. Then for each $\bs{e} \in \mathbb{R}^N$ the problem
%   \begin{equation}
%   \min_{\bs{u} \in \mathcal{U}} \ g(\bs{u}) + \langle \bs{e}, \bs{u} \rangle
%   \end{equation}
%   has a unique solution $\bs{u}^*(\bs{e})$, and the solution mapping $\bs{u}^*$ satisfies 
%     \[
%     \|\bs{u}^*(\bs{e}) - \bs{u}^*(\bs{e}')\|_2 \leq \mu^{-1} \|\bs{e} - \bs{e}'\|_2
%     \]
%   \end{lemma}
  % is a Lipschitz continuous function on $\mathbb{R}^n$ (globally) with Lipschitz constant $\mu^{-1}$.

%   \subsubsection*{Example 2: Linear Two-Stage Stochastic Programs with Quadratic Regularization}

%%%%%%%%%%%%%%%%%%%%%%%%%%%%%%%%%%%%%%%%%%%%%%%%%%%%%%%%%%%%%%%%%%%%%%%%%%%%%%
% Appendix: E joint-kernel view of contextual MMD losses
%%%%%%%%%%%%%%%%%%%%%%%%%%%%%%%%%%%%%%%%%%%%%%%%%%%%%%%%%%%%%%%%%%%%%%%%%%%%%%

\section{A joint-kernel interpretation of contextual MMD losses}\label{appx:joint_kernel_mmd}

This appendix complements the IPM/MMD background in the main text by showing how
the contextual MMD loss used in contextual scenario generation can be related to
an MMD between joint distributions on $\Omega^{K+1}$ via a suitably
constructed kernel $\kappa$. As discussed in the main text, the MMD between joint distributions is an IPM like the 1-Wasserstein distance appearing in Theorem~\ref{thm:appxlossguarantee} over a less rich class of functions. 

\subsection{Preliminaries: RKHS feature maps and MMD}\label{appx:prelim_mmd}
For convenience, we repeat the definition of the kernel $k$ and RKHS from Section~\ref{subsection:dsg}.
Let $k:\Omega\times \Omega\to\mathbb{R}$ be a positive semi-definite kernel with
associated RKHS $\mathcal{H}$.
Whenever $k$ is measurable and $\mathbb{E}_{\bs{\omega}\sim \mathbb{P}}[\sqrt{k(\bs{\omega},\bs{\omega})}]<\infty$,
the kernel mean embedding, introduced in Section \ref{subsection:dsg}, is given by $\mu_{\mathbb{P}}:=\mathbb{E}_{\bs{\omega}\sim\mathbb{P}}[k(\cdot,\bs{\omega})]\in\mathcal{H}$
and is well-defined. Furthermore, there exists an associated feature map
$\bs{\varphi}:\Omega\to\mathcal{H}$ such that
\[
k(\bs{\omega},\bs{\omega}')=\langle \bs{\varphi}(\bs{\omega}),\bs{\varphi}(\bs{\omega}')\rangle_{\mathcal{H}}
\quad \forall\,\bs{\omega},\bs{\omega}'\in\Omega.
\]
As stated in Section \ref{subsection:dsg}, for two distributions $\mathbb{P}_{\bs{\omega}}$ and $\mathbb{P}_{\bs{\xi}}$ on $\Omega$,
the squared MMD can be written as
\[
d^2_{\mathcal{G}_{\text{MMD}}}(\mathbb{P}_{\bs{\omega}},\mathbb{P}_{\bs{\xi}})
= \|\mu_{\mathbb{P}_{\bs{\omega}}}-\mu_{\mathbb{P}_{\bs{\xi}}}\|_{\mathcal{H}}^2,
\]
or an integral probability metric over functions with $\|\cdot\|_{\mathcal{H}} \leq 1$.
To distinguish between kernels, we write $d_k^2(\cdot,\cdot)$ for the MMD associated with kernel $k$.

\subsection{A kernel on $\Omega^{K+1}$ induced by $k$}\label{appx:kappa_def}
Fix $K\in\mathbb{N}$. For any $(\bs{\zeta}_{1 \hdots K},\bs{\omega})\in\Omega^{K+1}$, define the map
\begin{equation}\label{eq:psi_def}
\bs{\varsigma}(\bs{\zeta}_{1 \hdots K},\bs{\omega})
:= \frac{1}{K}\sum_{i=1}^K \bs{\varphi}(\bs{\zeta}_i)\;-\;\bs{\varphi}(\bs{\omega})
\;\in\;\mathcal{H}.
\end{equation}
Using $\bs{\varsigma}$ as a feature map on $\Omega^{K+1}$, define
\begin{equation}\label{eq:kappa_def}
\kappa\big((\bs{\zeta}_{1 \hdots K},\bs{\omega}),(\bs{\zeta}'_{1 \hdots K},\bs{\omega}')\big)
:= \big\langle \bs{\varsigma}(\bs{\zeta}_{1 \hdots K},\bs{\omega}),\,\bs{\varsigma}(\bs{\zeta}'_{1 \hdots K},\bs{\omega}')\big\rangle_{\mathcal{H}}.
\end{equation}

\begin{lemma}\label{lemma:kappa_psd}
If $k$ is positive semi-definite on $\Omega$, then $\kappa$ defined in \eqref{eq:kappa_def} is
positive semi-definite on $\Omega^{K+1}$. Moreover, using $k(\bs{\omega},\bs{\omega}')=\langle\bs{\varphi}(\bs{\omega}),\bs{\varphi}(\bs{\omega}')\rangle_{\mathcal{H}}$, $\kappa$ expands to
\begin{equation}\label{eq:kappa_expand}
\begin{aligned}
&\kappa\big((\bs{\zeta}_{1 \hdots K},\bs{\omega}),(\bs{\zeta}'_{1 \hdots K},\bs{\omega}')\big) \\
&\quad=
\frac{1}{K^2}\sum_{i=1}^K\sum_{j=1}^K k(\bs{\zeta}_i,\bs{\zeta}'_j)
\;+\; k(\bs{\omega},\bs{\omega}')
\;-\;\frac{1}{K}\sum_{i=1}^K k(\bs{\zeta}_i,\bs{\omega}')
\;-\;\frac{1}{K}\sum_{j=1}^K k(\bs{\omega},\bs{\zeta}'_j).
\end{aligned}
\end{equation}
Note that $\kappa$ is permutation-invariant in the first $K$ coordinates of each input since
it depends on $\bs{\zeta}_{1 \hdots K}$ only through sums over $i=1,\dots,K$.
\end{lemma}
\begin{proof}
Let $n\in\mathbb{N}$, $(\bs{\zeta}^{(a)}_{1 \hdots K},\bs{\omega}^{(a)})\in\Omega^{K+1}$ and $c_a\in\mathbb{R}$ for $a=1,\dots,n$.
Then
\[
\sum_{a=1}^n\sum_{b=1}^n c_a c_b\,
\kappa\big((\bs{\zeta}^{(a)}_{1 \hdots K},\bs{\omega}^{(a)}),(\bs{\zeta}^{(b)}_{1 \hdots K},\bs{\omega}^{(b)})\big)
=
\left\|\sum_{a=1}^n c_a\,\bs{\varsigma}(\bs{\zeta}^{(a)}_{1 \hdots K},\bs{\omega}^{(a)})\right\|_{\mathcal{H}}^2
\ge 0,
\]
so $\kappa$ is positive semi-definite. The expansion \eqref{eq:kappa_expand} follows from \eqref{eq:kappa_def} by expanding the inner product using $k(\bs{\omega},\bs{\omega}')=\langle\bs{\varphi}(\bs{\omega}),\bs{\varphi}(\bs{\omega}')\rangle_{\mathcal{H}}$.
\end{proof}

\subsection{Oracle and generated joint distributions}\label{appx:joint_distributions}
Recall $\mathbb{P}_{\bs{x}}$ and
$\mathbb{P}_{\bs{\omega}| \bs{x}}$ denote the distribution of $\bs{x}$ and the associated conditional distribution of $\bs{\omega}\in\Omega$
given $\bs{x}\in\mathcal{X}$. We define an \emph{oracle} joint distribution on $\Omega^{K+1}$ as follows.

\begin{assumption}[Oracle sampling model]\label{ass:oracle_sampling}
Conditioned on $\bs{x}\in\mathcal{X}$, the random variables
$\bs{\zeta}_1,\dots,\bs{\zeta}_K,\bs{\omega}$ are independent and identically distributed according to
$\mathbb{P}_{\bs{\omega}| \bs{x}}$.
\end{assumption}

\begin{definition}[Oracle joint distribution]\label{def:oracle_joint}
Under Assumption~\ref{ass:oracle_sampling}, let $\mathbb{P}^*_{\bs{\zeta}_{1 \hdots K},\bs{\omega}}$ denote
the distribution of $(\bs{\zeta}_{1 \hdots K},\bs{\omega})\in\Omega^{K+1}$ generated by the following sampling scheme:
\[
\bs{x}\sim\mathbb{P}_{\bs{x}},\quad
\bs{\zeta}_1,\ldots,\bs{\zeta}_K,\bs{\omega}\ \big|\ \bs{x}
\ \overset{\text{iid}}{\sim}\ \mathbb{P}_{\bs{\omega}\mid \bs{x}}.
\]
\end{definition}

We consider the oracle sampling model described above because it corresponds to the simplest baseline under the assumption of access to $\mathbb{P}_{\bs{\omega}\mid \bs{x}}$. Recall that $\bs{f}:\mathcal{X}\to\Omega^K$ is the scenario mapping, written
$\bs{f}(\bs{x})=(\bs{f}_1(\bs{x}),\dots,\bs{f}_K(\bs{x}))$, and the generated joint distribution
$\mathbb{P}_{\bs{f}(\bs{x}),\bs{\omega}}$ on $\Omega^{K+1}$ is defined by
$
\bs{x}\sim\mathbb{P}_{\bs{x}}, \ 
\bs{\zeta}_{1 \hdots K}=\bs{f}(\bs{x}), \ 
\bs{\omega}| \bs{x}\sim \mathbb{P}_{\bs{\omega}| \bs{x}}.$

\subsection{Joint MMD identity}\label{appx:identity}
Define the conditional mean embedding
$
\mu_{\bs{\omega}| \bs{x}} := \mathbb{E}_{\bs{\omega}\sim\mathbb{P}_{\bs{\omega}|\bs{x}}}\big[\bs{\varphi}(\bs{\omega})\big]\in\mathcal{H},$ 
and define the conditional embedding error of $\bs{f}$ at context $\bs{x}$:
\begin{equation*}
\Delta_{\bs{f}}(\bs{x})
:=\frac{1}{K}\sum_{i=1}^K \bs{\varphi}\big(\bs{f}_i(\bs{x})\big)\;-\;\mu_{\bs{\omega}| \bs{x}}
\;\in\;\mathcal{H}.
\end{equation*}
The following result shows that the MMD distance associated with the induced kernel $\kappa$ between
two joint distributions on $\Omega^{K+1}$ can be expressed in terms of the contextual embedding errors. Furthermore, in light of Theorem~\ref{thm:appxlossguarantee}, the MMD between joint distributions
$\mathbb{P}^*_{\bs{\zeta}_{1 \hdots K},\bs{\omega}}$ and $\mathbb{P}_{\bs{f}(\bs{x}),\bs{\omega}}$ provides a heuristic for controlling the
$1$-Wasserstein distance $d_{\text{Lip}}(\mathbb{P}_{\bs{\zeta}_{1 \hdots K},\bs{\omega}}, \mathbb{P}_{\bs{f}(\bs{x}),\bs{\omega}})$ that appears in the generalization guarantee. While the MMD function class is not as rich as the Wasserstein metric, the kernel-based IPM offers computational advantages
(see Appendix~\ref{appendix:comparing} for further discussion of these distances).

\begin{proposition}\label{prop:joint_mmd_identity}
Under Assumption~\ref{ass:oracle_sampling}, let $\mathbb{P}^*_{\bs{\zeta}_{1 \hdots K},\bs{\omega}}$ be as in
Definition~\ref{def:oracle_joint}, and let $\mathbb{P}_{\bs{f}(\bs{x}),\bs{\omega}}$ be the generated joint distribution.
Then, with $\kappa$ defined in \eqref{eq:kappa_def},
\begin{equation}\label{eq:joint_mmd_identity}
d^2_{\kappa}\!\Big(\mathbb{P}^*_{\bs{\zeta}_{1 \hdots K},\bs{\omega}},\,\mathbb{P}_{\bs{f}(\bs{x}),\bs{\omega}}\Big)  =\;\;
\big\|\mathbb{E}_{\bs{x}\sim\mathbb{P}_{\bs{x}}}\big[\Delta_{\bs{f}}(\bs{x})\big]\big\|_{\mathcal{H}}^2.
\end{equation}
\end{proposition}

\begin{proof}
For conditional expectations, we use the notation $(\bs{\zeta}_{1 \hdots K}, \bs{\omega})|\bs{x}$ to denote that $(\bs{\zeta}_{1 \hdots K}, \bs{\omega})$ is sampled from $\mathbb{P}_{\bs{\zeta}_{1 \hdots K}, \bs{\omega}|\bs{x}}$. Since $\kappa$ admits the feature map $\bs{\varsigma}$ in \eqref{eq:psi_def}, the squared MMD on $\Omega^{K+1}$ equals
\begin{align*}
d^2_{\kappa}\big(\mathbb{P}^*_{\bs{\zeta}_{1 \hdots K},\bs{\omega}},\,\mathbb{P}_{\bs{f}(\bs{x}),\bs{\omega}}\big)
& =
\big\|
\mathbb{E}_{(\bs{\zeta}_{1 \hdots K},\bs{\omega})\sim\mathbb{P}^*_{\bs{\zeta}_{1 \hdots K},\bs{\omega}}}\big[\bs{\varsigma}(\bs{\zeta}_{1 \hdots K},\bs{\omega})\big] \\
& -\mathbb{E}_{(\bs{f}(\bs{x}),\bs{\omega})\sim\mathbb{P}_{\bs{f}(\bs{x}),\bs{\omega}}}\big[\bs{\varsigma}(\bs{f}(\bs{x}),\bs{\omega})\big]
\big\|_{\mathcal{H}}^2.
\end{align*}
Under Assumption~\ref{ass:oracle_sampling}, for the oracle distribution, conditioning on $\bs{x}$ yields
\[
\mathbb{E}_{(\bs{\zeta}_{1 \hdots K},\bs{\omega})|\bs{x}}\big[\bs{\varsigma}(\bs{\zeta}_{1 \hdots K},\bs{\omega})\big]
=
\mathbb{E}_{(\bs{\zeta}_{1 \hdots K},\bs{\omega})|\bs{x}}\left[\frac{1}{K}\sum_{i=1}^K \bs{\varphi}(\bs{\zeta}_i)-\bs{\varphi}(\bs{\omega})\right]
=
\mu_{\bs{\omega}| \bs{x}}-\mu_{\bs{\omega}| \bs{x}}=0,
\]
hence $\mathbb{E}_{(\bs{\zeta}_{1 \hdots K},\bs{\omega})\sim\mathbb{P}^*_{\bs{\zeta}_{1 \hdots K},\bs{\omega}}}[\bs{\varsigma}(\bs{\zeta}_{1 \hdots K},\bs{\omega})]=0$.
For the generated distribution, conditioning on $\bs{x}$ gives
\[
\mathbb{E}_{(\bs{f}(\bs{x}),\bs{\omega})|\bs{x}}\big[\bs{\varsigma}(\bs{f}(\bs{x}),\bs{\omega})\big]
=
\frac{1}{K}\sum_{i=1}^K \bs{\varphi}\big(\bs{f}_i(\bs{x})\big) - \mu_{\bs{\omega}| \bs{x}}
=
\Delta_{\bs{f}}(\bs{x}),
\]
so $\mathbb{E}_{(\bs{f}(\bs{x}),\bs{\omega})\sim\mathbb{P}_{\bs{f}(\bs{x}),\bs{\omega}}}[\bs{\varsigma}(\bs{f}(\bs{x}),\bs{\omega})]=\mathbb{E}_{\bs{x}\sim\mathbb{P}_{\bs{x}}}[\Delta_{\bs{f}}(\bs{x})]$.
Substituting into the MMD expression proves \eqref{eq:joint_mmd_identity}.
\end{proof}

\subsection{Bias--variance decomposition}\label{appx:relationship_islip}
Recall from Section \ref{subsection:dsg} that the contextual scenario generation objective based on MMD is
\begin{equation*}
\mathcal{L}_{\text{dist}}(\bs{f})
:=
\mathbb{E}_{\bs{x}\sim\mathbb{P}_{\bs{x}}}\!\left[
d^2_{k}\!\Big(\widehat{\mu}_{\bs{f}(\bs{x})},\,\mathbb{P}_{\bs{\omega}| \bs{x}}\Big)
\right].
\end{equation*}
where the empirical measure associated with scenarios $\bs{f}(\bs{x})$ is denoted by $\widehat{\mu}_{\bs{f}(\bs{x})}$.
% \[
% \widehat{\mu}_{\bs{f}(\bs{x})}:=\frac{1}{K}\sum_{i=1}^K \delta_{\bs{f}_i(\bs{x})}.
% \]
Using the RKHS representation, for each fixed $\bs{x}$,
\[
d^2_{k}\!\Big(\widehat{\mu}_{\bs{f}(\bs{x})},\,\mathbb{P}_{\bs{\omega}| \bs{x}}\Big)
=
\left\|
\frac{1}{K}\sum_{i=1}^K \bs{\varphi}\big(\bs{f}_i(\bs{x})\big) - \mu_{\bs{\omega}| \bs{x}}
\right\|_{\mathcal{H}}^2
=
\|\Delta_{\bs{f}}(\bs{x})\|_{\mathcal{H}}^2.
\]
Hence
\begin{equation}\label{eq:ctx_loss_as_Delta}
\mathcal{L}_{\text{dist}}(\bs{f})=\mathbb{E}_{\bs{x}\sim\mathbb{P}_{\bs{x}}}\big[\|\Delta_{\bs{f}}(\bs{x})\|_{\mathcal{H}}^2\big].
\end{equation}

Combining \eqref{eq:joint_mmd_identity} and \eqref{eq:ctx_loss_as_Delta} yields the following decomposition.

\begin{corollary}[Bias--variance decomposition]\label{cor:bias_variance}
Under Assumption~\ref{ass:oracle_sampling},
\begin{equation}\label{eq:bias_variance}
\mathcal{L}_{\text{dist}}(\bs{f})
=
\underbrace{
d^2_{\kappa}\big(\mathbb{P}^*_{\bs{\zeta}_{1 \hdots K},\bs{\omega}},\,\mathbb{P}_{\bs{f}(\bs{x}),\bs{\omega}}\big)
}_{\text{joint MMD with kernel }\kappa}
+
\underbrace{
\mathbb{E}_{\bs{x}\sim\mathbb{P}_{\bs{x}}}\!\left[
\left\|\Delta_{\bs{f}}(\bs{x})-\mathbb{E}_{\bs{x}'\sim\mathbb{P}_{\bs{x}}}[\Delta_{\bs{f}}(\bs{x}')]\right\|_{\mathcal{H}}^2
\right]
}_{\text{variability term}}.
\end{equation}
In particular,
\begin{equation}\label{eq:joint_le_ctx}
d^2_{\kappa}\big(\mathbb{P}^*_{\bs{\zeta}_{1 \hdots K},\bs{\omega}},\,\mathbb{P}_{\bs{f}(\bs{x}),\bs{\omega}}\big)
\;\le\;
\mathcal{L}_{\text{dist}}(\bs{f}).
\end{equation}
\end{corollary}

\begin{proof}
Let $m:=\mathbb{E}_{\bs{x}\sim\mathbb{P}_{\bs{x}}}[\Delta_{\bs{f}}(\bs{x})]\in\mathcal{H}$. Then, adding and subtracting $m$ implies
\[
\|\Delta_{\bs{f}}(\bs{x})\|_{\mathcal{H}}^2
=
\|m\|_{\mathcal{H}}^2
+2\langle m,\Delta_{\bs{f}}(\bs{x})-m\rangle_{\mathcal{H}}
+\|\Delta_{\bs{f}}(\bs{x})-m\|_{\mathcal{H}}^2.
\]
Taking expectation over $\bs{x}\sim\mathbb{P}_{\bs{x}}$ makes the cross term vanish because
\hbox{$\mathbb{E}_{\bs{x}\sim\mathbb{P}_{\bs{x}}}[\Delta_{\bs{f}}(\bs{x})-m]=0$}, which yields \eqref{eq:bias_variance} after using
Proposition~\ref{prop:joint_mmd_identity} to identify $\|m\|_{\mathcal{H}}^2$ with the joint MMD term.
Inequality \eqref{eq:joint_le_ctx} follows since the variability term is nonnegative.
\end{proof}

% \end{appendices}

\bibliographystyleSM{sn-mathphys-num}
\bibliographySM{sn-bibliography-SM}

\end{document}